\documentclass[11pt]{amsart}
\usepackage{amscd,verbatim}
\usepackage{amssymb}
\usepackage[all]{xy}
\usepackage[inline]{enumitem}
\usepackage{hyphenat}
\usepackage[]{appendix}
\usepackage[colorlinks,linkcolor=blue,citecolor=blue,urlcolor=red]{hyperref}
\usepackage{cleveref}
\usepackage{xcolor}
\usepackage{mdframed}
\usepackage{bm}
\usepackage{amsmath}
\usepackage[normalem]{ulem}
\usepackage{mathtools}
\usepackage{tikz}
\usepackage{stmaryrd }
\allowdisplaybreaks
\usepackage[a4paper,bindingoffset=0in,%
            left=1in,right=1in,top=1in,bottom=1in,%
            footskip=.25in]{geometry}
\usepackage[french,english]{babel}
\usepackage[T1]{fontenc}
\makeatletter
\newenvironment{abstracts}{%
  \ifx\maketitle\relax
    \ClassWarning{\@classname}{Abstract should precede
      \protect\maketitle\space in AMS document classes; reported}%
  \fi
  \global\setbox\abstractbox=\vtop \bgroup
    \normalfont\Small
    \list{}{\labelwidth\z@
      \leftmargin3pc \rightmargin\leftmargin
      \listparindent\normalparindent \itemindent\z@
      \parsep\z@ \@plus\p@
      
      \itemsep\medskipamount
    }%
}{%
  \endlist\egroup
  \ifx\@setabstract\relax \@setabstracta \fi
}

\newcommand{\abstractin}[1]{%
  \otherlanguage{#1}%
  \item[\hskip\labelsep\scshape\abstractname.]%
}

\makeatletter
\newcommand\subsubsubsection{\@startsection{paragraph}{4}{\z@}{-1ex\@plus 0ex \@minus 2ex}{-1ex \@plus 2ex\@minus 2ex}{\normalfont\normalsize\itshape}}
\newcommand\subsubsubsubsection{\@startsection{subparagraph}{5}{\z@}{-2.5ex\@plus -1ex \@minus -.25ex}{1.25ex \@plus .25ex}{\normalfont\normalsize\itseries}}
\makeatother

\makeatletter
\expandafter\def\csname r@tocindent4\endcsname{0pt}
\makeatother

\let\oldtocsection=\tocsection
\let\oldtocsubsection=\tocsubsection
\let\oldtocsubsubsection=\tocsubsubsection
\let\oldtocparagraph=\tocparagraph
\let\oldtocsubparagraph=\tocsubparagraph
\renewcommand{\tocsection}[2]{\hspace{0em}\oldtocsection{#1}{#2}}
\renewcommand{\tocsubsection}[2]{\hspace{1em}\oldtocsubsection{#1}{#2}}
\renewcommand{\tocsubsubsection}[2]{\hspace{2em}\oldtocsubsubsection{#1}{#2}}
\renewcommand{\tocparagraph}[2]{\hspace{3em}\oldtocparagraph{#1}{#2}}
\renewcommand{\tocsubparagraph}[2]{\hspace{4em}\oldtocsubparagraph{#1}{#2}}

\newtheorem{lemma}{Lemma}[section]

\newtheorem{theorem}[lemma]{Theorem}
\newtheorem{cor-intro}{Corollary}

\newtheorem{proposition}[lemma]{Proposition}

\newtheorem{corollary}[lemma]{Corollary}

\theoremstyle{definition}

\newtheorem{definition}[lemma]{Definition}

\newtheorem{example}[lemma]{Example}
\newtheorem{remark}[lemma]{Remark}
\newtheorem{qn}[lemma]{Question}

\theoremstyle{remark}

\newtheorem*{claim*}{Claim}

\newcounter{zaehler}
\numberwithin{equation}{subsection}

\newcommand{\N}{\mathbb{N}}

\newcommand{\Z}{\mathbb{Z}}

\renewcommand{\P}{\mathbf{P}}
\newcommand{\A}{\mathbf{A}}
\newcommand{\F}{\mathbf{F}}

\newcommand{\cA}{\mathcal{A}}
\newcommand{\cB}{\mathcal{B}}

\newcommand{\cF}{\mathcal{F}}

\newcommand{\cH}{\mathcal{H}}
\newcommand{\cI}{\mathcal{I}}

\newcommand{\cL}{\mathcal{L}}
\newcommand{\cO}{\mathcal{O}}

\newcommand{\cK}{\mathcal{K}}
\newcommand{\cM}{\mathcal{M}}

\newcommand{\ssU}{{\mathfrak U}}

\newcommand{\Cone}{\operatorname{Cone}}

\newcommand{\Hom}{\operatorname{Hom}}

\newcommand{\Ker}{\operatorname{Ker}}

\newcommand{\Tr}{\operatorname{Tr}}
\newcommand{\Nm}{\operatorname{Nm}}

\newcommand{\Spec}{\operatorname{Spec}}

\newcommand{\Gal}{\operatorname{Gal}}

\newcommand{\tr}{{\operatorname{tr}}}

\newcommand{\Zar}{{\operatorname{Zar}}}

\newcommand{\et}{{\operatorname{\acute{e}t}}}

\newcommand{\Res}{\operatorname{Res}}

\newcommand{\id}{{\operatorname{id}}}

\newcommand{\Supp}{{\operatorname{Supp}\,}}

\newcommand{\CH}{{\operatorname{CH}}}

\newcommand{\oIm}{{\operatorname{Im}}}

\newcommand{\leftrarrows}{\mathrel{\raise.75ex\hbox{\oalign{%
  $\scriptstyle\leftarrow$\cr
  \vrule width0pt height.5ex$\hfil\scriptstyle\relbar$\cr}}}}
\newcommand{\lrightarrows}{\mathrel{\raise.75ex\hbox{\oalign{%
  $\scriptstyle\relbar$\hfil\cr
  $\scriptstyle\vrule width0pt height.5ex\smash\rightarrow$\cr}}}}
\newcommand{\Rrelbar}{\mathrel{\raise.75ex\hbox{\oalign{%
  $\scriptstyle\relbar$\cr
  \vrule width0pt height.5ex$\scriptstyle\relbar$}}}}

\makeatletter
\def\leftrightarrowsfill@{\arrowfill@\leftrarrows\Rrelbar\lrightarrows}
\newcommand{\xleftrightarrows}[2][]{\ext@arrow 3399\leftrightarrowsfill@{#1}{#2}}
\makeatother

\newcommand{\ev}{{\mathrm{ev}}}
\newcommand{\FRP}{{\mathrm{FRP}}}
\newcommand{\adj}{{\mathrm{adj}}}

\newcommand{\sep}{{\mathrm{sep}}}
\newcommand{\cHom}{{\mathcal Hom}}
\newcommand{\cExt}{{\mathcal Ext}}

\renewcommand{\tr}{\mathrm{tr}}

\newcommand{\trdeg}{\mathrm{trdeg}}

\newcommand{\qc}{{\mathrm{qc}}}

\newcommand{\ra}{\rightarrow}

\newcommand{\colim}{\operatorname{colim}}

\renewcommand{\tilde}{\widetilde}
\newcommand{\xra}{\xrightarrow}

\newcommand{\hra}{\hookrightarrow}

\newcommand{\xhra}{\xhookrightarrow}
\newcommand{\xhla}{\xhookleftarrow}

\renewcommand{\bar}{\overline}

\newcommand{\bthm}{\begin{theorem}}
\newcommand{\ethm}{\end{theorem}}
\newcommand{\bcor}{\begin{corollary}}
\newcommand{\ecor}{\end{corollary}}
\newcommand{\bprop}{\begin{proposition}}
\newcommand{\eprop}{\end{proposition}}
\newcommand{\ble}{\begin{lemma}}
\newcommand{\ele}{\end{lemma}}
\newcommand{\bex}{\begin{example}}
\newcommand{\eex}{\end{example}}
\newcommand{\bexc}{\begin{exercise}}
\newcommand{\eexc}{\end{exercise}}
\newcommand{\brmk}{\begin{remark}}
\newcommand{\ermk}{\end{remark}}
\newcommand{\bdefn}{\begin{definition}}
\newcommand{\edefn}{\end{definition}}
\newcommand{\bpf}{\begin{proof}}
\newcommand{\epf}{\end{proof}}
\newcommand{\benu}{\begin{enumerate}}
\newcommand{\eenu}{\end{enumerate}}
\newcommand{\bit}{\begin{itemize}}
\newcommand{\eit}{\end{itemize}}
\newcommand{\bqn}{\begin{qn}}
\newcommand{\eqn}{\end{qn}}
\newcommand{\beq}{\begin{equation}}
\newcommand{\eeq}{\end{equation}}

\begin{document}
\title[Bloch's cycle complex and coherent dualizing complexes]{Bloch's cycle complex and coherent dualizing complexes in positive characteristic
}

\author{Fei Ren}

\begin{abstracts}

\abstractin{english}
Let $X$ be a separated scheme of dimension $d$  of finite type over  a perfect field $k$ of positive characteristic $p$.
In this work, we show that Bloch's cycle complex $\Z^c_X$ of zero cycles mod $p^n$ is quasi-isomorphic to  the Cartier  operator fixed part of a certain dualizing complex from coherent duality theory. From this we obtain new vanishing results for the higher Chow groups of zero cycles with mod $p^n$ coefficients for singular varieties.
\end{abstracts}

\maketitle


\setcounter{tocdepth}{3}
\tableofcontents

\section*{Introduction}
Let $X$ be a separated scheme of dimension $d$  of finite type over  a perfect field $k$ of positive characteristic $p$.
In this work, we show that Bloch's cycle complex $\Z^c_X$ of zero cycles mod $p^n$ is quasi-isomorphic to  the Cartier  operator fixed part of a certain dualizing complex from coherent duality theory. From this we obtain new vanishing results for the higher Chow groups of zero cycles with mod $p^n$ coefficients for singular varieties.

As the first candidate for a motivic complex,
Bloch introduced his cycle complex $\Z^c_X$ in \cite{Bloch-CycleAndK} under the framework of Beilinson-Lichtenbaum. Let $i$ be an integer, and $\Delta^i=\Spec k[T_0,\dots, T_i]/\break(\sum T_j-1)$. 
Here $\Z^c_X:=z_0(-,-\bullet)$ is
a complex of sheaves in the Zariski or the \'etale topology. The global sections of its degree $(-i)$-term  $z_0(X,i)$ is the free abelian group generated by dimension $i$-cycles in $X\times \Delta^i$ intersecting all faces properly and the differentials are the alternating sums of the cycle-theoretic intersection of the cycle with each face 
(cf. \Cref{section bloch's cycle complex}).
Let $\pi:X\ra \Spec k$ be the structure morphism of $X$. Let $W_nX:=(|X|,W_n\cO_X)$, where $|X|$ is the underlying topological space of $X$, and $W_n\cO_X$ is the sheaf of length $n$ truncated Witt vectors. Let $W_n\pi:W_nX\ra \Spec W_nk$ be the morphism induced from $\pi$ via functoriality. 
In this article, our aim is to arrive at a triangle
$$\Z^c_{X}/p^n\ra (W_n\pi)^!W_nk \xra{C'-1}  (W_n\pi)^!W_nk \xra{+1}$$
in the derived category $D^b(X,\Z/p^n)$, in either the \'etale topology, or the Zariski topology with an extra $k=\bar k$ assumption.
Here $(-)^!$ is the extraordinary inverse image functor in the coherent setting as defined in \cite[VII.3.4]{Hartshorne-RD}\cite[(3.3.6)]{Conrad-GDBC}, and $(W_n\pi)^!W_nk$ is a dualizing complex for coherent sheaves on $W_nX$.
This is a generalization of the top degree case of \cite[ 8.3]{GeisserLevine}, 
which in particular implies the above triangle in the smooth case.
Our work is clearly inspired by Kato's paper \cite{Kato-DualitypPrimaryEtale-II}, but the proofs in this article do not use any 
results from {\em loc. cit.}

Let us briefly recall Kato's work in \cite{Kato-DualitypPrimaryEtale-II} and introduce our main object of interest, the complex $K_{n,X,log}$. 
According to Grothendieck's coherent duality theory, there exists an explicit Zariski complex $K_{n,X}$ of quasi-coherent sheaves representing the dualizing complex $(W_n\pi)^!W_nk$ (such a complex $K_{n,X}$ is called a residual complex, cf. \cite[ VI.3.1]{Hartshorne-RD}). 
There is a natural Cartier operator $C':K_{n,X}\ra K_{n,X}$, which is compatible with the classical Cartier operator $C:W_n\Omega^d_X\ra W_n\Omega_X^d$ in the smooth case via Ekedahl's quasi-isomorphism (see \Cref{Compatibility of $C'_n$ and $C$: Prop}). Here $W_n\Omega^d_X$ denotes the degree $d:=\dim X$ part of the $p$-typical de Rham-Witt complex. We define the complex $K_{n,X,log}$ to be the mapping cone of $C'-1$. 
Kato considered in \cite{Kato-DualitypPrimaryEtale-II} the $\FRP$ counterpart, where $\FRP$ is the "flat and relatively perfect" topology (this is a topology with \'etale coverings and with the underlying category lying in between the small and the big \'etale site). 
He  showed that $K_{n,X,log}$ serves in the $\FRP$ topology as a dualizing complex in a rather big triangulated subcategory of the derived category of $\Z/p^n$-sheaves, containing all coherent sheaves and the logarithmic de Rham-Witt sheaves \cite[ 0.1]{Kato-DualitypPrimaryEtale-II}. 
Kato also showed that in the smooth setting, $K_{n,X,log}$ is concentrated in one degree and its only nonzero cohomology sheaf is the top degree logarithmic de Rham-Witt sheaf \cite[ 3.4]{Kato-DualitypPrimaryEtale-II}. 
An analogoue of the latter statement holds naturally on the small \'etale site.
R\"ulling later observed that with a trick from $p^{-1}$-linear algebra, \cite[ 3.4]{Kato-DualitypPrimaryEtale-II} can be done on the Zariski site as well, as long as one assumes $k=\bar k$ (cf. \Cref{Raynaud-Illusie}). Comparing this with the Kato-Moser complex $\tilde \nu_{n,X}$ (cf. \Cref{section Gersten complex of logarithmic de Rham-Witt sheaves}), which is precisely the Gersten resolution of the logarithmic de Rham-Witt sheaf in the smooth setting, one gets an identification in the smooth setting $\tilde \nu_{n,X}\simeq K_{n,X,log}$ in the Zariski topology. Similar as in \cite[ 4.2]{Kato-DualitypPrimaryEtale-II} (cf. \Cref{Kato's complex thickening iso}), 
R\"ulling also built up the localization sequence for $K_{n,X,log}$ on the Zariski site in his unpublished notes (cf. \Cref{localization triangle KnXlog}). 
Compared with the localization sequence for $\Z^c_{X}$ \cite[ 1.1]{Bloch-MovingLemmaForHigherChow} and for $\tilde \nu_{n,X}$ (which trivially holds in the Zariski topology), it is reasonable to expect a chain map relating these objects in general.

The aim of this article is to construct a quasi-isomorphism $\bar\zeta_{log}:\tilde \nu_{n,X}\xra{\simeq} K_{n,X,log}$, for a possibly singular $k$-scheme $X$.
When pre-composed with Zhong's quasi-isomorphism $\bar \psi:\Z^c_{X}/p^n\ra \tilde \nu_{n,X}$ \cite[ 2.16]{Zhong14}, we therefore
obtain another perspective of Bloch's cycle complex with $\Z/p^n$-coefficients in terms of coherent dualizing complexes. More precisely, we prove the following result.

\bthm[{\Cref{Main theorem}, \Cref{Main theorem 2}}]
Let $X$ be a separated scheme of finite type over a perfect field $k$ of positive characteristic $p$. Then there exists a chain map
$$\bar \zeta_{log,\et}:\tilde \nu_{n,X,\et}\xra{\simeq} K_{n,X,log,\et},$$
which is a quasi-isomorphism. If moreover $k$ is algebraically closed, then this chain map induces a quasi-isomorphism on the Zariski site
$$\bar \zeta_{log,\Zar}:\tilde \nu_{n,X,\Zar}\xra{\simeq} K_{n,X,log,\Zar}.$$

Composition with  Zhong's quasi-isomorphism $\bar \psi$, yields the chain map
$$
\bar\zeta_{log,\et}\circ \bar \psi_{\et}: \Z^c_{X,\et}/p^n\xra{\simeq} K_{n,X,log,\et}
$$
which is a quasi-isomorphism. If moreover $k=\bar k$, then the composition
$$
\bar\zeta_{log,\Zar}\circ \bar\psi_{\Zar}:
\Z^c_{X,\Zar}/p^n\xra{\simeq} K_{n,X,log,\Zar}
$$
is  a quasi-isomorphism as well.
\ethm

We explain more on the motivation behind the definition of $K_{n,X,log}$. For  a smooth $k$-scheme $X$, the logarithmic de Rham-Witt sheaves can be defined in two ways: either as the subsheaves of $W_n\Omega_X^d$ generated by log forms, or as the invariant part under the Cartier operator $C$.
In the singular case, these two perspectives give two different (complexes of) sheaves.
The first definition can also be done in the singular case, and this was studied by Morrow \cite{Morrow-KandLogHWofFormal}.
For the second definition one has to replace $W_n\Omega_X^d$ by a dualizing complex on $W_nX$:
for this Grothendieck's duality theory yields a canonical and explicit choice, and this is what we have denoted by $K_{n,X}$. And then this method leads naturally to Kato's and also our construction of $K_{n,X,log}$. Now with our main theorem one knows that $\Z^c_X/p^n$ sits in a distinguished triangle
$$\Z^c_{X}/p^n\ra K_{n,X}\xra{C'-1} K_{n,X}\xra{+1}$$
in the derived category $D^b(X,\Z/p^n)$, in either the \'etale topology, or the Zariski topology with the extra  assumption $k=\bar k$.
In particular, if $X$ is Cohen-Macaulay of pure dimension $d$, then the triangle above becomes
$$\Z^c_{X}/p^n\ra W_n\omega_X[d]\xra{C'-1} W_n\omega_X[d]\xra{+1}.$$
where $W_n\omega_X$ is the only non-vanishing cohomology sheaf of $K_{n,X}$ (if $n=1$, $W_1\omega_X=\omega_X$ is the usual dualizing sheaf on $X$), and $\Z^c_{X}/p^n$ is concentrated at degree $-d$ (cf. \Cref{triangle prop ZcX KnX KnX}).

As corollaries, we arrive at some properties of the higher Chow groups of $0$-cycles with $p$-primary torsion coefficients. (We have specialized several statements here in the introduction part. Please see the main text for more general statements.)
\bcor[\Cref{CH_0=C'-1 invariant}, \Cref{semisimplicity cor}, \Cref{combine with JSS},
\Cref{Vanishing higher Chow}, 
\Cref{etale descent vanishing?}, \Cref{invariance under rational resolution},
\Cref{Galois descent cor}]\label{TotCor}
Let $X$ be a separated scheme of finite type over a perfect field $k$ of characteristic $p>0$ and $\pi:X\ra k$ be the structure map of $X$.
\benu
\item
(Cartier invariance)
Assume $k=\bar k$. Then
$$\CH_0(X,q;\Z/p^n)
=H^{-q}(W_nX,K_{n,X,\Zar})^{C'_\Zar-1}.$$
\item 
(Semisimplicity)
Assume $k=\bar k$. 
Let $X$ be proper over $k$. Then for any $q$,
$$H^{-q}(W_nX,K_{n,X})_{\mathrm{ss}}=\CH_0(X,q;\Z/p^n)\otimes_{\Z/p^n} W_nk.$$
(We refer to \Cref{def of sigma-linear generalized} and \Cref{def of sigma-linear generalized rmk}(2) for the definition of the 
semi-simple part.)
\item 
(Relation with $p$-torsion Poincar\'e duality)
There is an isomorphism in $D^b(X_\et, \Z/p^n)$
$$K_{n,X,log,\et} \simeq R\pi^!(\Z/p^n),$$
where $R\pi^!$ is the extraordinary inverse image functor defined in \cite[Expos\'e XVIII, Thm 3.1.4]{SGA4-3}.
\item
(Affine vanishing)
Assume $k=\bar k$. 
Suppose $X$ is affine and Cohen-Macaulay of pure dimension $d$. Then
$$\CH_0(X,q,\Z/p^n)=0$$ for $q\neq d$.
\item
(\'Etale descent)
Assume $k=\bar k$. 
Suppose $X$ is Cohen-Macaulay of pure relative dimension $d$. Then
$$R^i\epsilon_*(\Z^c_{X,\et}/p^n)
=R^i\epsilon_*\tilde \nu_{n,X,\et}=0,\quad i\neq -d.$$
\item
(Invariance under rational resolution)
Assume $k=\bar k$. 
For a rational resolution of singularities $f:\tilde X\ra X$ (cf. \Cref{def RatSing}) of an integral $k$-scheme $X$ of pure dimension, the trace map induces an isomorphism
$$\CH_0(\tilde X,q;\Z/p^n)\xra{\simeq}\CH_0(X,q;\Z/p^n)$$
for each $q$.
\item 
(Galois descent)
Assume $k=\bar k$. 
Let $f:Y\ra X$ be a finite \'etale Galois map with Galois group $G$. Then 
$$\CH_0(X,d;\Z/p^n) = \CH_0(Y,d;\Z/p^n)^G.$$
\eenu
\ecor

Now we give a more detailed description of the structure of this article.

The general setting is that $X$ is a separated scheme of finite type over a perfect field $k$ of positive characteristic $p$ (except in \Cref{subsection Residual complexes and GD}, where a scheme refers to a noetherian scheme of finite Krull dimension).
In Part 1, we review the basic properties of the chain complexes to appear. \Cref{section Kato's complex KnXlog} is devoted to the properties of the complex $K_{n,X,log}$, the most important object of our studies. We study the Zariski version in Sections 
\ref{subsection on Def of KnXlog} - \ref{subsection functoriality}. Following an idea in \cite{Kato-DualitypPrimaryEtale-II}, we define the Cartier operator $C'$ for the residual complex $K_{n,X}$, and then define the complex $K_{n,X,log}$ to be the mapping cone of $C'-1$ in \Cref{subsection on Def of KnXlog}. We compare our $C'$ with the classical definition of the Cartier operator $C$ for top degree de Rham-Witt sheaves in \Cref{subsection Comparison of WnOmegadXlog with KnXlog}. 
The necessary computations are presented in the subsections \ref{subsubsection: Cartier operator for top Witt differentials} and \ref{subsubsection: Trace map on topp Witt differentials}. The localization sequence is discussed in \Cref{subsection localization triangle associated to KnXlog}. 
The main ingredients are a surjectivity result for $C'-1$, which needs the base field $k$ to be algebraically closed, see \Cref{Raynaud-Illusie} 
(see also \Cref{appendix: sigma-linear algebra} for a short discussion of the necessary  semilinear algebra),
the trace map of a nilpotent thickening (cf. \Cref{Kato's complex thickening iso}), and the localization sequence (cf. \Cref{localization triangle KnXlog}). They were  already observed by R\"ulling and are only re-presented here by the author. After a short discussion on functoriality in \Cref{subsection functoriality}, we move to the \'etale case in \Cref{subsection etale counterpart KnXlog}.
Most of the properties mentioned above  continue to hold in a similar manner, moreover 
the surjectivity of $C_\et-1: W_n\Omega^d_{X,\et}\ra W_n\Omega_{X,\et}^d$ over a smooth $k$-scheme $X$ only requires $k$  to be perfect.
This enables us to build the quasi-isomorphism $\zeta_{log,\et}$ for any perfect field $k$ in the next part. 
In the remaining sections \ref{section bloch's cycle complex} - \ref{section Gersten complex of logarithmic de Rham-Witt sheaves} 
of  Part 1 we recall Bloch's cycle complex $\Z^c_X$, Kato's complex of Milnor $K$-theory $C^M_{X,t}$, 
and the Kato-Moser complex of logarithmic de Rham-Witt sheaves $\tilde \nu_{n,X,t}$.
There are no new results in these three short sections.

In Part 2 we construct the quasi-isomorphism $\bar \zeta_{log}:\tilde \nu_{n,X}\xra{\simeq} K_{n,X,log} $ and study its properties in \Cref{section construction of chain map zeta}. We first construct a chain map $\zeta_{}:C^M_{X}\ra K_{n,X}$ and then show that it induces a chain map $\zeta_{log}:C^M_{X}\ra K_{n,X,log}$. This map actually factors through the chain map $\bar \zeta_{log}:\tilde \nu_{n,X}\ra K_{n,X,log}$ via the Bloch-Gabber-Kato isomorphism \cite[2.8]{Bloch-Kato-pAdicEtaleCohomology}. We prove that $\bar \zeta_{log}$ is a quasi-isomorphism for $t=\et$, and also for $t=\Zar$ with an extra $k=\bar k$ assumption. In \Cref{section combine psi with zeta}, we review the main results of \cite[\S2]{Zhong14} and compose Zhong's quasi-isomorphism $\bar\psi:\Z^c_X/p^n\ra \tilde \nu_{n,X}$ with  $\bar\zeta_{log}$. 
This composite map enables us to use tools from the coherent duality theory in the calculation of certain higher Chow groups of $0$-cycles.

In Part 3 we discuss the applications. \Cref{section De Rham-Witt analysis} mainly serves as a preparatory  section for \Cref{section Higher Chow groups of zero cycles}. In \Cref{section Higher Chow groups of zero cycles} we arrive at several results for higher Chow groups of $0$-cycles with $p$-primary torsion coefficients: affine vanishing, finiteness (reproof of a theorem of Geisser), \'etale descent, and invariance under rational resolutions.

\textit{Acknowledgments}.
This paper is adapted from my PhD thesis. I'd like to express the deepest gratitude to my advisor, Kay R\"ulling,
for suggesting me this topic, for providing a lecture on the de Rham-Witt theory, for sharing his private manuscript and observations, and for the numerous discussions and guidance during the whole time of my PhD. 
In particular, \Cref{Raynaud-Illusie}, \Cref{Kato's complex thickening iso}, and \Cref{localization triangle KnXlog} were already contained in his unpublished manuscript.
Alexander Schmitt has read the preliminary version of this paper with great care and provided detailed comments on the mathematical contents, grammatical errors and typos. 
Thomas Geisser has also given detailed comments to the preliminary version of this work. In particular, he pointed out that \Cref{etale descent ZcX} is in fact a corollary of \cite[8.4]{GeisserLevine}. Yun Hao kindly provided me his notes on an elementary proof of \Cref{1-T surjective prop in appendix}. 
The anonymous referee has written a report with many detailed suggestions. In particular, \Cref{Cn computation for A^d} and \Cref{Trace map on top differentials: Example} are rewritten following the comments there.
S\'andor Kov\'acs informed me about an update in the normality condition in \cite{Kovacs-RatSing}. 
Accordingly, the statement of \Cref{invariance under rational resolution} was changed, and \Cref{invariance under rational resolution rmk}(3) was added. (In particular, the statement in \Cref{TotCor}(6) remains valid.)
I am indebted to them all.
I thank also the Berlin mathematical school for providing financial support during my studies in Berlin.

\part{The complexes}\label{part 1}
\section{Kato's complex $K_{n,X,log,t}$}\label{section Kato's complex KnXlog}

\subsection{Preliminaries: Residual complexes and Grothendieck's duality theory}\label{subsection Residual complexes and GD}
The general references for this topic are \cite{Hartshorne-RD} and \cite{Conrad-GDBC}.
All schemes in \Cref{subsection Residual complexes and GD} will be assumed to be noetherian of finite Krull dimension.
\subsubsection{Residual complexes}
A \emph{residual complex} (\cite[p.125]{Conrad-GDBC}, \cite[p.304]{Hartshorne-RD}) on a scheme $X$ is a complex $K$ such that
\begin{itemize}
\item
$K$ is bounded as a complex,
\item
all the terms of $K$ are quasi-coherent and injective $\cO_X$-modules,
\item
the cohomology sheaves are coherent, and
\item
there is an isomorphism of $\cO_X$-modules
$$\bigoplus_{q\in \Z}K^q \simeq \bigoplus_{x\in X}i_{x*}J(x),$$
where $i_x:\Spec \cO_{X,x}\ra X$ is the canonical map and $J(x)$ is the quasi-coherent sheaf on $\Spec \cO_{X,x}$ associated to an injective hull of $k(x)$ over $\cO_{X,x}$ (i.e. 
the unique injective $\cO_{X,x}$-module up to non-unique isomorphisms which contains $k(x)$ as a submodule and such that, for any $0\neq a\in J(x)$, there exists an element $b\in \cO_{X,x}$ with $0\neq ba\in k(x)$.)
\end{itemize}

Given a residual complex $K$ on $X$ and a point $x\in X$, there is a unique integer $d_{K}(x)$, such that $i_{x*}J(x)$ is a direct summand of $K^q$, i.e.,
$$K^q\simeq \bigoplus_{d_{K}(x)=q} i_{x*}J(x).$$
The assignment $x\mapsto d_{K}(x)$ is called the \emph{codimension function on $X$ associated to $K$} (cf. \cite[IV, 1.1(a)]{Hartshorne-RD}\cite[p.125]{Conrad-GDBC}).
We define the \emph{associated filtration}
$$Z^\bullet(K)=\{x\in X\mid d_{K}(x)\ge p \}.$$
On each irreducible component of $X$, this filtration equals the shifted codimension filtration.
By the \emph{codimension filtration} of a scheme $X$ we refer to the filtration $Z^\bullet$ with
$$Z^p=\{x\in X\mid \dim\cO_{X,x}\ge p\}.$$
If $Z^\bullet$ is a filtration on $X$, we denote by $Z^\bullet[n]$ the shifted filtration with $Z^\bullet[n]^p = Z^{p+n}$.

Let $Z^\bullet$ be a filtration on $X$ such that when restricted to each irreducible component, it is the shifted codimension filtration.
For any bounded below complex $\cF^\bullet$, choose a bounded below injective resolution $\cI^\bullet$ of $\cF^\bullet$.
Denote by $\uline\Gamma_{Z^p}$ the sheafified local cohomology functor with support in $Z^p$, cf. \cite[p223 5.]{Hartshorne-RD}.
Then one has a natural decreasing exhaustive filtration by subcomplexes of $\cI^\bullet$:
$$\dots \supset \uline\Gamma_{Z^p}(\cI^\bullet)
\supset \uline\Gamma_{Z^{p+1}}(\cI^\bullet)\supset\dots.$$
This filtration is stalkwise bounded below. Now consider the $E_1$-spectral sequence associated to this filtration 
$$E_1^{p,q}\Rightarrow H^{p+q}(\cF^\bullet).$$
The \emph{Cousin complex} (\cite[p.105]{Conrad-GDBC}) $E_{Z^\bullet}(\cF^\bullet)$ associated to $\cF^\bullet$ is defined to be the $0$-th line of the $E_1$-page, namely
$$E_{Z^\bullet}(\cF^\bullet):=(E_1^{p,0}=\cH^p_{Z^p/Z^{p+1}}(\cF), d_1^{p,0}).$$
Here $\cH^p_{Z^p/Z^{p+1}}(\cF):=R^p\uline\Gamma_{Z^p/Z^{p+1}}(\cF)$ and $\uline\Gamma_{Z^p/Z^{p+1}}(\cF):=\uline\Gamma_{Z^p}(\cF)/\uline\Gamma_{Z^{p+1}}(\cF)$ (cf. \cite[p.225 Variation 7]{Hartshorne-RD}).
We will also use the shortened notation $E$ for $E_{Z^\bullet}$ when the filtration $Z^\bullet$ is clear from the context.
Note that $E_{Z^\bullet}(\cF^\bullet)$ is indeed a Cousin complex in the sense of \cite[p.105]{Conrad-GDBC} by the canonical functorial isomorphism \cite[p.226]{Hartshorne-RD}\cite[(3.1.4)]{Conrad-GDBC}
$$
\cH^{i}_{Z^{p}/Z^{p+1}}(\cF^\bullet) \xra{\simeq} \bigoplus_{x\in Z^p-Z^{p+1}}i_{x*}( H^{i}_{x}(\cF^\bullet)),
$$
where $i_x:\Spec\cO_{X,x}\ra X$ is the canonical map, and $H^{i}_{x}(\cF^\bullet)$ is the local cohomology groups at $x$ as defined in \cite[p.225 Variation 8]{Hartshorne-RD}. By slight abuse of notation we denote by the same notation $H^{i}_{x}(\cF^\bullet)$ the quasi-coherent sheaf on $\Spec \cO_{X,x}$ associated to this local cohomology group, and it is a sheaf supported on the closed point if it is nonzero.

Let $X$ be a scheme and $Z^\bullet$ be a filtration on $X$ which is a shift of the codimension filtration on each irreducible component of $X$.
Denote by $Q$ the natural functor from the category of complexes of $\cO_X$-modules to the derived category of $\cO_X$-modules. Then $E_{Z^\bullet}$ and $Q$ induce quasi-inverses (\cite[3.2.1]{Conrad-GDBC})
\beq\label{cateq}
\left(
  \begin{array}{c}
    \text{dualizing complexes whose} \\
    \text{associated filtration is $Z^\bullet$}
  \end{array}
\right)
\xleftrightarrows[E_{Z^\bullet}]{Q}
\left(
  \begin{array}{c}
    \text{residual complexes whose} \\
    \text{associated filtration is $Z^\bullet$}
  \end{array}
\right).
\eeq
For the definition of a dualizing complex (as an object in the derived category) we refer the reader to \cite[p.118]{Conrad-GDBC}.
Since we have assumed that $X$ is noetherian and of finite Krull dimension, there always exists a residual complex on $X$.

\subsubsection{The functor $f^\triangle$}\label{section functor f upper triangle}
Let $f:X\ra Y$ be a finite type morphism between noetherian schemes of finite Krull dimension and let $K$ be a residual complex on $Y$ with associated filtration $Z^\bullet:=Z^\bullet(K)$ and codimension function $d_K$.
Define the function $d_{f^\triangle K}$ on $X$ to be (\cite[(3.2.4)]{Conrad-GDBC})
$$d_{f^\triangle K}(x):=d_{K}(f(x))-\trdeg(k(x)/k(f(x))$$ (so far the subscript $f^\triangle K$ is simply regarded as a formal symbol), and define $f^\triangle Z^\bullet$ accordingly
$$f^\triangle Z^\bullet=\{x\in X\mid d_{f^\triangle K}(x)\ge p \}.$$
Notice that if $f$ has constant fiber dimension $r$, $f^\triangle Z^\bullet$ is simply $f^{-1} Z^\bullet[r]$.

Following \cite[VI, 3.1]{Hartshorne-RD}\cite[3.2.2]{Conrad-GDBC}, we list some properties of the functor $f^\triangle$ below.

\bprop\label{f triangle}
There exists a functor
$$f^\triangle:
\left(
  \begin{array}{c}
    \text{residual complexes on $Y$} \\
    \text{with filtration $Z^\bullet$} \\
  \end{array}
\right)
\ra
\left(
  \begin{array}{c}
    \text{residual complexes on $X$} \\
    \text{with filtration $f^\triangle Z^\bullet$} \\
  \end{array}
\right)
$$
having the following properties (we assume all schemes are noetherian schemes of finite Krull dimension, and all morphisms are of finite type).
\begin{enumerate}
\item
If $f$ is finite, there is an isomorphism of complexes (\cite[VI.3.1]{Hartshorne-RD})
$$
\psi_f:f^\triangle K \xra{\simeq}
E_{f^{-1}Z^\bullet}(\bar f^*R\cHom_{\cO_Y}(f_*\cO_X,K))\simeq
\bar f^*\cHom_{\cO_Y}(f_*\cO_X,K),
$$
where $\bar f^*:=f^{-1}(-)\otimes_{f^{-1}f_*\cO_X}\cO_X$ is the pullback functor associated to the map of ringed spaces $\bar f:(X,\cO_X)\ra (Y,f_*\cO_X)$. Since $\bar f$ is flat, the pullback functor $\bar f^*$ is exact. The last isomorphism is due to the fact that $\bar f^*\cHom_{\cO_Y}(f_*\cO_X,K)$ is a residual complex with respect to filtration $f^{-1}Z^\bullet$ (see \cite[VI.4.1]{Hartshorne-RD}, \cite[(3.4.5)]{Conrad-GDBC}).

\item
If $f$ is smooth and separated of relative dimension $r$, there is an isomorphism of complexes (\cite[VI.3.1]{Hartshorne-RD})
$$\varphi_f: f^\triangle K\xra{\simeq}
E_{f^{-1}Z^\bullet[r]}(\Omega_{X/Y}^r[r]\otimes_{\cO_X}^L Lf^*K)
=E_{f^{-1}Z^\bullet[r]}(\Omega_{X/Y}^r[r]\otimes_{\cO_X} f^*K).
$$
The last equality is due to the flatness of $f$ and local freeness of $\Omega_{X/Y}^r$.

If $f$ is \'etale (or more generally residually stable, see (\ref{residually stable base change item}) below), this becomes
$$\varphi_f: f^\triangle K\xra{\simeq}
E_{f^{-1}Z^\bullet}(f^*K)
\simeq f^*K.
$$
The last isomorphism is due to \cite[VI.5.3]{Hartshorne-RD}.
In particular, if $f=j:X\hra Y$ is an open immersion, $j^\triangle K= j^*K$ is a residual complex with respect to filtration $X\cap Z^\bullet$ (\cite[p.128]{Conrad-GDBC}).

\item
If $f$ is finite \'etale, the chain maps $\psi_f, \varphi_f$ are compatible. Namely, for a given residual complex $K$ on $Y$, there exists an isomorphism of complexes
$\bar f^*\cHom_{\cO_Y}(f_*\cO_X,K)\xra{\simeq} f^*K$
as defined in \cite[(2.7.9)]{Conrad-GDBC}, such that the following diagram of complexes commutes
$$
\xymatrix@C=0.1cm@R=0.5cm{
&*+[r]{\bar f^*\cHom_{\cO_Y}(f_*\cO_X,K)
}
\ar[dd]^{\simeq}
\\
f^\triangle(K) \ar[ur]_{\simeq}^{\psi_f} \ar[dr]^{\simeq}_{\varphi_f} \\
&*+[r]{f^*K.}}
$$

\item
(Composition) If $f:X\ra Y$ and $g:Y\ra Z$ are two such morphisms, there is an natural isomorphism of functors (\cite[(3.2.3)]{Conrad-GDBC})
$$c_{f,g}:(gf)^\triangle \xra{\simeq} f^\triangle g^\triangle.$$

\item
\label{residually stable base change item}
(Residually stable base change) Following \cite[p.132]{Conrad-GDBC}, we say a (not necessarily locally finite type) morphism $f:X\ra Y$ between locally noetherian schemes is \emph{residually stable} if
\bit
\item
$f$ is flat,
\item
the fibers of $f$ are discrete and for all $x\in X$, the extension $k(x)/k(f(x))$ is algebraic, and
\item
the fibers of $f$ are Gorenstein schemes.
\eit
As an example, an \'etale morphism is residually stable. For more properties of residually stable morphisms, see \cite[VI, \S5]{Hartshorne-RD}. Let $f$ be a morphism of finite type, and $u$ be a residually stable morphism. Let
\begin{equation}\label{cartesion diagram residually stable}
\xymatrix{
  X' \ar[d]_{f'} \ar[r]^{u'}
                & X \ar[d]^{f}  \\
  Y' \ar[r]_{u}
                & Y   }
\end{equation}
be a cartesian diagram. Then there is an natural transformation between functors (\cite[VI.5.5]{Hartshorne-RD})
$$d_{u,f}: f'^\triangle u^*\xra{\simeq} u'^*f^\triangle.$$
\item
$f^\triangle$ is compatible with translation and tensoring with an invertible sheaf. More precisely, for an invertible sheaf $\cL$ on $Y$ and a locally constant $\Z$-valued function $n$ on $Y$, one has canonical isomorphisms of complexes \cite[(3.3.9)]{Conrad-GDBC}
$$f^\triangle(\cL[n]\otimes K)\simeq
(f^*K)[n]\otimes f^\triangle K\simeq
(f^*\cL\otimes f^\triangle K)[n].$$
\end{enumerate}
\eprop
More properties and compatibility diagrams can be found in \cite[\S3.3]{Conrad-GDBC} and \cite[VI, \S3, \S5]{Hartshorne-RD}.

\subsubsection{Trace map for residual complexes}
\label{subsubsection trace map for residual complexes}

\bprop\label{trace for residual}
Let $f:X\ra Y$ be a \emph{proper} morphism between noetherian schemes of finite Krull dimensions and let $K$ be a residual complex on $Y$. Then there exists a map of complexes
$$\Tr_f:f_*f^\triangle K\ra K,$$
such that the following properties hold (\cite[\S3.4]{Conrad-GDBC}).
\begin{enumerate}
\item
If $f$ is finite, 
$\Tr_f$ at a given residual complex $K$ agrees with the following composite as a map of complexes (\cite[(3.4.8)]{Conrad-GDBC}):
\begin{equation}
\label{trace map for residual complexes rmk (1)}
\Tr_f: f_*f^\triangle K\xra[\simeq]{\psi_f}
\cHom_{\cO_Y}(f_*\cO_X,K)\xra{\text{ev. at 1}} K.
\end{equation}
\item
If $f:\P^d_Y\ra Y$ is the natural projection, then the trace map $\Tr_f$ at $K$, as a map in the derived category $D_c^b(Y)$, agrees with the following composite (\cite[p.151]{Conrad-GDBC})
$$
f_*f^\triangle K
\xra[\simeq]{\varphi_f}
Rf_*(\Omega_{\P^n_X/X}^n[n])\otimes_{\cO_Y}K
\xra{} K.$$
The first map is induced from $\varphi_f$ followed by the projection formula (\cite[(2.1.10)]{Conrad-GDBC}), and the second map is induced by base change from the following isomorphism of groups(\cite[(2.3.3)]{Conrad-GDBC})
$$\Z\xra{\simeq}H^d(\P^d_\Z,\Omega_{\P^d_\Z/\Z}^d)
=\check{H}^d(\ssU,\Omega_{\P^d_\Z/\Z}^d),
\quad
1\mapsto (-1)^{\frac{d(d+1)}{2}}
\frac{dt_1\wedge\dots\wedge dt_d}{t_1\dots t_d},
$$
where $\ssU=\{U_0, ..., U_d\}$ is the standard covering of $\P_\Z^d$ and the $t_i$'s are the coordinate functions on $U_0$.
\item
(Functoriality, \cite[3.4.1(1)]{Conrad-GDBC})
$\Tr_f$ is functorial with respect to residual complexes with the same associated filtration.
\item
(Composition, \cite[3.4.1(2)]{Conrad-GDBC})
If $g:Y\ra Z$ is another proper morphism, then
$$\Tr_{gf} = \Tr_g \circ g_*(\Tr_f )\circ(gf)_*c_{f,g}.$$

\item
(Residually stable base change, \cite[VI.5.6]{Hartshorne-RD})
Notations are the same as in diagram (\ref{cartesion diagram residually stable}), and we assume $f$ proper and $u$ residually stable. Then the diagram
$$
\xymatrix@C=5em{
u^*Rf_*f^\triangle\ar[r]^{u^*\Tr_f}\ar[d]^{\simeq}& u^*\\
Rf'_*u'^*f^\triangle\ar[r]_{\simeq}^{Rf'_*(d_{u,f})}&
Rf'_*f'^\triangle u^*\ar[u]^{\Tr_{f'}}
}
$$
commutes.
\item
$\Tr_f$ is compatible with translation and tensoring with an invertible sheaf (\cite[p.148]{Conrad-GDBC}).
\item
(Grothendieck-Serre duality, \cite[3.4.4]{Conrad-GDBC}) If $f:X\ra Y$ is proper, then for any $\cF\in D^-_\qc(X)$, the composition
$$Rf_*R\cHom_X(\cF,f^\triangle K)
\ra R\cHom_Y(Rf_*\cF,Rf_*f^\triangle K)\xra{\Tr_f}
R\cHom_Y(Rf_*\cF,K)$$
is an isomorphism in $D^+_\qc(Y)$.
\end{enumerate}
\eprop

More properties and compatibility diagrams can be found in \cite[\S3.4]{Conrad-GDBC} and \cite[VI, \S4-5; VII, \S2]{Hartshorne-RD}.

\subsection{Definition of $K_{n,X,log}$}\label{subsection on Def of KnXlog}
Let $k$ be a perfect field of characteristic $p$. Let $W_nk$ be the ring of Witt vectors of length $n$ of $k$.
Notice that $W_nk$ is an injective $W_nk$-module by Baer's criterion. So $\Spec W_nk$ is a Gorenstein scheme by \cite[V. 9.1(ii)]{Hartshorne-RD}, and its structure sheaf placed at degree 0 is a residual complex 
(with codimension function being the zero function and the associated filtration being $Z^\bullet(W_nk)=\{Z^0(W_nk)\}$, where $Z^0(W_nk)$ is the set of the unique point in $\Spec W_nk$) 
by \cite[p299 1.]{Hartshorne-RD} and the categorical equivalence (\ref{cateq}) (note that in this case the Cousin functor $E_{Z^\bullet(W_nk)}$ applied to $W_nk$ is still $W_nk$). This justifies the symbol $(W_nF_k)^\triangle$ to appear.
To avoid possible confusion we will distinguish the source and target of the absolute Frobenius by using the symbols $k_1=k_2=k$. Absolute Frobenius is then written as $F_k:(\Spec k_1,k_1)\ra (\Spec k_2,k_2)$, and the $n$-th Witt lift is written as $W_nF_k:(\Spec W_nk_1,W_nk_1)\ra (\Spec W_nk_2,W_nk_2)$.
There is a natural isomorphism of $W_nk_1$-modules (the last isomorphism is given by \Cref{f triangle}(1))
\begin{align}\label{(1.0.1)}
W_nk_1 &\xra{\simeq}
\bar{W_nF_k}^*\Hom_{W_nk_2}((W_nF_{k})_*(W_nk_1),W_nk_2)
\simeq
(W_nF_k)^\triangle (W_nk_2),
\\
a&\mapsto a^{}\otimes [(W_nF_{k})_*1\mapsto 1]
\quad(=[(W_nF_{k})_*a\mapsto 1]),
\nonumber
\end{align}
where $\bar{W_nF_k}:(\Spec W_nk_1,W_nk_1)\ra (\Spec W_nk_2,(W_nF_{k})_* (W_nk_1))$ is the natural map of ringed spaces, and the Hom set is given the $(W_nF_{k})_* (W_nk_1)$-module structure via the first place. 
In fact, 
it is clearly a bijection: identify the target with $W_nk_2$ via the evaluate-at-1 map, then one can see that the map (\ref{(1.0.1)}) is identified with $a\mapsto (W_nF_k)^{-1}(a)$.

Let $X$ be a separated scheme of finite type over $k$ of dimension $d$ with structure map $\pi:X\ra k$. 
Since $W_nk$ is a Gorenstein scheme as we recalled in the last paragraph,
$$K_{n,X}:=(W_n\pi)^\triangle W_nk$$
is a residual complex on $W_nX$, associated to the codimension function $d_{K_{n,X}}$ with $$d_{K_{n,X}}(x)=-\dim \bar{\{x\}},$$ 
and the filtration $Z^\bullet(K_{n,X})=\{Z^p(K_{n,X})\}$ with $$Z^p(K_{n,X})=\{x\in X\mid \dim\bar{\{x\}}\le -p \}.$$ 
In particular, $K_{n,X}$ is a bounded complex of injective quasi-coherent $W_n\cO_X$-modules with coherent cohomologies sitting in degrees 
$$[-d,0].$$ 
If $n=1$, we write $K_X:=K_{1,X}$. 
Now we turn to the definition of $C'$. 
Denote the level $n$ Witt lift of the absolute Frobenius $F_X$ by $W_nF_X:(W_nX_1,W_n\cO_{X_1})\ra (W_nX_2,W_n\cO_{X_2})$.
The structure maps of $W_nX_1,W_nX_2$ are $W_n\pi_1,W_n\pi_2$ respectively. These schemes fit into a commutative diagram
$$\xymatrix{
W_nX_1\ar[d]_{W_n\pi_1}\ar[r]^{W_nF_X}&
W_nX_2\ar[d]^{W_n\pi_2}\\
\Spec W_nk_1\ar[r]^{W_nF_k}& \Spec W_nk_2.
}$$
Denote
$$K_{n,X_i}:=(W_n\pi_i)^\triangle (W_nk_i), \quad i=1,2.$$
Via functoriality, one has a $W_n\cO_{X_1}$-linear map
\begin{align}\label{(1.0.2)}
K_{n,X_1}= (W_n\pi_1)^\triangle (W_nk_1)
&\xra[\simeq]{(W_n\pi_1)^\triangle(\ref{(1.0.1)})}
(W_n\pi_1)^\triangle (W_nF_{k}) ^\triangle (W_n k_2)
\\
&\qquad\simeq (W_nF_X)^\triangle (W_n\pi_2)^\triangle (W_n k_2)
\simeq(W_nF_X)^\triangle K_{n,X_2}.
\nonumber
\end{align}
Here the isomorphism at the beginning of the second line is given by \Cref{f triangle}(4). 
Then via the adjunction with respect to the morphism $W_nF_X$, one has a $W_n\cO_{X_2}$-linear map
\begin{align}\label{C'_n def}
C':=C'_n:(W_nF_X)_*K_{n,X_1}
\xra[\simeq]{(W_nF_X)_*(\ref{(1.0.2)})} (W_nF_X)_*(W_nF_X)^\triangle K_{n,X_2}
\xra{\Tr_{W_nF_X}} K_{n,X_2},
\end{align}
where the last map is the trace map of $W_nF_X$ for residual complexes.
We call it the (level $n$) \emph{Cartier operator} for residual complexes.
We sometimes omit the $(W_nF_X)_*$-module structure of the source and write simply as $C':K_{n,X}\ra K_{n,X}$.

Now we come to the construction of $K_{n,X,log}$ (cf. \cite[\S3]{Kato-DualitypPrimaryEtale-II}). Define
\begin{equation}\label{def of K log}
K_{n,X,log}:=\Cone (K_{n,X}\xra{C'-1}K_{n,X})[-1].
\end{equation}
This is a complex of abelian sheaves sitting in degrees 
$$[-d,1].$$ 
If $n=1$, we set $K_{X,log}:=K_{1,X,log}$.
Writing more explicitly, $K_{n,X,log}$ is the following complex
$$(K_{n,X}^{-d}\oplus 0) \ra
(K_{n,X}^{-d+1} \oplus K_{n,X}^{-d})\ra \dots \ra
(K_{n,X}^{0}\oplus K_{n,X}^{-1})\ra
(0\oplus K_{n,X}^0).$$
The differential of $K_{n,X,log}$ at degree $i$ is given by
\begin{align*}
d_{log}=d_{n,log}:K_{n,X,log}^{i}&\ra K_{n,X,log}^{i+1}\\
(K_{n,X}^{i} \oplus K_{n,X}^{i-1})&\ra (K_{n,X}^{i+1} \oplus K_{n,X}^{i}) \\
(a,b)&\mapsto (d(a),-(C'-1)(a)-d (b)),
\end{align*}
where $d$ is the differential in $K_{n,X}$.
The sign conventions we adopt here for shifted complexes and the cone construction are the same as in \cite[p6, p8]{Conrad-GDBC}.  And naturally, one has a distinguished triangle
\begin{equation}\label{distinguished triangle of KX C-1}
K_{n,X,log}\xra{}K_{n,X}\xra{C'-1}K_{n,X}\xra{+1} K_{n,X,log}[1].
\end{equation}
Explicitly, the first map is in degree $i$ given by
\begin{align*}
K_{n,X,log}^i=K_{n,X}^i\oplus K_{n,X}^{i-1}&\xra{}K_{n,X}^i,\\
(a,b)&\mapsto a.
\end{align*}
The "$+1$" map is given by
\begin{align*}
K_{n,X}^i&\xra{}(K_{n,X,log}[1])^i=K_{n,X,log}^{i+1}=
(K_{n,X}^{i+1}\oplus K_{n,X}^{i}),\\
b&\mapsto (0,b).
\end{align*}
Both maps are indeed maps of chain complexes.

\subsection{Comparison of $W_n\Omega^d_{X,log}$ with $K_{n,X,log}$}\label{subsection Comparison of WnOmegadXlog with KnXlog}
Recall the following result from the classical Grothendieck duality theory \cite[IV.3.4]{Hartshorne-RD}\cite[3.1.3]{Conrad-GDBC} and Ekedahl \cite[\S1]{Ekedahl-MultPropertiesOfDeRhamWittI} (see also \cite[proof of 1.10.3 and Rmk. 1.10.4]{CR12}).
\bprop[Ekedahl]\label{Ekedahl qis}
If $X$ is smooth and of pure dimension $d$ over $k$, then there is a canonical quasi-isomorphism
$$W_n\Omega_X^d[d]\xra{\simeq}K_{n,X}.$$
\eprop

\brmk\label{partial Ekedahl for singular}
Suppose $X$ is a separated scheme of finite type over $k$ of dimension $d$. Denote by $U$ the smooth locus of $X$, and 
suppose that the complement $Z$ of $U$ is of dimension $e$. Suppose moreover that $U$ is non-empty and equidimensional (it is satisfied for example, if $X$ is integral). 
Then Ekedahl's quasi-isomorphism \Cref{Ekedahl qis} gives a quasi-isomorphism of dualizing complexes
\beq\label{WnOmegaUd qis to KnU} 
W_n\Omega_U^d[d]\xra{\simeq}K_{n,U}.
\eeq 
Note that the associated filtrations of quasi-isomorphic dualizing complexes are the same (cf. \cite[V.3.4]{Hartshorne-RD}). 
Let $Z^\bullet$ be the codimension filtration of $U$.
As explained above, the associated filtration of $K_{n,U}$ is the shifted codimension filtration, 
i.e., $Z^\bullet[d]$.
Apply the Cousin functor associated to the shifted codimension filtration $Z^\bullet[d]$ to the quasi-isomorphism (\ref{WnOmegaUd qis to KnU}) between dualizing complexes, we have an isomorphism of residual complexes
$$E_{Z^\bullet[d]}(W_n\Omega_U^d[d])\xra{\simeq}K_{n,U}$$
with the same filtration $Z^\bullet[d]$ by (\ref{cateq}).
Since $W_nj$ is an open immersion, we can canonically identify the residual complexes $(W_nj)^*K_{n,X}\simeq K_{n,U}$ by \Cref{f triangle}(2). 
Since $K_{n,X}$ is a residual complex and in particular is a Cousin complex (cf. \cite[p. 105]{Conrad-GDBC}), the adjunction map $K_{n,X}\ra (W_nj)_*(W_nj)^*K_{n,X}\simeq (W_nj)_*K_{n,U}$ is an isomorphism at degrees $[-d,-e-1]$. Thus the induced chain map
$$K_{n,X}\ra (W_nj)_*E_{Z^\bullet[d]}(W_n\Omega_U^d[d])$$ 
is an isomorphism at degrees $[-d,-e-1]$. 
\ermk

\subsubsection{Compatibility of $C'$ with the classical Cartier operator $C$}
\label{subsubsection compatibility of C' with Cartier C}

We review the absolute Cartier operator in the classical literature (see e.g. \cite[Chapter 1 \S3]{BrionKimar-BookFrobSplit}, \cite[\S0.2]{Illusie-DeRhamWitt}, \cite[7.2]{Katz-NilpConn}, \cite[III \S1]{Illusie-Raynaud}).
Let $X$ be a $k$-scheme. The \emph{(absolute) inverse Cartier operator $\gamma_X$} of degree $i$ on a scheme $X$ is affine locally, say, on $\Spec A\subset X$, given additively by the following expression ($\cH^i(-)$ denotes the cohomology sheaf of the complex)
\begin{align}
\gamma_A:\quad\qquad\Omega_{A/k}^i\quad\quad
&\xra{}\quad\quad
\cH^i(F_{A,*}\Omega_{A/k}^\bullet)
\label{def of C inverse: local expression}
\\
ada_1 \dots   da_i &
\mapsto a^p a^{p-1}_1da_1 \dots  a_i^{p-1}da_i,
\nonumber
\end{align}
where $a,a_1,\dots a_i\in A$. Here $\cH^i(F_{A,*}\Omega_{A/k}^\bullet)$ denotes the $A$-module structure on $\cH^i(\Omega_{A/k}^\bullet)$ via the absolute Frobenius $F_A:A\ra A, a\mapsto a^p$ (note that $F_{A,*}\Omega_{A/k}^\bullet$ is a complex of $A$-modules in positive characteristics). For each degree $i$, $\gamma_A$ thus defined is an $A$-linear map. These local maps patch together and give rise to a map of sheaves
\begin{equation}\label{inverse Cartier operator}
\gamma_X:
\Omega_X^i
\xra{}
\cH^i(F_{X,*}\Omega_X^\bullet)
\end{equation}
which is $\cO_X$-linear. 
If $X$ is smooth of dimension $d$, $\gamma_X$ is an isomorphism of $\cO_X$-modules, which is called the \emph{(absolute) Cartier isomorphism}. See \cite[1.3.4]{BrionKimar-BookFrobSplit} for a proof (note that although the authors there assumed the base field to be algebraically closed, the proof of this theorem works for any perfect field $k$ of positive characteristics).

This can be generalized to the de Rham-Witt case.
\ble[{cf. \cite[4.1.3]{Kato-DualitypPrimaryEtale-I}}]\label{bar F is an iso for de rham witt}
Denote by $W_n\Omega_X'^i$ the abelian sheaf $F(W_{n+1}\Omega_X^i)$ regarded as a $W_n\cO_X$-submodule of $(W_nF_X)_*W_n\Omega_X^i$.
If $X$ is smooth of dimension $d$, the map
$$\bar F:W_n\Omega_X^i\ra W_n\Omega_X'^i/dV^{n-1}\Omega_X^{i-1}$$
induced by Frobenius $F:W_{n+1}\Omega_X^i\ra R_*(W_nF_X)_*W_n\Omega_X^i$ is an isomorphism of $W_n\cO_X$-modules.

In particular, if $i=d$,
$$\bar F:W_n\Omega_X^d\ra (W_nF_X)_* W_{n}\Omega_X^d/dV^{n-1}\Omega_X^{d-1}$$
is an isomorphism of $W_n\cO_X$-modules.
\ele
\bpf
Since
$$\Ker (R:W_{n+1}\Omega^i\ra W_n\Omega^i)=V^n\Omega^i+dV^n\Omega^{i-1},$$
$FV^n\Omega^i=0$ and $FdV^n\Omega^{i-1}=dV^{n-1}\Omega^{i-1}$, $F:W_{n+1}\Omega^i\ra W_n\Omega^i$ reduces to
$$\bar F:W_{n}\Omega^i\ra W_n\Omega^i/dV^{n-1}.$$
The surjectivity is clear. We show the injectivity. Suppose
$x\in W_{n+1}\Omega^i$, $y\in \Omega^{i-1}$, such that $F(x)=dV^{n-1}y$. Then $F(x-dV^ny)=0$, which implies by \cite[I (3.21.1.2)]{Illusie-DeRhamWitt} that $x-dV^ny\in V^n\Omega^i$.

The second claim follows from the fact that $F:W_{n+1}\Omega^d\ra R_*(W_nF_X)_*W_n\Omega^d$ is surjective on top degree $d$ \cite[I (3.21.1.1)]{Illusie-DeRhamWitt}, and therefore $W_n\Omega'^d=(W_nF_X)_*W_n\Omega^d$ as $W_n\cO_X$-modules.
\epf

\bdefn[(absolute) Cartier operator]\label{def of C}
Let $X$ be a smooth scheme of dimension $d$ over $k$.
\begin{enumerate}
\item
The composition
\begin{align}\label{Cartier operator}
C:=C_X:&Z^i(F_{X,*}\Omega^\bullet_X)\ra \cH^i(F_{X,*}\Omega_X^\bullet)\xra{(\gamma_X)^{-1}} \Omega_X^i
\\
\big(\text{with } &Z^i(F_{X,*}\Omega^\bullet_X):=\Ker(F_{X,*}\Omega^i_X\xra{d} F_{X,*}\Omega^{i+1}_X)\big)\nonumber
\end{align}
is called the \emph{(absolute) Cartier operator} of degree $i$, denoted by $C$ or $C_X$.
\item
(cf. \cite[4.1.2, 4.1.4]{Kato-DualitypPrimaryEtale-I})
More generally,
for $n\ge 1$, define the \emph{(absolute) Cartier operator $C_n:=C_{n,X}$ of level $n$} to be the composite
\begin{equation}\label{level n cartier def}
C_n:W_n\Omega_X'^i\twoheadrightarrow W_n\Omega_X'^i/dV^{n-1}\Omega_X^{i-1}\xra[\simeq]{\bar F^{-1}}
W_n\Omega_X^i,
\end{equation}
where
$\bar F:W_{n}\Omega_X^i\xra{\simeq} W_n\Omega_X'^i/dV^{n-1}\Omega_X^{i-1}$
is the map in \Cref{bar F is an iso for de rham witt}. 
If $i=d$ is the top degree we obtain the $W_n\cO_X$-linear
map
\begin{equation}\label{level n cartier def, top degree}
C_n: (W_nF_X)_*W_{n}\Omega_X^d\twoheadrightarrow (W_nF_X)_*W_{n}\Omega_X^d/dV^{n-1}\Omega_X^{d-1}\xra{\bar F^{-1}}
W_n\Omega_X^d.
\end{equation}
\end{enumerate}
\edefn

\brmk
\begin{enumerate}
\item\label{illusie I.3.3}
According to the explicit formula for $F$, we have $C=C_1$ \cite[I 3.3]{Illusie-DeRhamWitt}. For this reason we will simply write $C$ for $C_n$ sometimes.
\item
$C_n$ (for all $n$) are compatible with \'etale pullbacks. Actually any de Rham-Witt system (e.g. $(W_n\Omega^\bullet_X,F,V,R,\uline p, d)$) is compatible with \'etale base change \cite[1.3.2]{CR12}.
\item
The $n$-th power of Frobenius $F$ induces a map
$$\bar F^n:W_n\Omega_X^i\xra{\simeq} \cH^i((W_nF_X)^n_*W_n\Omega_X^\bullet),$$
which is the same as \cite[III (1.4.1)]{Illusie-Raynaud}.
\item
Notice that on $\Spec W_nk$, $C_n:W_nk\ra W_nk$ is simply the map $(W_nF_k)^{-1}$, because $F:W_{n+1}k\ra W_nk$ equals $R\circ W_{n+1}F_k$ in characteristic $p$.
\item
We sometimes omit "$(W_nF_{X})_*$" in the source. But one should always keep that in mind and be careful with the module structure.
\end{enumerate}
\ermk

\brmk\label{etale over WnX}
Before we move on, we state a remark on \'etale schemes over $W_nX$.
\benu
\item
Notice that every \'etale $W_nX$-scheme is of the form $W_ng:W_nU\ra W_nX$, where $g:U\ra X$ is an \'etale $X$-scheme. In fact, there are two functors 
\begin{align*}
    F:\{\text{\'etale $W_nX$-schemes}\}
    &\leftrightarrows\{\text{\'etale $X$-schemes}\}:G
    \\
    V&\mapsto V\times_{W_nX} X
    \\
    W_nU&\mapsfrom U
\end{align*}
The functor $F$ is a categorical equivalence according to \cite[Ch. IV, 18.1.2]{EGAIV-4}. The functor $G$ is well-defined (i.e. produces \'etale $W_nX$-schemes) and is a right inverse of $F$ by \cite[Thm. 1.25]{Hesselholt-BigDeRhamWittComplex}. We want to show that there is a natural isomorphism $GF\simeq id$, and this is the consequence of the following purely categorical statement: If $F:\cA\ra \cB$ and $G:\cB\ra \cA$ are two functors satisfying both $F$ being a categorical equivalence and $FG\simeq id$, then $G$ is a quasi-inverse of $F$, i.e., there exists a canonical natural isomorphism $GF\simeq id$. We leave this as an easy exercise for the reader.

\item
The square
$$\xymatrix{
W_nU\ar[d]_{W_ng}\ar[r]^{W_nF_U}&
W_nU\ar[d]^{W_ng}\\
W_nX\ar[r]^{W_nF_X}&W_nX.
}$$
is a cartesian square. This is because for any \'etale map $g:U\ra X$, the relative Frobenius $F_{U/X}$ is an isomorphism by \cite[10.3.1]{Fu-EtaleBook}. Thus $W_nF_{U/X}$ is an also isomorphism and the claim follows.
\eenu
\ermk

We shall now state the main result in this subsection, which seems to be known by experts (cf. proof of \cite[3.4]{Kato-DualitypPrimaryEtale-II}) but we cannot find a proof in the literature.
To eliminate possible sign inconsistency of the Cartier operator with the Grothendieck trace map calculated via residue symbols \cite[Appendix A]{Conrad-GDBC}, we give a proof by explicit calculations (see \Cref{subsubsection: Cartier operator for top Witt differentials}-\Cref{subsubsection: Trace map on topp Witt differentials}). 
At the same time, this result justifies our notation for $C'$: The classical Cartier operator $C$ is simply the $(-d)$-th cohomology of our $C'$ in the smooth case.
\bthm[Compatibility of $C'$ with $C$]\label{Compatibility of $C'_n$ and $C$: Prop}
Suppose that $X$ is a smooth scheme of dimension $d$ over a perfect field $k$ of characteristic $p>0$. Then the top degree classical Cartier operator
$$C:(W_nF_X)_*W_n\Omega^d_{X/k}\ra W_n\Omega^d_{X/k}$$
as defined in \Cref{def of C}, agrees with the $(-d)$-th cohomology of the Cartier operator for residual complexes
$$C': (W_nF_X)_*W_n\Omega^d_{X/k}\ra W_n\Omega^d_{X/k}$$
as defined in (\ref{C'_n def}) via Ekedahl's quasi-isomorphism \Cref{Ekedahl qis}.
\ethm 

\bpf
The Cartier operator is stable under \'etale base change, i.e., for any \'etale morphism $W_ng:W_nX\ra W_nY$ (which must be of this form according to \Cref{etale over WnX}(1)), we have 
$$C_X\simeq (W_ng)^*C_Y: (W_nF_X)_*W_n\Omega_X^d\ra W_n\Omega_X^d.$$
We claim that the map $C'$ defined in (\ref{C'_n def}) is also compatible with \'etale base change. That is, whenever we have an \'etale morphism $W_ng: W_nX\ra W_nY$, there is a canonical isomorphism 
$$C'_X\simeq (W_ng)^*C'_Y: (W_nF_X)_*K_{n,X}\ra K_{n,X}.$$
First of all, the Grothendieck trace map $\Tr_{W_nF_X}$ for residual complexes is compatible with \'etale base change by \Cref{trace for residual}(5), i.e., 
$$\Tr_{W_nF_X}\simeq g^*\Tr_{W_nF_Y}: (W_nF_X)_*(W_nF_X)^\triangle K_{n,X}\ra K_{n,X}.$$
Secondly, because of the cartesian square in \Cref{etale over WnX}(2) and the flat base change theorem
$$(W_ng)^*(W_nF_X)_*\simeq (W_nF_X)^*(W_ng)_*,$$
we are reduced to show that (\ref{(1.0.2)}) is compatible with \'etale base change. And this is true, because we have
$$(W_ng)^*\simeq (W_ng)^\triangle$$
by \Cref{f triangle}(2), and the compatibility of $(-)^\triangle$ with composition by \Cref{f triangle}(4). This finishes the claim.

Note that the question is local on $W_nX$.
Thus to prove the statement for smooth $k$-schemes $X$, using the compatibility of $C$ and $C'$ with respect to \'etale base change, it suffices to prove it for $X=\A^d_k$. That is, we need to check that the expression given in \Cref{Trace map on top differentials: Example} for $C'$ agrees with the expression for $C$ given in \Cref{Cn computation for A^d}. This is apparent.
\epf

\subsubsection{Proof of \Cref{Compatibility of $C'_n$ and $C$: Prop}: $C$ for the top Witt differentials on the affine space}
\label{subsubsection: Cartier operator for top Witt differentials}
Let $k$ be a perfect field of positive characteristic $p$. The aim of this subsection is to provide the formula for the Cartier operator on the top degree de Rham-Witt sheaf over the affine space (\Cref{Cn computation for A^d}).

Consider the polynomial ring $k[X_1,...,X_d]$.
Let $h:\{1,...,d\}\ra \N[\frac{1}{p}]\setminus\{0\}$ be a function such that $\mathrm{Im}(p^{n-1}h)\subset \N$.
Write $h_i:=h(i)$. 
Let $\{i_1,...,i_d\}$ be a reordering of $\{1,...,d\}$, such that 
$$v_p(h_{i_1})\le v_p(h_{i_2})\le...\le v_p(h_{i_d}).$$
This order depends on $h$. 
Since $\{1,...,d\}$ is a finite set, we can also choose a uniform order for elements in $\Supp h$ and $\Supp p^ah$ for any integer $a$ and any function $h$.
If $\mathrm{Im}(h)\not\subset \N$, let $r\in [1,d]$ be the unique integer such that
$$v_p(h_{i_1})\le...\le v_p(h_{i_r})<0\le
v_p(h_{i_{r+1}})\le...\le v_p(h_{i_d}).$$
For all $j\in [1,d]$, write
$$v_j:=v_p(h_{i_j}),\quad h_j':=h_{i_j}p^{-v_j}.$$
According to \cite[2.17]{LangerZink}, any element in $W_n\Omega_{k[X_1,...,X_d]}^d$ is uniquely written as a sum of (\ref{FormNonint}) and (\ref{FormInt}):
\bit
\item 
$h$ is a function such that $\mathrm{Im}(h)\not\subset \N$, $\alpha\in W_{n+v_1}k$,
\begin{align}\label{FormNonint}
dV^{-v_1}(\alpha[X_1]^{h_1'})
\cdots 
dV^{-v_r}([X_r]^{h_r'})\cdot 
F^{v_{r+1}}d[X_{r+1}]^{h_{r+1}'}
 \cdots 
F^{v_d}d[X_d]^{h_d'},
\end{align}
and 
\item 
$h$ is a function such that $\mathrm{Im}(h)\subset \N$, $\beta\in W_nk$,
\begin{align}\label{FormInt}
& \beta 
F^{v_1}d[X_1]^{h_1'} 
 \cdots 
F^{v_d}d[X_d]^{h_d'}.
\end{align}
\eit

\ble[$C_n$ on $\A^d$]
\label{Cn computation for A^d}
The Cartier operator (cf. \Cref{def of C})
$$C:=C_n:W_n\Omega_{k[X_1,...,X_d]}^d\ra W_n\Omega_{k[X_1,...,X_d]}^d$$
is the map uniquely determined by the following assignment by taking products: 
for $j\in[1,d]$,
\begin{align*}
\begin{array}{c}
\\
dV^{-v_j}(\alpha[X_j]_{n+v_j}^{h_j'})\\
(\alpha\in W_{n+v_j}k, v_j<0)
\end{array}
&\mapsto 
\left\{
  \begin{array}{ll}
    dV^{1-v_j}(R(\alpha)[X_j]_{n-1+v_j}^{h_j'}), &
       \hbox{$-v_j\in [1,n-2]$;} \\
    0, & \hbox{$-v_j=n-1$.}
  \end{array}
\right.
\\
\begin{array}{c}
\\
\beta\\
(\beta\in W_nk)
\end{array}
&\mapsto 
(W_nF_k)^{-1}(\beta),
\\
\begin{array}{c}
\\
F^{v_{j}}d[X_{j}]_{n+v_j}^{h_j'}\\
(v_j\ge 0)
\end{array}
&\mapsto
\left\{
  \begin{array}{ll}
    F^{v_{j}-1}d[X_{j}]_{n+v_j-1}^{h_j'}, & \hbox{$v_j\ge 1$;} \\
    dV[X_{j}]_{n-1}^{h_{i_j}}, & \hbox{$v_j=0$.}
  \end{array}
\right.
\nonumber
\end{align*}
\ele

\bpf
For any $\gamma,\delta\in W_n\Omega_{k[X_1,...,X_d]}^\bullet$, $C(F(\gamma)\cdot F(\delta))=C(F(\gamma))\cdot C(F(\delta)).$
Hence it suffices to check the formulae on each factor.
For $\alpha\in W_{n+v_j}k, -v_j\in[1,n-2]$, the formulae  $CF=R$, $d=FdV$ imply
\begin{align*}
    C\Big(dV^{-v_j}(\alpha[X_j]_{n+v_j}^{h_j'})\Big)
    &=C\Big(FdV^{1-v_j}(\alpha[X_j]_{n+v_j}^{h_j'})\Big)
    \\
    &=dV^{1-v_j}(R(\alpha)[X_j]_{n-1+v_j}^{h_j'}).
\end{align*}
For $\alpha \in k, -v_j=n-1$, the formulae  $CF=R$, $d=FdV$ and $\Ker R=V^n+dV^n$ imply
\begin{align*}
    C\Big(dV^{-v_j}(\alpha[X_j]_{n+v_j}^{h_j'})\Big)
    &=C\Big(FdV^{n}(\alpha X_j^{h_{i_j}})\Big)
    \\
    &=R\Big(dV^n(\alpha X_j^{h_{i_j}})\Big)
    \\
    &=0.
\end{align*}
For $\beta\in W_nk$,
\begin{align*}
    C(\beta)=(W_nF_k)^{-1}(\beta).
\end{align*}
For $v_j\ge 1$, the formula $CF=R$ implies
\begin{align*}
    C\Big(F^{v_{j}}d[X_{j}]_{n+v_j}^{h_j'}\Big)
    &=R\Big(F^{v_{j}-1}d[X_{j}]_{n+v_j}^{h_j'}\Big)
    \\&=F^{v_{j}-1}d[X_{j}]_{n+v_j-1}^{h_j'}
\end{align*}
For $v_j=0$, the formulae $CF=R$ and $d=FdV$ imply
\begin{align*}
    C\Big(F^{v_{j}}d[X_{j}]_{n+v_j}^{h_j'}\Big)
    &=C\Big(FdV[X_{j}]_{n}^{h_{i_j}}\Big)\\
    &=R\Big(dV[X_{j}]_{n}^{h_{i_j}}\Big)\\
    &=dV[X_{j}]_{n-1}^{h_{i_j}}
\end{align*}
\epf

\subsubsection{Proof of \Cref{Compatibility of $C'_n$ and $C$: Prop}: $C'$ for the top Witt differentials on the affine space}\label{subsubsection: Trace map on topp Witt differentials}
The aim of \Cref{subsubsection: Trace map on topp Witt differentials} is to calculate $C'$ for the top de Rham-Witt sheaves on the affine space (\Cref{Trace map on top differentials: Example}). To do this, one needs to first calculate the trace map of the canonical lift of the absolute Frobenius.

\subsubsubsection{Trace map of the canonical lift $\tilde F_{\tilde X}$ of absolute Frobenius $F_X$}
\label{Trace map of the canonical lift of Frobenius} 
 
Before we start with the computation we recall some properties of the residue symbol from \cite[\S A]{Conrad-GDBC}. Let $X\hra P$ be a closed immersion of affine schemes with the sheaf of ideals generated by  $t_1,...,t_d\in \Gamma(P,\cO_P)$. Let $P\ra Y$ be a separated smooth morphism of affine schemes with pure relative dimension $d$. 
Suppose $X\ra Y$ is finite flat.  
For any $\omega\in \Gamma(P,\omega_{P/Y})$, there is a well-defined element
$$\Res_{P/ Y}
\left[
  \begin{array}{c}
   \omega
   \\
    t_1,\dots, t_d \\
  \end{array}
\right]
\in \Gamma( Y,\cO_{Y})$$
which is called the residue symbol (cf. \cite[(A.1.4)]{Conrad-GDBC}). It satisfies the following properties (we use the same numbering as in \cite[\S A.1]{Conrad-GDBC}):
\bit 
\item 
Suppose $h:Y'\ra Y$ is any morphism of schemes, and $P'=P\times_Y Y'$. Then 
\beq\tag{R5}
\Res_{P'/ Y'}
\left[
  \begin{array}{c}
   h^*\omega
   \\
    h^*t_1,\dots, h^*t_d \\
  \end{array}
\right]
=
h^*
\Res_{P/ Y}
\left[
  \begin{array}{c}
   \omega
   \\
    t_1,\dots, t_d \\
  \end{array}
\right]
\eeq 
\item 
For any $\varphi\in \Gamma(P,\cO_P)$,
\beq\tag{R6}
\Res_{P/ Y}
\left[
  \begin{array}{c}
 \varphi\cdot dt_1\dots dt_d
   \\
    t_1,\dots, t_d \\
  \end{array}
\right]
=\Tr_{X/Y}(\varphi|_X).
\eeq 
Here the notation $\Tr_{{X}/{Y}}$ denotes the classical trace map associated to the finite locally free ring extension $\Gamma(Y,\cO_{Y})\ra \Gamma(X,\cO_{X})$.
\item 
For $\eta\in \Gamma(P,\Omega_{P/Y}^{n-1})$, and $k_1,\dots, k_d$ positive integers, 
\beq\tag{R9}
\Res_{P/Y}
\left[
  \begin{array}{c}
    d\eta \\
    t_1^{k_1},\dots, t_d^{k_d} \\
  \end{array}
\right]
=
\sum_{i=1}^n
k_i\cdot 
\Res_{P/Y}
\left[
  \begin{array}{c}
   dt_i\wedge \eta
   \\
   t_1^{k_1},\dots, t_i^{k_i+1},\dots, t_d^{k_d}  \\
  \end{array}
\right].
\eeq 

\eit 

Let $k$ be a perfect field of positive characteristic $p$. Let $X=\A^d_k$, and denote by $\tilde X:=\Spec W_n(k)[X_1,\dots, X_d]$  the canonical smooth lift of $X$ over $W_n(k)$. 
To make the module structures in the following discussion explicit, we distinguish the source and the target of the absolute Frobenius of $\Spec k$ and write it as
$$F_k:\Spec k_1\ra \Spec k_2.$$
Similarly, write the absolute Frobenius on $X$ as 
$$F_X:X=\Spec k_1[X_1,\dots, X_d]\ra Y=\Spec k_2[Y_1,\dots, Y_d].$$
There is a canonical lift  of $F_X$ over $\tilde X$
$$\tilde F_{\tilde X}: \tilde X=\Spec W_n(k_1)[X_1,\dots,X_d] \ra \tilde Y:=\Spec W_n(k_2)[Y_1,\dots,Y_d].$$
given by
\begin{align*}
\tilde F_{\tilde X}^*:\quad
\Gamma(\tilde Y,\cO_{\tilde Y})=W_n(k_2)[Y_1,\dots,Y_d]&\ra W_n(k_1)[X_1,\dots,X_d]=\Gamma(\tilde X,\cO_{\tilde X}),\\
W_nk_2\ni
\alpha &\mapsto W_n(F_k)(\alpha),\\
Y_i&\mapsto X_i^p.
\end{align*}
Let 
$$\pi_X:X\ra \Spec k_1,\quad\pi_Y:Y\ra\Spec k_2,\quad\pi_{\tilde X}: \tilde X\ra W_nk_1, \quad\pi_{\tilde Y}:\tilde Y\ra W_nk_2$$ 
be the structure maps. 
The composition $\pi_{\tilde Y}\circ \tilde F_{\tilde X} : \tilde X\ra \Spec W_nk_2$ gives $\tilde X$ a $W_nk_2$-scheme structure, and the map $\tilde F_{\tilde X}$ is then a map of $W_nk_2$-schemes. Therefore the trace map 
$$ 
\Tr_{\tilde F_{\tilde X}}:
\tilde F_{\tilde X,*}\tilde F_{\tilde X}^\triangle K_{\tilde Y}
\ra K_{\tilde Y}$$ 
makes sense. Consider the following map of complexes
\begin{align*}
\tilde F_{\tilde X,*}K_{\tilde X}
\simeq \tilde F_{\tilde X,*}\pi_{\tilde X}^{\triangle} W_nk_1
\xra[\sim]{\tilde F_{\tilde X,*}\pi_{\tilde X}^\triangle(\ref{(1.0.1)})}
&\tilde F_{\tilde X,*} 
\pi_{\tilde X}^{\triangle} W_nF_{k}^\triangle  W_nk_2
\simeq
\\
&\tilde F_{\tilde X,*}\tilde F_{\tilde X}^\triangle 
\pi_{\tilde Y}^{\triangle} W_nk_2
\simeq
\tilde F_{\tilde X,*}\tilde F_{\tilde X}^\triangle K_{\tilde Y}
\xra{\Tr_{\tilde F_{\tilde X}}} 
K_{\tilde Y}.
\nonumber
\end{align*}
Taking the $(-d)$-th cohomology, it induces a map
\beq\label{eq:trace of lift Frob on Coh}
\tilde F_{X,*}\Omega_{\tilde X/W_nk_1}^d\ra \Omega_{\tilde Y /W_nk_2}^d
\eeq

\ble\label{Trace map on top differentials: Example on lifting}
The notations are the same as above. The map (\ref{eq:trace of lift Frob on Coh}) has the following expression:
\begin{align}
\label{Trace map on top differentials: in particular for Ad}
\quad
\Omega_{\tilde X/W_nk_1}^d\quad 
&\xra{(\ref{eq:trace of lift Frob on Coh})} \quad \Omega_{\tilde Y /W_nk_2}^d\\
\alpha \bm X^{\bm \lambda+p\bm\mu}d \bm X
&\mapsto
\left\{
  \begin{array}{ll}
    (W_nF_k)^{-1}(\alpha) \bm Y^{\bm\mu}d \bm Y, &
    \hbox{if $\lambda_i = p-1$ for all $i$;} \\
    0,
&\hbox{if $\lambda_i\neq p-1$ for some $i$.}
  \end{array}
\right.
\nonumber
\end{align}
\ele

\bpf
Consider the closed immersion $i:\tilde X\hra \tilde P=\A^d_{\tilde Y}$ associated to the following homomorphism of rings:
\begin{align*}
\Gamma(\tilde P,\cO_{\tilde P})=
W_n(k_2)[Y_1,\dots, Y_d,T_1,\dots, T_d]&\ra W_n(k_1)[X_1,\dots, X_d]=\Gamma(\tilde X,\cO_{\tilde X}), \\
\alpha&\mapsto (W_nF_k)(\alpha),\quad \alpha\in W_n(k_2),\\
Y_i&\mapsto X_i^p,\quad i=1,\dots, d,\\
T_i&\mapsto X_i,\quad i=1,\dots, d.
\end{align*}
Its kernel is
$$I=(T_1^p-Y_1,\dots, T_d^p-Y_d).$$
Denote
$$t_i=T_i^p-Y_i , \quad i=1,\dots, d.$$
Obviously the $t_i$'s form a regular sequence in $\Gamma(\tilde P,\cO_{\tilde P})$, and hence $i$ is a regular immersion. Then one has a factorization of $\tilde F_{\tilde X}$:
\begin{equation}\label{Trace map on top differentials: diagram in example}
\xymatrix{
{\tilde X=\Spec W_n(k_1)[X_1,\dots, X_d]}
\ar[dr]_{\tilde F_{\tilde X}}\ar@{^(->}[r]^(0.45){i} &
{\tilde P=\Spec W_n(k_2)[Y_1,\dots, Y_d,T_1,\dots, T_d]}\ar[d]^{\pi}\\
&{\tilde Y=\Spec W_n(k_2)[Y_1,\dots,Y_d]}.
}
\end{equation}
Regarding $\tilde X$ as a $W_nk_2$-scheme via the composite map $\tilde F_{\tilde X}\circ \pi_{\tilde Y}$, the diagram (\ref{Trace map on top differentials: diagram in example}) is then a diagram in the category of $W_nk_2$-schemes.

A general element in $\Gamma(\tilde X, \Omega_{\tilde X/W_nk_1}^d)$ is a sum of expressions of the form
\beq\label{Trace map on top differentials: element}
\alpha\bm X^{\bm \lambda+p\bm\mu} d\bm X,
\quad \alpha\in W_nk_1, \bm\lambda\in[0,p-1]^d, \bm\mu\in\N^d.
\eeq
Here $\bm \lambda=\{\lambda_1,\dots,\lambda_d\}$, $\bm\mu=\{\mu_1,\dots,\mu_d\}$ are multi-indices, and $\bm X^{\bm \lambda}:=X_1^{\lambda_1}\dots X_d^{\lambda_d}$ (similar for $\bm Y^{\bm \mu}$, $\bm X^{\bm \lambda+p\bm\mu}$, etc.), $d\bm X:=dX_1\dots dX_d$ (similar for $d\bm T$, etc.). 
The element (\ref{Trace map on top differentials: element}) in $\Gamma(\tilde X, \Omega_{\tilde X/W_nk_1}^d)$  corresponds to
\beq\label{Trace map on top differentials: element2}
(W_nF_k)^{-1}(\alpha)\bm X^{\bm \lambda+p\bm\mu} d\bm X,
\quad \alpha\in W_nk_2, \bm\lambda\in[0,p-1]^d, \bm\mu\in\N^d,
\eeq
in $\Gamma(\tilde X, \Omega_{\tilde X/W_nk_2}^d)$ under $(-d)$-th cohomology of the map $\tilde F_{\tilde X,*}\pi_{\tilde X}^\triangle(\ref{(1.0.1)})$, and
$$(W_nF_k)^{-1}(\alpha)\bm T^{\bm \lambda}\bm Y^{\bm\mu} d\bm T,
\quad \alpha\in W_nk_2, \bm\lambda\in[0,p-1]^d, \bm\mu\in\N^d$$
is a lift of (\ref{Trace map on top differentials: element2}) to $\Gamma(\tilde P, \Omega_{\tilde P/W_nk_2}^d)$. 
Write
\begin{align*}
\beta:&=dt_d\wedge\dots\wedge dt_1\wedge
(W_nF_k)^{-1}(\alpha)\bm T^{\bm \lambda}\bm Y^{\bm\mu} d\bm T
\\&=
(-1)^d
dY_d\wedge\dots\wedge dY_1\wedge
(W_nF_k)^{-1}(\alpha)\bm T^{\bm \lambda}\bm Y^{\bm\mu} d\bm T
\end{align*}
in $\Gamma(\tilde P, \omega_{\tilde P/W_nk_2})$, where $\omega_{\tilde P/W_nk_2}$ denotes the dualizing sheaf with respect to the smooth morphism $\tilde P\ra W_nk_2$. 
Recall that there is a natural isomorphism (\cite[p.30 (a)]{Conrad-GDBC})
\begin{align*}
\omega_{\tilde P/W_nk_2}&\simeq
\omega_{\tilde P/\tilde Y}\otimes_{\cO_{\tilde P}} \pi^*\omega_{\tilde Y/W_nk_2},
\end{align*}
where $\omega_{\tilde P/\tilde Y}$ and $\omega_{\tilde Y/W_nk_2}$ denote the dualizing sheaves with respect to the smooth morphisms $\pi:\tilde P\ra \tilde Y$ and $\tilde Y\ra W_nk_2$.
This isomorphism maps 
$\beta$ to 
$$(-1)^{\frac{d(3d+1)}{2}}
(W_nF_k)^{-1}(\alpha)\bm T^{\bm \lambda} d\bm T
\otimes
\pi^*\bm Y^{\bm\mu}d\bm Y.$$
It is easily seen that $\tilde F_{\tilde X}$ is a finite flat morphism between smooth $W_nk_2$-schemes.  Applying \cite[Lemma A.3.3]{CR11}, one has
$$\Tr_{\tilde F_{\tilde X}}((W_nF_k)^{-1}(\alpha)\bm X^{\bm \lambda+p\bm\mu} d\bm X)
=
(W_nF_k)^{-1}(\alpha)
\Res_{\tilde P/\tilde Y}
\left[
  \begin{array}{c}
   \bm T^{\bm\lambda} d\bm T
   \\
    t_1,\dots, t_d \\
  \end{array}
\right]
\bm Y^{\bm\mu}d\bm Y,
$$
where $\Res_{\tilde P/\tilde Y}
\left[
  \begin{array}{c}
   \bm T^{\bm\lambda} d\bm T
   \\
    t_1,\dots, t_d \\
  \end{array}
\right]
\in \Gamma(\tilde Y,\cO_{\tilde Y})$
is the residue symbol defined in \cite[(A.1.4)]{Conrad-GDBC}, and $\Tr_{\tilde F_{\tilde X}}$ is the trace map on top differentials of the $W_nk_2$-morphism $\tilde F_{\tilde X}$ \cite[(2.7.36)]{Conrad-GDBC}.

We consider the following cases:
\begin{itemize}
\item
If $(\lambda_1,\dots,\lambda_d)\neq (p-1,\dots,p-1)$, $\bm T^{\bm\lambda} d\bm T=d\eta$ for some $\eta\in \Omega_{\tilde P/\tilde Y}^{d-1}$. Suppose without loss of generality that $\lambda_1\neq p-1$. Then we can take
$$\eta=\frac{1}{\lambda_1+1}T_1^{\lambda_1+1}T_2^{\lambda_2}\dots T_d^{\lambda_d}dT_2\dots dT_d.$$
Noticing that $$dt_i=d(T_i^p-Y_i)=pT_i^{p-1}dT_i$$ in $\Omega^1_{\tilde P/\tilde Y}$, and that $\lambda_1+mp+1$ ($m\in \Z_{>0}$) is not divisible by $p$ if $\lambda_1+1$ is so. Now we calculate
\begin{align*}
\Res_{\tilde P/\tilde Y}
\left[
  \begin{array}{c}
    \bm T^{\bm\lambda} d\bm T \\
    t_1,\dots, t_d \\
  \end{array}
\right]
&=
\frac{1}{\lambda_1+1}
\Res_{\tilde P/\tilde Y}
\left[
  \begin{array}{c}
    d(T_1^{\lambda_1+1}T_2^{\lambda_2}\dots T_d^{\lambda_d}dT_2\dots dT_d) \\
    t_1,t_2,\dots, t_d \\
  \end{array}
\right]
\\
&=
\frac{p}{\lambda_1+1}
\Res_{\tilde P/\tilde Y}
\left[
  \begin{array}{c}
   T_1^{\lambda_1+p}T_2^{\lambda_2}\dots T_d^{\lambda_d}
   dT_1
   dT_2\dots dT_d \\
   t_1^2,t_2,\dots, t_d \\
  \end{array}
\right]
\\
&=
\frac{p}{(\lambda_1+1)(\lambda_1+p+1)}
\Res_{\tilde P/\tilde Y}
\left[
  \begin{array}{c}
   d(T_1^{\lambda_1+p+1}T_2^{\lambda_2}\dots T_d^{\lambda_d}
   dT_2\dots dT_d) \\
   t_1^2,t_2,\dots, t_d \\
  \end{array}
\right]
\\
&=
\frac{2p^2}{(\lambda_1+1)(\lambda_1+p+1)}
\Res_{\tilde P/\tilde Y}
\left[
  \begin{array}{c}
   T_1^{\lambda_1+2p}T_2^{\lambda_2}\dots T_d^{\lambda_d}
   dT_1
   dT_2\dots dT_d \\
    t_1^3,t_2,\dots, t_d \\
  \end{array}
\right]
\\
&=
\frac{2p^2}{\prod_{i=0}^{2}(\lambda_1+ip+1)}
\Res_{\tilde P/\tilde Y}
\left[
  \begin{array}{c}
   d(T_1^{\lambda_1+2p+1}T_2^{\lambda_2}\dots T_d^{\lambda_d}
   dT_2\dots dT_d) \\
   t_1^3,t_2,\dots, t_d \\
  \end{array}
\right]
\\
&=
\frac{6p^3}{\prod_{i=0}^{2}(\lambda_1+ip+1)}
\Res_{\tilde P/\tilde Y}
\left[
  \begin{array}{c}
   T_1^{\lambda_1+3p}T_2^{\lambda_2}\dots T_d^{\lambda_d}
   dT_1
   dT_2\dots dT_d \\
    t_1^4,t_2,\dots, t_d \\
  \end{array}
\right]
\\
&=
\dots
\\
&=
\frac{(\prod_{i=1}^{n}i)\cdot p^n}{\prod_{i=0}^{n-1}(\lambda_1+ip+1)}
\Res_{\tilde P/\tilde Y}
\left[
  \begin{array}{c}
   T_1^{\lambda_1+np}T_2^{\lambda_2}\dots T_d^{\lambda_d}
   dT_1
   dT_2\dots dT_d \\
    t_1^{n+1},t_2,\dots, t_d \\
  \end{array}
\right]
\\
&=
0.
\end{align*}
We have used (R9) on the second, the fourth, the sixth, and the eighth equality signs.
The last equality is because $p^n=0$ in $\Gamma(\tilde Y,\cO_{\tilde Y})$.

\item
If $(\lambda_1,\dots,\lambda_d)= (p-1,\dots,p-1)$,
consider
\begin{equation}\label{Trace map on top differentials: lifted diagram}
\xymatrix{
X':=\Spec \frac{\Z[Y'_1,\dots, Y'_d,T'_1,\dots,T'_d]}{({T'_1}^p-Y'_1,\dots, {T'_d}^p-Y'_d)}
\ar[dr]_{f}\ar@{^(->}[r]^{} &
{\Spec \Z[Y'_1,\dots, Y'_d,T'_1,\dots,T'_d]}=:P'\ar[d]\\
&{\Spec \Z[Y'_1,\dots, Y'_d]}=:Y'.
}
\end{equation}
The map $f$ is given by $f^*(Y'_i)=Y'_i={T'_i}^p$ in $\Gamma(X',\cO_{X'})$. This is a finite locally free morphism of rank $p^d$.
Consider the map $h:\tilde Y\ra Y'$ given by 
\begin{align*}
\Gamma(Y',\cO_{Y'})= \Z[Y'_1,\dots, Y'_d]&\ra W_n(k_2)[Y_1,\dots,Y_d]=\Gamma(\tilde Y,\cO_{\tilde Y}),\\
Y'_i&\mapsto Y_i\text{ for all $i$,}
\end{align*}
that relates the two diagrams (\ref{Trace map on top differentials: lifted diagram}) and (\ref{Trace map on top differentials: diagram in example}).
In $\Gamma(Y',\cO_{Y'})$, we have
\begin{align*}
p^d\cdot
\Res_{P'/Y'}
\left[
  \begin{array}{c}
    {T'_1}^{p-1}\dots {T'_d}^{p-1} dT'_1 \dots   dT'_d \\
    {T'_1}^p-Y'_1,\dots, {T'_d}^p-Y'_d
  \end{array}
\right]
&=
\Res_{P'/Y'}
\left[
  \begin{array}{c}
    d({T'_1}^p-Y'_1) \dots   d({T'_d}^p-Y'_d) \\
     {T'_1}^p-Y'_1,\dots, {T'_d}^p-Y'_d \\
  \end{array}
\right]
\\
&\stackrel{\text{(R6)}}{=}
\Tr_{{X'}/{Y'}}(1)
\\
&=p^d.
\end{align*}
The notation $\Tr_{{X'}/{Y'}}$ denotes the classical trace map associated to the finite locally free ring extension $\Gamma(Y',\cO_{Y'})\ra \Gamma(X',\cO_{X'})$.
As for the last equality, $\Tr_{{X'}/{Y'}}(1)=p^d$ because $f$ is a finite locally free map of rank $p^d$. Since $p^d$ is a non-zerodivisor in $\Gamma(Y',\cO_{Y'})$, one deduces
$$\Res_{P'/Y'}
\left[
  \begin{array}{c}
    {T'_1}^{p-1}\dots {T'_d}^{p-1} dT'_1 \dots   dT'_d \\
    {T'_1}^p-Y'_1,\dots, {T'_d}^p-Y'_d
  \end{array}
\right]
=1.$$
Set
$$\bm T^{\bm{p-1}}=T_1^{p-1}\dots T_d^{p-1},$$
which is the canonical lift of $\bm X^{\bm \lambda}$ via the map $i:\tilde X\hra \tilde P$ in our current case. 
Pulling back to $\Gamma(\tilde Y,\cO_{\tilde Y})$ via $h$, one has
\begin{align}
\Res_{\tilde P/\tilde Y}
\left[
  \begin{array}{c}
   \bm T^{\bm{p-1}}d\bm T \\
    t_1,\dots, t_d \\
  \end{array}
\right]
&\stackrel{\text{(R5)}}{=}
h^*
\Res_{P'/Y'}
\left[
  \begin{array}{c}
    {T'_1}^{p-1}\dots {T'_d}^{p-1} dT'_1 \dots   dT'_d \\
    {T'_1}^p-Y'_1,\dots, {T'_d}^p-Y'_d
  \end{array}
\right]
\stackrel{}{=}
1.
\end{align}
\end{itemize}
Altogether, we know that the map (\ref{eq:trace of lift Frob on Coh}) takes the following expression
\begin{align*}
\Omega_{\tilde X/W_nk_1}^d\quad&\ra \quad \Omega_{\tilde Y /W_nk_2}^d\\
\alpha \bm X^{\bm \lambda+p\bm\mu}d \bm X
&\mapsto
\left\{
  \begin{array}{ll}
    (W_nF_k)^{-1}(\alpha) \bm Y^{\bm\mu}d \bm Y, &
    \hbox{if $\lambda_i = p-1$ for all $i$;} \\
    0,
&\hbox{if $\lambda_i\neq p-1$ for some $i$.}
  \end{array}
\right.
\nonumber
\end{align*}
\epf

\subsubsubsection{$C'$ for top Witt differentials}
Now we turn to the $W_n$-version. The aim of this subsection is to calculate $C'$ for top Witt differentials on $\A^d_k$ (\Cref{Trace map on top differentials: Example}).

Let $f : X\ra Y$ be a finite
morphism between smooth, separated and equidimensional $k$-schemes of dimension $d$. Same as before, we denote by $\pi_X: X\ra k$ and $\pi_Y:Y\ra k$ the respective structure maps. 
The complexes $K_{n,X}:=(W_n\pi_X)^\triangle W_nk, K_{n,Y}=(W_n\pi_Y)^\triangle W_nk$ are residual complexes on $X$ and $Y$. Then we define the trace map
\begin{equation}\label{Trace map on top Witt differentials - just  this}
\Tr_{W_nf}:(W_nf)_*(W_n\Omega^d_{X})\ra W_n\Omega^d_{Y}
\end{equation}
to be the $(-d)$-th cohomology map of the composition
\begin{equation}\label{Trace map on top Witt differentials: where it is induced from}
\Tr_{W_nf}:
(W_nf)_*K_{n,X}\simeq
\cHom_{W_n\cO_Y}((W_nf)_*W_n\cO_X, K_{n,Y})\xra{\text{ev. at 1}} K_{n,Y}
\end{equation}
via Ekedahl's isomorphism $W_n\Omega_{X}^d\simeq \cH^{-d}(K_{n,X})$ in
\Cref{Ekedahl qis}.

It suffices to compute the trace map locally on $Y$. Thus by possibly shrinking $Y$ we can assume that $Y$ and (therefore also $X$) is affine. 
In this case, there exist smooth affine $W_nk$-schemes $\tilde X$ and $\tilde Y$ which lift $X$ and $Y$. Denote the structure morphisms of $\tilde X,\tilde Y$ by $\pi_{\tilde X}$ and $\pi_{\tilde Y}$, respectively. Then there exists a finite $W_nk$-morphism
$\tilde f:\tilde X\ra \tilde Y$
lifting $f:X\ra Y$ by the formal smoothness property of $\tilde Y$.

Consider the map of abelian sheaves \cite[I (2.3)]{Ekedahl-MultPropertiesOfDeRhamWittI}
\begin{align}
\label{theta_Y cO eq}
\varrho_Y^*:\quad
W_n\cO_Y&\xra[\simeq]{\vartheta_Y} \cH^0(\Omega_{\tilde Y/W_nk}^\bullet)\hra
\cO_{\tilde Y},\qquad\qquad
\\
\sum_{i=0}^{n-1}V^i([a_i])&\mapsto \tilde a_0^{p^n}+  p\tilde a_1^{p^{n-1}}+\dots +p^{n-1} \tilde a_{n-1}^{p},
\nonumber
\end{align}
where $a_i\in\cO_Y$, and $\tilde a_i\in \cO_{\tilde Y}$ are arbitrary liftings of $a_i$.
The map $\vartheta_Y$ appearing above is the $i=0$ case of the canonical 
isomorphism defined in \cite[III. 1.5]{Illusie-Raynaud}
\begin{align}\label{theta_Y}
\vartheta_Y:W_n\Omega_Y^i&\xra{\simeq}
\cH^i(\Omega_{\tilde Y/W_nk}^\bullet).
\end{align}
Note that the map $\varrho_Y^*:W_n\cO_Y\ra \cO_{\tilde Y}$ is a morphism of sheaves of rings, and it induces a finite morphism $\varrho_Y:W_nY\ra \tilde Y$ (cf. \cite[I, paragraph after (2.4)]{Ekedahl-MultPropertiesOfDeRhamWittI}).
Altogether we have the following commutative diagram of schemes (cf. \cite[I. (2.4)]{Ekedahl-MultPropertiesOfDeRhamWittI})
$$
\xymatrix{
&
\tilde X\ar[dd]^{\tilde f}\ar[rrr]^(.4){\varrho_X}&&&
W_nX\ar[dd]^{W_nf}\\
X\ar[dd]_(.3){f}\ar@{^(->}[ru]^{\tilde i_X}&&&
X\ar[dd]_(.3){f}\ar@{^(->}[ru]^{i_X}&\\
&
\tilde Y\ar[dd]^{\pi_{\tilde Y}}\ar[rrr]^(.4){\varrho_Y}&&&
W_nY\ar[dd]^{W_n\pi_Y}\\
Y\ar[dd]_(.3){\pi_Y}\ar@{^(->}[ru]^{\tilde i_Y}&&&
Y\ar[dd]_(.3){\pi_Y}\ar@{^(->}[ru]^{i_Y}&\\
&
W_nk\ar[rrr]^(.3){\simeq}_(.3){W_n(F_k^n)}&&&W_nk.\\
k\ar@{^(->}[ru]\ar[rrr]^{\simeq}_{F_k^n}&&&k\ar@{^(->}[ru]&
}
$$

\ble\label{C.3}
Set $K_{\tilde X}=\pi_{\tilde X}^{\triangle}W_nk$, and $K_{\tilde Y}= \pi_{\tilde Y}^\triangle W_nk$.
The $(-d)$-th cohomology of the map $\Tr_{\tilde f}: \tilde f_*K_{\tilde X}\ra K_{\tilde Y}$ gives a map $\tilde f_*\Omega_{\tilde X/W_nk}^d\ra \Omega_{\tilde Y/W_nk}^d$, which we again denote by $\Tr_{\tilde f}$. 
Then by passing to quotients, this map $\Tr_{\tilde f}$ induces a well-defined map
$$\tau_{\tilde f}:\cH^d(\tilde f_*\Omega_{\tilde X/W_nk}^\bullet)\ra \cH^d(\Omega_{\tilde Y/W_nk}^\bullet).$$
Moreover, the map $\tau_{\tilde f}$ is compatible with $\Tr_{W_nf}$ defined in (\ref{Trace map on top Witt differentials: where it is induced from}), i.e., the following diagram commutes:
$$\xymatrix{
(W_nf)_*W_n\Omega_X^d
\ar[r]^(.6){\Tr_{W_nf}}
\ar[d]_{(W_nf)_*\vartheta_X}^{\simeq}&
W_n\Omega_Y^d
\ar[d]_(.4){\vartheta_Y}^(.4){\simeq}
\\
(\varrho_Y\tilde f)_*\cH^d(\Omega_{\tilde X/W_nk}^\bullet)
\ar[r]^{(\varrho_Y)_*\tau_{\tilde f}}&
(\varrho_Y)_*\cH^d(\Omega_{\tilde Y/W_nk}^\bullet)
}$$

\ele

\bpf
We argue the other way around, namely we define the map $\tau_{\tilde f}:
\cH^{d}(\tilde f_*\Omega_{\tilde X/W_nk}^\bullet)\ra \cH^d(\tilde \Omega_{\tilde Y/W_nk}^\bullet)$ via 
$\Tr_{W_nf}:(W_nf)_*W_n\Omega^d_X\ra W_n\Omega_Y^d$, 
and then show that this is the reduction of $\Tr_{\tilde f}: \tilde f_*\Omega_{\tilde X/W_nk}^d\ra \Omega_{\tilde Y/W_nk}^d$. 

First of all, via isomorphisms $\vartheta_X$, $\vartheta_Y$, the map $\Tr_{W_nf}:(W_nf)_*W_n\Omega_X^d\ra W_n\Omega_Y^d$ defined in (\ref{Trace map on top Witt differentials: where it is induced from}) induces a well-defined map $\tau_{\tilde f}:\cH^d(\tilde f_*\Omega_{\tilde X/W_nk}^\bullet)\ra \cH^d(\Omega_{\tilde Y/W_nk}^\bullet)$. 
To show the compatibility with $\Tr_{\tilde f}$, one needs the observation of Ekedahl
that the composite map
\begin{align*}
t_Y:
(\varrho_Y)_*\Omega_{\tilde Y/W_nk}^d[d]
\xra{\simeq}
(\varrho_Y)_*K_{\tilde Y}
&\simeq
(\varrho_Y)_*\pi_{\tilde Y}^\triangle W_nk
\xra[\simeq]{(\varrho_Y)_*\pi_{\tilde Y}^\triangle (\ref{(1.0.1)})^n}
\\
&(\varrho_Y)_*\pi_{\tilde Y}^\triangle (W_nF_k^n)^\triangle W_nk
\simeq (\varrho_Y)_*(\varrho_Y)^\triangle (W_n\pi_{Y})^\triangle W_nk
\xra{\Tr_{\varrho_Y}}
K_{n,Y}
\nonumber
\end{align*}
factors through $\bar t_Y: (\varrho_{Y})_*\cH^d(\Omega_{\tilde Y/W_nk}^\bullet)[d]\ra K_{n,Y}$ (cf. \cite[\S1 (2.6)]{Ekedahl-MultPropertiesOfDeRhamWittI}). Then he defined the map $W_n\Omega_Y^d[d]\ra K_{n,Y}$ to be the composite
\beq\label{WnOmega^d -> KnY}
s_Y:
W_n\Omega_Y^d[d]\xra[\simeq]{\vartheta_Y}\cH^d(\Omega_{\tilde Y/W_nk}^\bullet)[d]\xra{\bar t_Y} K_{n,Y}.
\eeq

Now consider the following diagram of complexes of sheaves
$$\resizebox{\displaywidth}{!}{
\xymatrix@C=5ex{
&(W_nf)_*K_{n,X}\ar[rrr]^(.4){\Tr_{W_nf}}
&&&
K_{n,Y}
\\
(W_nf)_*W_n\Omega_X^d[d]
\ar[ur]^{(W_nf)_*s_X}\ar[rrr]^(.6){\Tr_{W_nf}}
\ar[ddd]_(.4){(W_nf)_*\vartheta_X}^(.4){\simeq}&&&
W_n\Omega_Y^d[d]\ar[ur]^{s_Y}\ar[ddd]_(.4){\vartheta_Y}^(.4){\simeq}&
\\
\\
&(\varrho_Y\tilde f)_*\Omega_{\tilde X/W_nk}^d[d]
\ar[dl]\ar[uuu]_(.4){(W_nf)_* t_X }
\ar[rrr]^(.4){(\varrho_Y)_*\Tr_{\tilde f}}
&&&
(\varrho_Y)_*\Omega_{\tilde Y/W_nk}^d[d]
\ar[dl]\ar[uuu]_(.4){t_Y}.
\\
(\varrho_Y\tilde f)_*\cH^d(\Omega_{\tilde X/W_nk}^\bullet)[d]
\ar[rrr]^{(\varrho_Y)_*\tau_{\tilde f}}\ar[uuuur]^{(W_nf)_*\bar t_X} &&&
(\varrho_Y)_*\cH^d(\Omega_{\tilde Y/W_nk}^\bullet)[d]
\ar[uuuur]^{\bar t_Y}&
}}
$$
The unlabeled arrows are given by the natural quotient maps.
The front commutes by the definition of $\tau_{\tilde f}$. The top commutes by the definition of $\Tr_{W_nf}: (W_nf)_*W_n\Omega_X^d\ra W_n\Omega_Y^d$. The triangles in the right (resp. the left) side commute due to the definition of $\bar t_Y$ and $s_Y$ (resp. $\bar t_X$ and $s_X$).
The back square commutes, because the trace map $\Tr_{\tilde f}$ is functorial with respect to maps of residual complexes with the same associated filtration by \Cref{trace for residual}(3).
We want to show that the bottom square commutes. 
To this end, it suffices to show 
$(\varrho_Y)_*\Tr_{\tilde f}:(\varrho_Y\tilde f)_*\Omega_{\tilde X/W_nk}^d \ra (\varrho_Y)_*\Omega_{\tilde Y/W_nk}^d$
is compatible with 
$\Tr_{W_nf}: (W_nf)_*W_n\Omega_X^d\ra W_n\Omega_Y^d$ 
via $\vartheta_X$ and $\vartheta_Y$. 
Because the map 
$\Tr_{W_nf}: (W_nf)_*W_n\Omega_X^d\ra W_n\Omega_Y^d$
is determined by the degree $-d$ part of the map 
$\Tr_{W_nf}: (W_nf)_*K_{n,X}\ra K_{n,Y}$, 
we are reduced to show compatibility of
$(\varrho_Y)_*\Tr_{\tilde f}:(\varrho_Y\tilde f)_*\Omega_{\tilde X/W_nk}^d \ra (\varrho_Y)_*\Omega_{\tilde Y/W_nk}^d$ 
with
$\Tr_{W_nf}: (W_nf)_*K_{n,X}\ra K_{n,Y}$ via $(W_nf)_*(s_X\circ \vartheta_X^{-1})$ and $s_Y\circ\vartheta_Y^{-1}$. 
By the commutativity of the left and right squares, this is reduced to the commutativity of the square on the back, which is known. Therefore the bottom square commutes as a result.
\epf

The notation $\tau_{\tilde f}$ is only temporarily used in the lemma above. Later we will denote $\tau_{\tilde f}$ by $\Tr_{\tilde f}$.

\ble[{\cite[8.4(ii)]{BER12}}]\label{vartheta lemma}
Let $W_n\Omega_{\tilde Y/W_nk}^\bullet$ denote the relative de Rham-Witt complex defined by \cite{LangerZink}.
The rest of the notations are the same as above.
There is a commutative diagram
$$\xymatrix{
W_{n+1}\Omega_{\tilde Y/W_nk}^q\ar[r]^{F^n}\ar[d]^{}&
Z^q(\Omega_{\tilde Y/W_nk}^\bullet)\ar[d]
\\
W_n\Omega_Y^q\ar[r]^{\vartheta_Y}&\cH^q(\Omega_{\tilde Y/W_nk}^\bullet).
}$$

\ele

Recall that any element in $W_n\Omega_{k[X_1,...,X_d]}^d$ is uniquely written as a sum of (\ref{FormNonint}) and (\ref{FormInt}).
\ble\label{Trace map on top differentials: Example}
Let 
$$
C'=C_n': W_n\Omega_{k[X_1,...,X_d]}^d
\ra
W_n\Omega_{k[X_1,...,X_d]}^d.
$$ 
be the map given by the $-d$-th cohomology of the level $n$ Cartier operator for residual complexes (cf. (\ref{C'_n def})).
Let $\alpha=\sum_{j=0}^{n+v_1-1}V^j[\alpha_j]\in W_{n+v_1}k$ with each $\alpha_j\in k$.
Let $\beta=\sum_{j=0}^{n-1} V^j[\beta_j]\in W_nk$ with each $\beta_j\in k$.
\benu 
\item 
If $v_1=1-n$, 
$$C'\Big(dV^{-v_1}(\alpha[X_1]^{h_1'})
\cdots 
dV^{-v_r}([X_r]^{h_r'})\cdot 
F^{v_{r+1}}d[X_{r+1}]^{h_{r+1}'}
 \cdots 
F^{v_d}d[X_d]^{h_d'}\Big)=0.$$

\item 
If $1-n<v_1<0$, $v_{r+1}=...=v_{r+s}=0$ ($s$ can be zero),
\begin{align*} 
&C'\Big(dV^{-v_1}(\alpha[X_1]^{h_1'})
\cdots 
dV^{-v_r}([X_r]^{h_r'})\cdot 
F^{v_{r+1}}d[X_{r+1}]^{h_{r+1}'}
 \cdots 
F^{v_d}d[X_d]^{h_d'}\Big)
\\ 
&=dV^{1-v_1}(R(\tilde \alpha) [X_1]^{h_1'})\cdots 
dV^{1-v_r}[X_r]^{h_r'}
\cdot 
dV[X_{r+1}]^{h_{r+1}'}
\cdots 
dV[X_{r+s}]^{h_{r+s}'}
\\
&\omit\hfill$
\cdot
F^{v_{r+s+1}-1}d[X_{r+s+1}]^{h_{r+s+1}'}
\cdots 
F^{v_{d}-1}d[X_{d}]^{h_{d}'}
$.
\end{align*}
Here
$$\tilde\alpha=\sum_{j=0}^{n+v_1-1}V^j[\tilde \alpha_j]_{n+v_1+1-j}\in W_{n+v_1+1}(W_nk),$$
where $\tilde \alpha_j:=[\alpha_j]_n\in W_nk$ is the Teichm\"uller lift of $\alpha_j\in k$.
\item 
If 
$v_1\ge 0$, $v_{1}=...=v_{s}=0$ ($s$ can be zero),
\begin{align*} 
&C'\Big( 
\beta 
F^{v_1}d[X_1]^{h_1'} 
 \cdots 
F^{v_d}d[X_d]^{h_d'}
\Big)
\\
&=
(W_nF_k)^{-1}(\beta)\cdot
dV[X_{1}]^{h_{1}'}\cdots 
dV[X_{s}]^{h_{s}'}\cdot
F^{v_{s+1}-1}d[X_{s+1}]^{h_{s+1}'}
 \cdots 
F^{v_d-1}d[X_d]^{h_d'}.
\end{align*} 
Here 
$$\tilde \beta=\sum_{j=0}^{n-1} V^j[\tilde \beta_j]_{n+1-j}\in W_{n+1}(W_nk),$$
where $\tilde \beta_j:=[\beta_j]_n\in W_nk$ is the Teichm\"uller lifts of $\beta_j\in k$. 
\eenu 

\ele

\bpf
Consider the map $W_nF_X:W_nX\ra W_nX$ with $X:=\A^d_k$. It is not a map of $W_nk$-schemes a priori, but after labeling the source by $W_nX:=W_n\A^d_{k_1}$ and the target by $W_nY:=W_n\A^d_{k_2}$, one can realize $W_nF_X$ as a map of $W_nk_2$-schemes (the $W_nk_2$-scheme structure of $W_nX$ is given by $W_nF_X\circ W_n\pi_Y$, where $\pi_Y:Y\ra k_2$ denotes the structure morphism of the scheme $Y$).
Write
$$\tilde X=\A^d_{W_nk_1}=\Spec W_nk_1[X_1,\dots, X_d] \quad (\text{resp. $\tilde Y=\A^d_{W_nk_2}=\Spec W_nk_2[X_1,\dots, X_d]$}),$$
and take the canonical lift $\tilde F_{\tilde X}$ of $F_X$ as in \Cref{Trace map on top differentials: Example on lifting}. Consider
\begin{equation}
\label{compute Tr_Wnf via Tr_tilde f}
\resizebox{\displaywidth}{!}{
\xymatrix@C=8em{
(W_nF_X)_*W_n\Omega_{X/k_1}^d
\ar[r]^{(W_nF_X)_*(\ref{(1.0.2)})}_{\simeq}
\ar[d]^{(W_nF_X)_*\vartheta_X}&
(W_nF_X)_*W_n\Omega_{X/k_2}^d
\ar[r]^{\Tr_{W_nF_X}}\ar[d]^{(W_nF_X)_*\vartheta_X}&
W_n\Omega_{Y/k_2}^d\ar[d]^{\vartheta_Y}
\\
(\varrho_Y\tilde F_{\tilde X})_* \cH^{d}(\Omega_{\tilde X/W_nk_1}^\bullet)
\ar[r]^{(\varrho_Y\tilde F_{\tilde X})_*\pi_{\tilde X}^\triangle (\ref{(1.0.1)})}_{\simeq}
&(\varrho_Y\tilde F_{\tilde X})_* \cH^{d}(\Omega_{\tilde X/W_nk_2}^\bullet)
\ar[r]^{(\varrho_Y)_*( \Tr_{\tilde F_{\tilde X}})}&
(\varrho_Y)_* \cH^{d}(\Omega_{\tilde Y/W_nk_2}^\bullet).
}}
\end{equation}
The composite map of the top row is $C'$ (cf. (\ref{C'_n def}) and Ekedahl's quasi-isomorphism \Cref{Ekedahl qis}). The composite of the bottom row is induced from $\varrho_{Y,*}$(\ref{eq:trace of lift Frob on Coh}).
The right side commutes due to \Cref{C.3}. The left side commutes by the naturality.
Hence we can decompose $C'$ in the following way:
$$
C'=\vartheta_Y^{-1}\circ 
 (\ref{eq:trace of lift Frob on Coh}) \circ \vartheta_X: W_n\Omega_{X/k_1}^d
\ra
W_n\Omega_{Y/k_2}^d.
$$

Consider the first two cases. Suppose $v_1<0$ and suppose there are $s$ many $v_j$'s being zero,
$$v_{1}\le\cdots \le v_r<0=\cdots =0<v_{r+s+1}\le\cdots \le  v_d.$$
Note that in $W_{n+1}\Omega_{(W_nk)[X_1,...,X_d]/W_nk}^1$:
\begin{align} 
\label{3eq}
F^{n+v_1}d(\tilde\alpha[X_1]^{h_1'})
&=F^{n+v_1}d\Big([X_1]^{h_1'}\cdot\sum_{j=0}^{n+v_1-1}V^j[\tilde \alpha_j ]\Big)
\\
&=F^{n+v_1}d
\Big(
\sum_{j=0}^{n+v_1-1}V^j[\tilde \alpha_j X_1^{h_1'p^j}]
\Big)
\nonumber\\
&=
\sum_{j=0}^{n+v_1-1}
F^{n+v_1-j}d[\tilde \alpha_j X_1^{h_1'p^j}]
\nonumber\\
&=\sum_{j=0}^{n+v_1-1}(\tilde \alpha_j X_1^{h_1'p^j})^{p^{n+v_1-j}-1}d(\tilde \alpha_j X_1^{h_1'p^j})
\nonumber\\
&=\sum_{j=0}^{n+v_1-1}h_1'p^j\cdot\tilde \alpha_j^{p^{n+v_1-j}}\cdot  X_1^{h_1'p^{n+v_1}-1} dX_1
\nonumber
\end{align}
Similarly we have
\beq\label{5eq}\sum_{j=0}^{n+v_1-1} 
p^jh_1'\cdots h_d' \cdot (W_nF_k)^{-1}(\tilde\alpha_j^{p^{n+v_1-j}})
X_1^{h_1'p^{n+v_1-1}-1} dX_1
=F^ndV^{1-v_1}(R(\tilde \alpha) [X_1]^{h_1}).
\eeq 
Here $R$ is the restriction map  $R:W_{n+v_1+1}(W_nk)\ra W_{n+v_1}(W_nk)$.
Now according to the formula (\ref{Trace map on top differentials: in particular for Ad}) and \Cref{vartheta lemma}, we carry out the following calculations.
\begin{align}
\label{PartNonInt}
&C'\Big(dV^{-v_1}(\alpha[X_1]^{h_1'})
\cdots 
dV^{-v_r}([X_r]^{h_r'})\cdot 
F^{v_{r+1}}d[X_{r+1}]^{h_{r+1}'}
 \cdots 
F^{v_d}d[X_d]^{h_d'}\Big)
\\
&=\vartheta_Y^{-1}\circ 
(\ref{eq:trace of lift Frob on Coh}) 
\circ F^n
\Big(dV^{-v_1}(\tilde\alpha[X_1]^{h_1'})
\cdots 
dV^{-v_r}([X_r]^{h_r'})\cdot 
F^{v_{r+1}}d[X_{r+1}]^{h_{r+1}'}
 \cdots 
F^{v_d}d[X_d]^{h_d'}
\Big)
\nonumber\\
&=\vartheta_Y^{-1}\circ 
(\ref{eq:trace of lift Frob on Coh}) 
\Big(
F^{n+v_1}d(\tilde\alpha[X_1]^{h_1'})
\cdots 
F^{n+v_r}d([X_r]^{h_r'})\cdot 
F^{n+v_{r+1}}d[X_{r+1}]^{h_{r+1}'}
 \cdots 
F^{n+v_d}d[X_d]^{h_d'}
\Big)
\nonumber\\
&=\vartheta_Y^{-1}\circ 
(\ref{eq:trace of lift Frob on Coh}) 
\Big(
\Big(\sum_{j=0}^{n+v_1-1}h_1'p^j\cdot\tilde \alpha_j^{p^{n+v_1-j}}\cdot  X_1^{h_1'p^{n+v_1}-1} dX_1\Big)
\cdot 
\Big(h_2'\cdot X_2^{h_2'p^{n+v_2}-1} dX_2\Big)\cdot
\nonumber\\
&\omit\hfill$\cdots 
\Big(h_d'\cdot X_d^{h_d'p^{n+v_d}-1} dX_d\Big)
\Big)$
\quad 
(\text{by (\ref{3eq})})
\nonumber\\
&=\vartheta_Y^{-1}\Big(
\sum_{j=0}^{n+v_1-1} 
p^jh_1'\cdots h_d' \cdot (W_nF_k)^{-1}(\tilde\alpha_j^{p^{n+v_1-j}})
X_1^{h_1'p^{n+v_1-1}-1} dX_1\cdots
X_d^{h_d'p^{n+v_d-1}-1} dX_d
\Big)
\nonumber\\
&=\vartheta_Y^{-1}\Big(
F^ndV^{1-v_1}(R(\tilde \alpha) [X_1]^{h_1'})\cdots 
F^ndV^{1-v_r}[X_r]^{h_r'}
\cdot 
F^{n}dV[X_{r+1}]^{h_{r+1}'}
\cdots 
F^{n}dV[X_{r+s}]^{h_{r+s}'}
\nonumber\\
&\omit\hfill$
\cdot
F^{n+v_{r+s+1}-1}d[X_{r+s+1}]^{h_{r+s+1}'}
\cdots 
F^{n+v_{d}-1}d[X_{d}]^{h_{d}'}
$\Big)
\quad 
(\text{by (\ref{5eq})})
\nonumber 
\end{align}
If $v_1=1-n$, $(\ref{PartNonInt})=0$ because 
$$F^ndV^{1-v_1}(R(\tilde \alpha) [X_1]^{h_1'})=d(R(\tilde \alpha) [X_1]^{h_1'})=0$$
in $\cH^1(\Omega^\bullet_{(W_nk)[X_1,...,X_d]/W_nk})$.
If $v_1\neq 1-n$ (hence $v_j>1-n$ for all $j$),
\begin{align*}
(\ref{PartNonInt})
&=
dV^{1-v_1}(R(\tilde \alpha) [X_1]^{h_1'})\cdots 
dV^{1-v_r}[X_r]^{h_r'}
\\
&\omit\hfill$\cdot 
dV[X_{r+1}]^{h_{r+1}'}
\cdots 
dV[X_{r+s}]^{h_{r+s}'}
\cdot
F^{v_{r+s+1}-1}d[X_{r+s+1}]^{h_{r+s+1}'}
\cdots 
F^{v_{d}-1}d[X_{d}]^{h_{d}'}
$
\end{align*}
and this is the same as what our lemma claims.

Now we check the  third case. 
If all the $v_j\ge 0$, suppose the first $s$ $v_j$'s are zero, 
$$0=\cdots=0<v_{s+1}\le \cdots\le v_d.$$
Note that in $W_nk$,
\begin{align}
\label{coeffbeta}
&(W_nF_k)^{-1}(F^n(\tilde \beta))
=(W_nF_k)^{-1}(F^n(\sum_{j=0}^{n-1} V^j[\tilde \beta_j]_{n+1-j}))
=(W_nF_k)^{-1}(\sum_{j=0}^{n-1}F^{n-j} [\tilde \beta_j]_{n+1-j}))
\\
&=(W_nF_k)^{-1}(\sum_{j=0}^{n-1}\tilde \beta_j^{p^{n-j}})
=\sum_{j=0}^{n-1}\tilde \beta_j^{p^{n-j-1}}
=F^n\Big(\sum_{j=0}^{n-1}
V^j[(W_nF_k)^{-1}\tilde \beta_j]_{n+1-j}
\Big)
=F^n\Big((W_{n+1}(W_nF_k))^{-1}(\beta)\Big).
\nonumber
\end{align}
We carry out the computation
\begin{align*}
&C'\Big( 
\beta 
F^{v_1}d[X_1]^{h_1'} 
 \cdots 
F^{v_d}d[X_d]^{h_d'}
\Big)
\\
&=\vartheta_Y^{-1}\circ 
(\ref{eq:trace of lift Frob on Coh}) \circ F^n
\Big(\tilde \beta\cdot
F^{v_1}d[X_1]^{h_1'} 
 \cdots 
F^{v_d}d[X_d]^{h_d'}
\Big)
\\
&=\vartheta_Y^{-1}\circ 
(\ref{eq:trace of lift Frob on Coh}) 
\Big(F^n(\tilde \beta)\cdot
F^{n+v_1}d[X_1]^{h_1'} 
 \cdots 
F^{n+v_d}d[X_d]^{h_d'}
\Big)
\\
&=\vartheta_Y^{-1}\circ 
(\ref{eq:trace of lift Frob on Coh}) 
\Big(F^n(\tilde \beta)\cdot 
h_1'\cdots h_d'
X_1^{h_1'p^{n+v_1}-1}dX_1 
 \cdots 
F^{n+v_d}d[X_d]^{h_d'}
\Big)
\\
&=\vartheta_Y^{-1}
\Big((W_nF_k)^{-1}(F^n(\tilde \beta))\cdot 
h_1'
X_1^{h_1'p^{n+v_1-1}-1}dX_1 
 \cdots 
 h_d'
X_d^{h_d'p^{n+v_d-1}-1}dX_d 
\Big)
\\
&=\vartheta_Y^{-1}
\Big(
F^n\Big((W_{n+1}(W_nF_k))^{-1}(\beta)\Big)\cdot
F^ndV[X_{1}]^{h_{1}'}\cdots 
F^ndV[X_{s}]^{h_{s}'}\cdot
\\
&\omit\hfill
$F^{n+v_{s+1}-1}d[X_{s+1}]^{h_{s+1}'}
 \cdots 
F^{n+v_d-1}d[X_d]^{h_d'}$
\Big)
\quad\text(by (\ref{coeffbeta}))
\\
&=(W_nF_k)^{-1}(\beta)\cdot
dV[X_{1}]^{h_{1}'}\cdots 
dV[X_{s}]^{h_{s}'}\cdot
\\
&\omit\hfill
$F^{v_{s+1}-1}d[X_{s+1}]^{h_{s+1}'}
 \cdots 
F^{v_d-1}d[X_d]^{h_d'}$
\end{align*}
In the last equality we have used that 
$$\vartheta_Y^{-1}\Big(
F^n\Big((W_{n+1}(W_nF_k))^{-1}(\beta)\Big)
\Big)=(W_nF_k)^{-1}(\beta).$$
We hence proved the lemma.

\epf

\subsubsection{Criterion for surjectivity of $C'-1$}
The following proposition is proven in the smooth case by Illusie-Raynaud-Suwa \cite[2.1]{Suwa-NoteOnGerstenForlogSheaves}. The proof presented here is due to R\"ulling.
\bprop[Raynaud-Illusie-Suwa]\label{Raynaud-Illusie}
Let $k=\bar k$ be an algebraically closed field of characteristic $p>0$ and let $X$ be a separated scheme of finite type over $k$. 
Then for every $i$, 
$C'-1$ induces a surjective map on global cohomology groups
$$H^i(W_nX,K_{n,X}):=R^i\Gamma(W_nX,K_{n,X})\xra{C'-1} H^i(W_nX,K_{n,X}).$$
In particular,
$$R^i\Gamma(W_nX,K_{n,X,log})\simeq 
H^i(W_nX,K_{n,X})^{C'-1}.$$
\eprop
\bpf
Take a Nagata compactification of $X$, i.e., an open immersion
$$j:X\hra \bar X$$
such that $\bar X$ is proper over $k$.
The boundary $\bar X\setminus X$ is a closed subscheme in $\bar X$. By blowing up in $\bar X$ one can assume $\bar X\setminus X$ is the closed subscheme associated to an effective Cartier divisor $D$ on $\bar X$. We can thus assume $j$ is an affine morphism. 
Therefore 
$$W_nj:W_nX\hra W_n\bar X$$ 
is also an affine morphism. 

For any quasi-coherent sheaf $\cM$ on $W_n\bar X$, the difference between $\cM$ and $(W_nj)_*(W_nj)^*\cM$ are precisely those sections that have poles (of any order) at $\Supp D=W_n\bar X\setminus W_nX$.
Suppose that the effective Cartier divisor $D$ is represented by $(U_i,f_i)_i$, where $\{U_i\}_i$ is an affine cover of $\bar X$, and $f_i\in \Gamma(U_i,\cO_X)$. 
Recall that $\cO_{\bar X}(mD)$ denotes the line bundle on $\bar X$ which is the inverse (as line bundles) of the $m$-th power of the ideal sheaf of $\bar X\setminus X\hra \bar X$. Locally, one has an isomorphism
$$\cO_{\bar X}(mD)\mid_{U_i} \simeq \cO_{U_i}\cdot \frac{1}{f_i^m}$$ 
for each $i$.
Denote by $W_n\cO_{\bar X}(mD)$ the line bundle on $W_n\bar X$ such that 
$$W_n\cO_{\bar X}(mD)\mid_{U_i} \simeq W_n\cO_{U_i}\cdot \frac{1}{[f_i]^m},$$
where $[-]=[-]_n$ denotes the Teichm\"uller lift. Denote 
$$\cM(mD):=\cM\otimes_{W_n\cO_{\bar X}} W_n\cO_{\bar X}(mD).$$ 
The natural map
\begin{equation}\label{eq:ind lim M}
\cM(*D):=\colim_m \cM(mD)\xra{\simeq} (W_nj)_*(W_nj)^*(\cM(mD))= (W_nj)_*(W_nj)^*\cM
\end{equation}
is an isomorphism of sheaves. 
Here the inductive system on the left hand side is given by the natural map $$\cM(mD):=\cM\otimes_{W_n\cO_{\bar X}} W_n\cO_{\bar X}(mD)\ra \cM\otimes_{W_n\cO_{\bar X}} W_n\cO_{\bar X}((m+1)D)$$
induced from the inclusion $W_n\cO_{\bar X}(mD)\hra W_n\cO_{\bar X}((m+1)D)$, i.e., locally on $U_i$, this inclusion is the map 
\begin{align*}
W_n\cO_{\bar X}(mD)\mid_{U_i}&
\hra W_n\cO_{\bar X}((m+1)D)\mid_{U_i}\\
\frac{a}{[f_i]^m}&\mapsto \frac{a[f_i]}{[f_i]^{m+1}}.
\end{align*}
where $a\in W_n\cO_{U_i}$. As a result,
\begin{align}
\label{eq:ind lim M coh}
H^i(W_nX,(W_nj)^*\cM)&=H^i(R\Gamma(W_n\bar X, R(W_nj)_*(W_nj)^*\cM))
\\
&=H^i(R\Gamma(W_n\bar X, (W_nj)_*(W_nj)^*\cM))
\qquad\quad\text{($W_nj$ is affine)}
\nonumber
\\
&=H^i(R\Gamma(W_n\bar X,\colim_m \cM(mD))\qquad\qquad\text{(\ref{eq:ind lim M})}
\nonumber
\\
&=\colim_m H^i(W_n\bar X,\cM(mD)).
\nonumber
\end{align}

Apply this to the bounded complex $K_{n,\bar X}$ of injective quasi-coherent $W_n\cO_{\bar X}$-modules.  
Taking into account $K_{n,X}\simeq (W_nj)^* K_{n,\bar X}$ by \Cref{f triangle}(2),
(\ref{eq:ind lim M}) gives an isomorphism of complexes
\begin{equation}\label{eq:ind lim KX}
K_{n,\bar X}(*D):=\colim_m K_{n,\bar X}(mD)\xra{\simeq} (W_nj)_*K_{n,X},
\end{equation}
and (\ref{eq:ind lim M coh}) gives an isomorphism of $W_nk$-modules
$$\colim_m H^i(W_n\bar{X}, K_{n,\bar X}(mD))=H^i(W_nX,K_{n,X}).$$
Via the projection formula \cite[II.5.6]{Hartshorne-RD}
and tensoring
$$C':(W_nF_X)_*K_{n,\bar X}\ra K_{n,\bar X}$$
with $W_n\cO_{\bar X}(mD)$, one gets a map
\begin{align*}
(W_nF_X)_*(K_{n,\bar X}&\otimes_{W_n\cO_{\bar X}} W_n\cO_{\bar X}(pmD))
\simeq (W_nF_X)_*(K_{n,\bar X}\otimes_{W_n\cO_{\bar X}} (W_nF_X)^*W_n\cO_{\bar X}(mD))\\
&\simeq
((W_nF_X)_*K_{n,\bar X})\otimes_{W_n\cO_{\bar X}} W_n\cO_{\bar X}(mD)
\xra{C'\otimes \id_{W_n\cO_{\bar X}(mD)}} 
K_{n,\bar X}\otimes_{W_n\cO_{\bar X}} W_n\cO_{\bar X}(mD).
\end{align*}
Precomposing with the natural map
$$(W_nF_X)_*(K_{n,\bar X}\otimes_{W_n\cO_{\bar X}} W_n\cO_{\bar X}(mD))\ra (W_nF_X)_*(K_{n,\bar X}\otimes_{W_n\cO_{\bar X}} W_n\cO_{\bar X}(pmD)),$$
and taking the global section cohomologies, one gets
$$
C':H^i(W_n\bar{X}, K_{n,\bar X}(mD))\ra H^i(W_n\bar{X}, K_{n,\bar X}(mD)).
$$
To show the surjectivity of
$$C'-1:H^i(W_nX,K_{n,X})\xra{} H^i(W_nX,K_{n,X}),$$
it suffices to show the surjectivity for
$$C'-1:H^i(W_n\bar{X}, K_{n,\bar X}(mD))\ra H^i(W_n\bar{X}, K_{n,\bar X}(mD)).$$
Because $\cH^q(K_{n,\bar X})$ are coherent sheaves on the proper scheme $\bar X$ for all $q$,  $\cH^q(K_{n,\bar X}\otimes_{W_n\cO_{\bar X}} W_n\cO_{\bar X}(mD))=\cH^q(K_{n,\bar X})\otimes_{W_n\cO_{\bar X}} W_n\cO_{\bar X}(mD)$
are also coherent, therefore the local-to-global spectral sequence implies that 
$$M:=H^i(W_n\bar{X}, K_{n,\bar X}(mD))$$
is a finite $W_nk$-module. 
Now $M$ is equipped with an endomorphism $C'$ which acts $p^{-1}$-linearly (cf. \Cref{def of sigma-linear generalized}). The proposition is then a direct consequence of \Cref{1-T surjective prop in appendix for W_n}.
\epf

The following proposition is a corollary of \cite[Lemma 2.1]{Suwa-NoteOnGerstenForlogSheaves}. We restate it here as a convenient reference.
\bprop[Raynaud-Illusie-Suwa]\label{surjectivities 1-C C^-1 -1 etc}
Assume $k=\bar k$.
If $X$ is separated smooth over $k$ of pure dimension $d$,
$$C-1:W_n\Omega^d_X\ra W_n\Omega^d_X$$
is surjective.
\eprop
\bpf
Apply affine locally the $H^{-d}$-case of \Cref{Raynaud-Illusie}. Then Ekedahl's quasi-isomorphism $W_n\Omega^d_X[d]\simeq K_{n,X}$ from \Cref{Ekedahl qis}  together with compatibility of $C'$ and $C$ from \Cref{Compatibility of $C'_n$ and $C$: Prop} gives the claim.
\epf

\brmk\label{surjectivities 1-C C^-1 -1 etc CM}
If $X$ is Cohen-Macaulay of pure dimension $d$, $W_nX$ is also Cohen-Macaulay by Serre's $S_k$-criterion (\cite[(5.7.3)(i)]{EGAIV-2}) of the same pure dimension, and thus the complex $K_{n,X}$ is concentrated at degree $-d$ for all $n$ \cite[3.5.1]{Conrad-GDBC}. Denote by $W_n\omega_X$ the only nonzero cohomology sheaf of $K_{n,X}$ in this case. Then the same reasoning as in \Cref{surjectivities 1-C C^-1 -1 etc} shows that if $k=\bar k$ and $X$ is Cohen-Macaulay over $k$ of pure dimension, the map
$$C'-1:W_n\omega_X\ra W_n\omega_X$$
is surjective.
\ermk

\subsubsection{Comparison between $W_n\Omega^d_{X,log}$ and $K_{n,X,log}$}\label{subsubsection Comparison between W_nOmega^d log and KnXlog}
Let $X$ be a $k$-scheme.
Denote by $d\log$ the following map of abelian \'etale sheaves
\begin{align*}
d\log: (\cO_{X,\et}^{*})^{\otimes q}&\ra W_n\Omega_{X,\et}^q,\\
a_1\otimes\dots\otimes a_q&\mapsto d\log[a_1]_n\dots d\log[a_q]_n,
\end{align*}
where $a_1,\dots, a_q\in \cO_{X,\et}^*$, $[-]_n:\cO_{X,\et}\ra W_n\cO_{X,\et}$ denotes the Teichm\"{u}ller lift, and $d\log[a_{i}]_n:=\frac{d[a_{i}]_n}{[a_{i}]_n}$. We will denote its sheaf theoretic image by
$W_n\Omega^q_{X,log,\et}$ and call it the \'etale sheaf of log forms. We denote by $W_n\Omega^q_{X,log}:=W_n\Omega^q_{X,log,\Zar}:=\epsilon_*W_n\Omega^q_{X,log,\et}$, and call it the Zariski sheaf of log forms.

\ble[{\cite[lemme 2]{CTSS-TorDanGroupeDeChow}, \cite[1.6(ii)]{GrosSuwa-AbelJacobi}}]
\label{ses log forms are inv}
Let $X$ be a separated smooth $k$-scheme. Then we have the following left exact sequences
\beq\label{1-bar F ses}
0\ra W_n\Omega_{X,log}^q \ra W_n\Omega_{X}^q \xra{1-\bar F} W_n\Omega_{X}^q/dV^{n-1},
\eeq
\beq\label{C-1 ses}
0\ra W_n\Omega_{X,log}^q\ra W_n\Omega_X'^q \xra{C-1} W_n\Omega_X^q,
\eeq
where $W_n\Omega_X'^q:=F(W_{n+1}\Omega_X^q)$.
The right hand maps are also surjective if $t=\et$.
\ele

The following proposition collects what we have done so far.
\bprop[{cf. \cite[Prop. 4.2]{Kato-DualitypPrimaryEtale-II}}]\label{log qis complex log}
Let $X$ be a separated smooth scheme of pure dimension $d$ over a perfect field $k$. Then
\begin{enumerate}
\item
we have
$\cH^{-d}(K_{n,X,log})= W_n\Omega^d_{X,log}$, and $\cH^i(K_{n,X,log})=0$ for all $i\neq -d,-d+1$.
\item
If $k=\bar k$, the natural map
$$W_n\Omega^d_{X,log}[d]\xra{} K_{n,X,log}$$
is a quasi-isomorphism of complexes of abelian sheaves.
\end{enumerate}
\eprop

\bpf
\begin{enumerate}
\item
Since $C$ is compatible with $C'$ by \Cref{Compatibility of $C'_n$ and $C$: Prop},
the natural map
$\Cone(W_n\Omega^d_X[d]\xra{C-1}W_n\Omega_X^d[d])[-1]\ra K_{n,X,log}$ is a quasi-isomorphism by the five lemma and the Ekedahl quasi-isomorphism \Cref{Ekedahl qis}. The claim thus follows from the exact sequence (\ref{C-1 ses}).
\item
\Cref{surjectivities 1-C C^-1 -1 etc}+(1) above.
\end{enumerate}
\epf

\subsection{Localization triangle associated to $K_{n,X,log}$}\label{subsection localization triangle associated to KnXlog}
\subsubsection{Definition of $\Tr_{W_nf,log}$}
\bprop[Proper pushforward, {cf. \cite[(3.2.3)]{Kato-DualitypPrimaryEtale-II}}]
\label{def of Tr f log}
Let $f:X\ra Y$ be a proper map between separated schemes of finite type over $k$. Then so is $W_nf:W_nX\ra W_nY$, and we have a map $$\Tr_{W_nf,log}:(W_nf)_*K_{n,X,log}\ra K_{n,Y,log}$$ of complexes that fits into the following commutative diagram of complexes, where the two rows are distinguished triangles in $D^b(W_nX,\Z/p^n)$
\begin{displaymath}
\xymatrix{
(W_nf)_*K_{n,X,log}\ar[r]\ar[d]^{\Tr_{W_nf,log}}&(W_nf)_*K_{n,X}
\ar[r]^{C'-1}\ar[d]^{\Tr_{W_nf}}&(W_nf)_*K_{n,X}\ar[r]\ar[d]^{\Tr_{W_nf}}&
\\
K_{n,Y,log}\ar[r]&K_{n,Y}\ar[r]^{C'-1}&K_{n,Y}\ar[r]&.
}
\end{displaymath}
Moreover $\Tr_{W_nf,log}$ is compatible with compositions and open restrictions.
\eprop
This is the covariant functoriality of $K_{n,X,log}$ with respect to proper morphisms. Thus we also denote $\Tr_{W_nf,log}$ by $f_*$.
\bpf
It suffices to show the following diagrams commute.
$$\xymatrix@C=8em{
(W_nF_Y)_*(W_nf)_*K_{n,X}
\ar[r]_{\simeq}^(.44){(W_nF_Y)_*(W_nf)_*(\ref{(1.0.2)})}
\ar[d]^{(W_nF_Y)_*\Tr_{W_nf}}&
(W_nF_Y)_*(W_nf)_*(W_nF_X)^\triangle K_{n,X}
\ar[d]^{(W_nF_Y)_*\Tr_{W_nf}}
\\
(W_nF_Y)_*K_{n,Y}
\ar[r]_{\simeq}^(.4){(W_nF_Y)_*(\ref{(1.0.2)})}&
(W_nF_Y)_*(W_nF_Y)^\triangle K_{n,Y},
}$$
$$\resizebox{\displaywidth}{!}{\xymatrix@C=6em{
(W_nF_Y)_*(W_nf)_*(W_nF_X)^\triangle K_{n,X}
\ar[d]^{(W_nF_Y)_*\Tr_{W_nf}}
\ar[r]^{\simeq}  &
(W_nf)_*(W_nF_X)_*(W_nF_X)^\triangle K_{n,X}
\ar[r]^(.6){(W_nf)_*\Tr_{W_nF_X}}&
(W_nf)_*K_{n,X}\ar[d]^{\Tr_{W_nf}}
\\
(W_nF_Y)_*(W_nF_Y)^\triangle K_{n,Y}
\ar[rr]^(.6){\Tr_{W_nF_Y}}&&
K_{n,Y},
}}$$
where $\Tr_{W_nf}$ on the right of the first diagram and the left of the second diagram denote the trace map of the residual complex $(W_nF_Y)^\triangle K_{n,Y}$:
$$\Tr_{W_nf}:(W_nf)_*(W_nF_X)^\triangle K_{n,X}\simeq (W_nf)_*(W_nf)^\triangle (W_nF_Y)^\triangle K_{n,Y}\ra (W_nF_Y)^\triangle K_{n,Y}.$$
The commutativity of the first diagram is due to the functoriality of the trace map with respect to residual complexes with the same associated filtration (\Cref{trace for residual}(3)). The commutativity of the second is because of the compatibility of the trace map with compositions of morphisms (\Cref{trace for residual}(4)). 
\epf
\subsubsection{$\Tr_{W_nf,log}$ in the case of a nilpotent immersion}

\bprop[R\"ulling. {Cf. \cite[4.2]{Kato-DualitypPrimaryEtale-II}}]\label{Kato's complex thickening iso}
Let $i:X_0\hra X$ be a nilpotent immersion (thus so is $W_ni:W_n(X_0)\ra W_nX$). Then the natural map
$$\Tr_{W_ni,log}:(W_ni)_*K_{n,X_0,log}\ra K_{n,X,log}$$
is a quasi-isomorphism.
\eprop
\bpf
Put ${I_n}:=\Ker(W_n\cO_X\ra (W_ni)_*W_n\cO_{X_0})$.
Applying $\cHom_{W_n\cO_X}(-,K_{n,X})$ to the sequence of $W_n\cO_X$-modules
\begin{equation}\label{ses before Hom(-,KX)}
0\ra {I_n}\ra W_n\cO_X\ra (W_ni)_*W_n\cO_{X_0}\ra 0,
\end{equation}
we get again a short exact sequence of complexes of $W_n\cO_X$-modules
$$0 \ra (W_ni)_*K_{n,X_0}\xra{\Tr_{W_ni}} K_{n,X}\ra Q_n:=\cHom_{W_n\cO_X}({I_n},K_{n,X})\ra 0.$$
The first map is clearly $\Tr_{W_ni}$ by duality. The restriction of the map $(W_nF_X)^*:W_n\cO_X\ra (W_nF_X)_*W_n\cO_X$ to ${I_n}$ gives a map
\begin{align*}
(W_nF_X)^*\mid_{I_n}:\quad {I_n} &\ra (W_nF_X)_*{I_n},\\
\sum_{i=0}^{n-1} V([a_i])&\mapsto \sum_{i=0}^{n-1} V([a_i^p]).
\end{align*}
Define
\begin{align}
C'_{{I_n}}:
(W_nF_X)_*Q_n&= (W_nF_X)_*\cHom_{W_n\cO_X}({I_n},K_{n,X})
\label{C'In}
\\
\nonumber
&\ra \cHom_{W_n\cO_X}((W_nF_X)_*{I_n},(W_nF_X)_* K_{n,X})\\
\nonumber
&\xra[\simeq]{(W_nF_X)_*(\ref{(1.0.2)})\circ } \cHom_{W_n\cO_X}((W_nF_X)_*{I_n},(W_nF_X)_*(W_nF_X)^\triangle K_{n,X})\\
\nonumber
&\xra{\Tr_{W_nF_X}\circ}\cHom_{W_n\cO_X}((W_nF_X)_*{I_n}, K_{n,X})\\
\nonumber
&\xra{((W_nF_X^*)\mid_{I_n})^\vee} \cHom_{W_n\cO_X}({I_n}, K_{n,X})=Q_n.
\end{align}
According to the definition of $C'$ in (\ref{C'_n def}), $C'$ is compatible with $C'_{I_n}$. Thus one has the following commutative diagram
$$
\xymatrix@C=4em{
0 \ar[r]&
(W_nF_X)_*(W_ni)_*K_{n,X_0}\ar[r]^{(W_nF_X)_*\Tr_{W_ni}}\ar[d]^{C'}&
(W_nF_X)_*K_{n,X}\ar[r]\ar[d]^{C'}&
(W_nF_X)_*Q_n\ar[r]\ar[d]^{C'_{{I_n}}}&
0
\\
0 \ar[r]&
(W_ni)_*K_{n,X_0}\ar[r]^{\Tr_{W_ni}}&
K_{n,X}\ar[r]&
Q_n\ar[r]&
0.
}$$
Replacing $C'$ by $C'-1$, and $C'_{I_n}$ by $C'_{I_n}-1$, we arrive at the two lower rows of the following diagram.
Denote
$$Q_{n,log}:=\Cone (Q_n\xra{C'_{{I_n}}-1}Q_n)[-1].$$
Taking into account the shifted cones of $C'-1$ and $C'_{I_n}-1$, we get the first row of the following diagram which is naturally a short exact sequence. Now we have the whole commutative diagram of complexes, where all the three rows are exact, and all the three columns are distinguished triangles in the derived category:
\begin{displaymath}
\xymatrix@C=3em{
0\ar[r]&
(W_ni)_*K_{n,X_0,log}
\ar[r]^{\Tr_{W_ni,log}}\ar[d]&
K_{n,X,log}
\ar[d]\ar[r]&
Q_{n,log}
\ar[d]\ar[r]&
0\\
0\ar[r]&(W_ni)_*K_{n,X_0}\ar[d]^{C'-1}\ar[r]^{\Tr_{W_ni}}&K_{n,X}\ar[d]^{C'-1}\ar[r]&Q_n\ar[d]^{C'_{{I_n}}-1}\ar[r]&0
\\
0\ar[r]&(W_ni)_*K_{n,X_0}\ar[d]^{+1}\ar[r]^{\Tr_{W_ni}}&
K_{n,X}\ar[d]^{+1}\ar[r]&Q_n\ar[d]^{+1}\ar[r]&0.
\\
&&&
}
\end{displaymath}

We want to show that $\Tr_{W_ni,log}$ is a quasi-isomorphism. By the exactness of the first row, it suffices to show that $Q_{n,log}$ is an acyclic complex.
Because the right column is a distinguished triangle, it suffices to show that $C'_{{I_n}}-1:Q_n\ra Q_n$ is a quasi-isomorphism. Actually it is even an isomorphism of complexes: since $(W_nF_X)^*\mid_{I_n}:{I_n}\ra (W_nF_X)_*{I_n}$ is nilpotent (because $I_1=\Ker(\cO_X\ra i_*\cO_{X_0})$ is a finitely generated nilpotent ideal of $\cO_X$), the map $C'_{{I_n}}:Q_n\ra Q_n$ is therefore nilpotent (because one can alter the order of the three labeled maps in (\ref{C'In}) in the obvious sense), and $C'_{{I_n}}-1$ is therefore an isomorphism of complexes. 
\epf

\subsubsection{Localization triangles associated to $K_{n,X,log}$}

Let $i:Z\hra X$ be a closed immersion with $j:U\hra X$ being its open complement. Recall 
\begin{equation}\label{Gamma_Z zar def}
\underline\Gamma_Z(\cF):=\Ker(\cF\ra j_*j^{-1}\cF)
\end{equation}
for any abelian sheaf $\cF$. Denote its $i$-th derived functor by $\cH^i_Z(\cF)$.
Notice that
\bit
\item
$\Gamma_{Z'}(\cF)=\Gamma_Z(\cF)$ for any nilpotent thickening $Z'$ of $Z$ (e.g. $Z'=W_nZ$),
\item
$\cF\ra j_*j^{-1}\cF$ is surjective whenever $\cF$ is flasque, and
\item
flasque sheaves are $\uline\Gamma_Z$-acyclic (\cite[1.10]{Hartshorne-LocCoh}) and $f_*$-acyclic for any morphism $f$.
\eit
Therefore, for any complex of flasque sheaves $\cF^\bullet$ of $\Z/p^n$-modules on $W_nX$,
$$0\ra \uline\Gamma_{Z}(\cF^\bullet)\ra \cF^\bullet\ra (W_nj)_*(\cF^\bullet|_{W_nU})\ra 0$$
is a short exact sequence of complexes. Thus the induced triangle
\begin{equation}\label{Localization triangle KnXlog: preliminary}
\uline\Gamma_{Z}(\cF^\bullet)\ra \cF^\bullet\ra (W_nj)_*(\cF^\bullet|_{W_nU})\xra{+1}
\end{equation}
is a distinguished triangle in $D^b(W_nX,\Z/p^n)$, whenever $\cF^\bullet$ is a flasque complex with bounded cohomologies. In particular, since $K_{n,X,log}$ is a bounded complex of flasque sheaves, this is true for $\cF^\bullet=K_{n,X,log}$.

The following proposition is proven in the smooth case by Gros-Milne-Suwa \cite[2.6]{Suwa-NoteOnGerstenForlogSheaves}. The proof presented here comes from an unpublished manuscript of R\"ulling.

\bprop[R\"{u}lling]\label{localization triangle KnXlog}
Let $i:Z\hra X$ be a closed immersion with $j:U\hra X$ its open complement. Then
\begin{enumerate}
\item
(Purity)
The map
$$(W_ni)_*K_{n,Z,log}=\underline\Gamma_Z((W_ni)_*K_{n,Z,log})
\xra{\Tr_{W_ni,log}}\underline\Gamma_Z(K_{n,X,log})$$
is a quasi-isomorphism of complexes of sheaves.
\item (Localization triangle)
The following
\begin{equation}\label{adj dis triangle asso to Kato's complex}
(W_ni)_*K_{n,Z,log}\xra{\Tr_{W_ni,log}}K_{n,X,log}\ra (W_nj)_*K_{n,U,log}\xra{+1}
\end{equation}
is a distinguished triangle in $D^b(W_nX,\Z/p^n)$.
\end{enumerate}
\eprop
Note that we are working on the Zariski site and abelian sheaves on $W_nX$ can be identified with abelian sheaves on $X$ canonically. Thus we can replace $(W_ni)_*K_{n,Z,log}$ by $i_*K_{n,Z,log}$, and $(W_nj)_*K_{n,U,log}$ by $j_*K_{n,U,log}$ freely.
\bpf
\begin{enumerate}
\item
Let $I_n$ be the ideal sheaf associated to the closed immersion $W_ni:W_nZ\hra W_nX$, and let $Z_{n,m}$ be the closed subscheme of $W_nX$ determined by $m$-th power ideal $I_n^m$. In particular, $Z_{n,1}=W_nZ$. Denote by $i_{n,m}:Z_{n,m}\hra W_nX$ and by $j_{n,m}:W_nZ\hra Z_{n,m}$ the associated closed immersions.
In this way, for each $m$, one has a decomposition of $W_ni$ as maps of $W_nk$-schemes:
$$\xymatrix@C=6em{
W_nZ\ar@{^(->}[r]^{j_{n,m}}\ar[drr]_(.4){W_n\pi_Z} &
Z_{n,m}\ar@{^(->}[r]^{i_{n,m}}\ar[dr]^(.65){\pi_{Z_{n,m}}}
&W_nX\ar[d]^{W_n\pi_X}
\\
&&W_nk.
}$$

Denote $K_{Z_{n,m}}:=(\pi_{Z_{n,m}})^\triangle (W_nk)$, where $\pi_{Z_{n,m}}:Z_{n,m}\ra W_nk$ is the structure morphism. 
We have a canonical isomorphism
\beq\label{ind sys}
i_{n,m,*}\cH^i(K_{Z_{n,m}})
\simeq
\cExt^i_{W_n\cO_X}(i_{n,m,*}\cO_{Z_{n,m}},K_{n,X})
\eeq 
by \Cref{f triangle}(4) and \Cref{f triangle}(1) associated to the closed immersion $i_{n,m}$. 
The trace maps associated to the closed immersions
$$Z_{n,m}\hra Z_{n,m+1}$$
for different $m$ make the left hand side of (\ref{ind sys}) an inductive system.
The right hand side also lies in an inductive system when $m$ varies: the canonical surjections
$$i_{n,m+1,*}\cO_{Z_{n,m+1}}\ra i_{n,m,*}\cO_{Z_{n,m}}$$
induce the maps
\beq\label{ind sys hom}
\cHom_{W_n\cO_X}(i_{n,m,*}\cO_{Z_{n,m}},K_{n,X})\ra \cHom_{W_n\cO_X}(i_{n,m+1,*}\cO_{Z_{n,m+1}},K_{n,X})
\eeq
whose $i$-th cohomologies are the connecting homomorphisms of the inductive system. By duality, the map (\ref{ind sys hom}) is the trace map associated to the closed immersion $Z_{n,m}\hra Z_{n,m+1}$, and thus is compatible with the inductive system on the left hand side of (\ref{ind sys}).

Consider the trace map associated to the closed immersion $i_{n,m}:Z_{n,m}\hra W_nX$, i.e., the evaluation-at-1 map
$$\cHom_{W_n\cO_X}(i_{n,m,*}\cO_{Z_{n,m}},K_{n,X})\ra 
K_{n,X}.
$$
Its image naturally lies in $\Gamma_{W_nZ}(K_{n,X})$.
It induces an isomorphism on cohomology sheaves after taking the colimit on $m$ 
$$\colim_{m} \cExt^i_{W_n\cO_X}(i_{n,m,*}\cO_{Z_{n,m}},K_{n,X})
\xra[\simeq]{\ev_1}
\cH^i_Z(K_{n,X})$$
by \cite[V.4.3]{Hartshorne-RD}.

Now we consider
\begin{align}\label{eq: RD V S4 applied to KX}
\colim_{m}i_{n,m,*}\cH^i(K_{Z_{n,m}})
&\simeq
\colim_{m} \cExt^i_{W_n\cO_X}(i_{n,m,*}\cO_{Z_{n,m}},K_{n,X})
\\
&\xra[\simeq]{\ev_1}
\cH^i_Z(K_{n,X}).
\nonumber
\end{align}
The composite map of (\ref{eq: RD V S4 applied to KX}) is $\colim_m\Tr_{i_{n,m}}$. 
On the other hand, consider the log trace associated to the closed immersion $i_{n,m}$
(cf. \Cref{def of Tr f log})
\begin{align}\label{Tr in,log}
\Tr_{i_{n,m},log}:
\cH^i(i_{n,m,*}K_{Z_{n,m},log})
=\cH^i(\uline\Gamma_Z&(i_{n,m,*}K_{Z_{n,m},log}))\\
&\xra{}\cH^i(\uline\Gamma_Z(K_{n,X,log}))
=\cH^i_Z(K_{n,X,log}).
\nonumber
\end{align}
The maps (\ref{eq: RD V S4 applied to KX}), (\ref{Tr in,log}) give the vertical maps in the following diagram (due to formatting reason we omit $i_{n,m,*}$ from every term of the first row) which are automatically compatible by \Cref{def of Tr f log}:
$$
\xymatrix{
\cH^{i-1}(K_{Z_{n,m}})\ar[r]^{C'-1}\ar[d]_{\Tr_{i_{n,m}}}
&
\cH^{i-1}(K_{Z_{n,m}})\ar[r]\ar[d]_{\Tr_{i_{n,m}}}
&
\cH^i(K_{Z_{n,m},log})\ar[r]\ar[d]^{}_{\Tr_{i_{n,m},log}}
&
\cH^i(K_{Z_{n,m}})\ar[r]^{C'-1}\ar[d]_{\Tr_{i_{n,m}}}
&
\cH^i(K_{Z_{n,m}})\ar[d]_{\Tr_{i_{n,m}}}
\\
\cH^{i-1}_Z(K_{n,X})\ar[r]^{C'-1}
&\cH^{i-1}_Z(K_{n,X})\ar[r]
&\cH^i_Z(K_{n,X,log})\ar[r]
&\cH^i_Z(K_{n,X})\ar[r]^{C'-1}
&\cH^i_Z(K_{n,X}).
}
$$
Taking the colimit with respect to $m$, the five lemma immediately gives that $\colim_{m} \Tr_{i_{n,m},log}$ is an isomorphism. Then $\Tr_{W_ni,log}$, which is the composition of
$$
(W_ni)_*\cH^i(K_{n,Z,log})
\xra[\text{\Cref{Kato's complex thickening iso}},\simeq]{\colim_{m} \Tr_{j_{n,m},log}}
\colim_{m}i_{n,m,*} \cH^i(K_{Z_{n,m},log})
\xra[\simeq]{\colim_{m} \Tr_{i_{n,m},log}}
\cH^i_Z(K_{n,X,log}),
$$
is an isomorphism. This proves the statement.
\item
Since $\underline\Gamma_Z(K_{n,X,log})\ra K_{n,X,log}\ra (W_nj)_*K_{n,U,log}\xra{+1}$ is a distinguished triangle, the second part follows from the first part.
\end{enumerate}

\epf

\subsection{Functoriality}\label{subsection functoriality}
The push-forward functoriality of $K_{n,X,log}$ has been done in \Cref{def of Tr f log} for proper $f$. Now we define the pullback map for an \'etale morphism $f$. Since $W_nf$ is then also \'etale, we have an isomorphism of functors $(W_nf)^*\simeq (W_nf)^\triangle$ by \Cref{f triangle}(2). Define a chain map of complexes of $W_n\cO_Y$-modules
\begin{equation}\label{f* Kn,X}
f^*:K_{n,Y}\xra{\adj} (W_nf)_*(W_nf)^*K_{n,Y}\simeq (W_nf)_*(W_nf)^\triangle K_{n,Y}\simeq (W_nf)_*K_{n,X}.
\end{equation}
Here $\adj$ stands for the adjunction map of the identity map of $(W_nf)^*K_{n,Y}$.

\bprop[\'Etale pullback]\label{f* KnXlog}
Suppose $f:X\ra Y$ is an \'etale morphism. Then
$$f^*:K_{n,Y,log}\ra (W_nf)_*K_{n,X,log},$$
defined by termwise applying (\ref{f* Kn,X}), is a chain map between complexes of abelian sheaves.
\eprop
\bpf
It suffices to prove that $C'$ is compatible with $f^*$ defined above. Consider the following diagram in the category of complexes of $W_n\cO_Y$-modules
$$
\resizebox{\displaywidth}{!}{
\xymatrix@C=2em{
(W_nF_Y)_*K_{n,Y}
\ar@{}[dr]|{a)}
\ar[d]^{\adj}\ar[r]^(.45){(\ref{(1.0.2)})}_(.45){\simeq} &
(W_nF_Y)_*(W_nF_Y)^\triangle K_{n,Y}
\ar@{}[dr]^(.6){d)}
\ar[d]^{\adj}\ar[r]^(.6){\Tr_{W_nF_Y}} &
K_{n,Y}
\ar[d]^{\adj}
\\
(W_nf)_*(W_nf)^*(W_nF_Y)_* K_{n,Y}
\ar@{}[dr]|{b)}
\ar[r]^(.45){(\ref{(1.0.2)})}_(.45){\simeq}\ar[d]^{\simeq}_{\alpha}&
(W_nf)_*(W_nf)^*(W_nF_Y)_*(W_nF_Y)^\triangle K_{n,Y}
\ar@{}[ddr]^(.6){e)}
\ar[r]^(.6){\Tr_{W_nF_Y}}\ar[d]^{\simeq}_{\beta}&
(W_nf)_*(W_nf)^* K_{n,Y}
\ar[dd]^{\simeq}
\\
(W_nf)_*(W_nF_X)_*(W_nf)^* K_{n,Y}
\ar@{}[dr]|{c)}
\ar[d]^(.55){\simeq}
\ar[r]^(.45){(\ref{(1.0.2)})}_(.45){\simeq}&
(W_nf)_*(W_nF_X)_*(W_nf)^*(W_nF_Y)^\triangle K_{n,Y}
\ar[d]^(.55){\simeq}
&
\\
(W_nf)_*(W_nF_X)_*K_{n,X}
\ar[r]^(.45){(\ref{(1.0.2)})}_(.45){\simeq} &
(W_nf)_*(W_nF_X)_*(W_nF_X)^\triangle K_{n,X}
\ar[r]^(.6){\Tr_{W_nF_X}}&
(W_nf)_*K_{n,X}.
}}
$$
In this diagram we use shortened notations for the maps due to formatting reasons, e.g. we write $(\ref{(1.0.2)})$ instead of $(W_nf)_*(W_nF_X)_*(\ref{(1.0.2)})$, etc..
The maps labelled $\alpha$ and $\beta$ are base change maps, and they are isomorphisms because $W_nf$ is flat (actually $W_nf$ is \'etale because $f$ is \'etale) \cite[II.5.12]{Hartshorne-RD}.
The composite of the maps on the very left and very right are $(W_nF_Y)_*(f^*)$ and $f^*$ (where $f^*$ is as defined in (\ref{f* Kn,X})). The composite of the maps on the very top and very bottom are $C'_Y$ and $(W_nf)_*C'_X$.
Diagrams $a), b), c), d)$ commute due to naturality. 
Diagram $e)$ commutes, because we have a cartesian square
$$
\xymatrix@C=4em{ 
W_nX\ar[r]^{W_nF_X}\ar[d]_{W_nf}&W_nX\ar[d]^{W_nf}\\
W_nY\ar[r]^{W_nF_Y}&W_nY
}$$
by \Cref{etale over WnX}(2), and then the base change formula of the Grothendieck trace map as given in \Cref{trace for residual}(5) gives the result.
\epf

\ble\label{pull-push compatibility KnX zar}
Consider the following cartesian diagram
$$\xymatrix{
W\ar[d]^{g'}\ar[r]^{f'}&Z\ar[d]^g
\\
X\ar[r]^f&Y
}$$
with $g$ being proper, and $f$ being \'etale.
Then we have a commutative diagram of residual complexes
$$\xymatrix{
(W_ng)_*K_{n,Z}\ar[r]^{f'^*}\ar[d]^{\Tr_{W_ng}}&
(W_ng)_*(W_nf')_*K_{n,W}\ar[r]^{\simeq}&
(W_nf)_*(W_ng')_*K_{n,W}\ar[d]^{\Tr_{W_ng'}}
\\
K_{n,Y}\ar[rr]^{f^*}&&
(W_nf)_*K_{n,X}.
}$$
\ele
\bpf
We decompose the diagram into the following two diagrams and show their commutativity one by one. First we consider
$$\xymatrix{
(W_ng)_*K_{n,Z}\ar[r]^(.4){\adj}\ar[ddd]^{\Tr_{W_ng}}&
(W_ng)_*(W_nf')_*(W_nf')^*K_{n,Z}\ar[d]^(.4){\simeq}
\\
&(W_nf)_*(W_ng')_*(W_nf')^*(W_ng)^\triangle K_{n,Y}
\\
&
(W_nf)_* (W_nf)^*(W_ng)_*(W_ng)^\triangle K_{n,Y}
\ar[d]^{\Tr_{W_ng}}\ar[u]^{\simeq}_{\alpha}
\\
K_{n,Y}\ar[r]^(.4){\adj}&
(W_nf)_* (W_nf)^*K_{n,Y}.
}$$
Here $\alpha$ denotes the base change map, it is an isomorphism because $W_nf$ is flat \cite[II.5.12]{Hartshorne-RD}.  This diagram commutes by the naturality. Next we consider
$$\xymatrix{
(W_ng)_*(W_nf')_*(W_nf')^* K_{n,Z}\ar[d]^{\simeq}
\ar[r]^(.6){\simeq}&
(W_ng)_*(W_nf')_*K_{n,W}\ar[d]^{\simeq}
\\
(W_nf)_*(W_ng')_*(W_nf')^*(W_ng)^\triangle K_{n,Y}
\ar[r]^(.6){\simeq}
&
(W_nf)_*(W_ng')_*K_{n,W}\ar[dd]^{\Tr_{W_ng'}}
\\
(W_nf)_* (W_nf)^* (W_ng)_*(W_ng)^\triangle K_{n,Y}
\ar[d]^{\Tr_{W_ng}}\ar[u]^{\simeq}_{\alpha}
&
\\
(W_nf)_*(W_nf)^* K_{n,Y}
\ar[r]^{\simeq}&
(W_nf)_*K_{n,X}.}$$
The top part commutes by the naturality. The bottom part commutes by the base change formula of the Grothendieck trace maps with respect to \'etale morphisms (\Cref{trace for residual}(5)). 
\epf

Since both $f^*$ for log complexes in \Cref{f* KnXlog} and $g_*:=\Tr_{W_ng,log}$ are defined termwise, we arrive immediately the following compatibility as a consequence of \Cref{pull-push compatibility KnX zar}.

\bprop\label{pull-push compatibility KnXlog zar}
Notations are the same as \Cref{pull-push compatibility KnX zar}. One has a commutative diagram of complexes
$$\xymatrix{
(W_ng)_*K_{n,Z,log}\ar[r]^{f'^*}\ar[d]^{g_*}&
(W_ng)_*(W_nf')_*K_{n,W,log}\ar[r]^{\simeq}&
(W_nf)_*(W_ng')_*K_{n,W,log}\ar[d]^{g'_*}
\\
K_{n,Y,log}\ar[rr]^{f^*}&&
(W_nf)_*K_{n,X,log}.
}$$
\eprop

\subsection{\'Etale counterpart $K_{n,X,log,\et}$}\label{subsection etale counterpart KnXlog}
Let $X$ be a separated scheme of finite type over $k$ of dimension $d$.
In this subsection we will use $t=\Zar,\et$ to distinguish objects, morphisms on different sites. If $t$ is omitted, it means $t=\Zar$ unless otherwise stated.

Denote the structure sheaf on the small \'etale site $(W_nX)_\et$ by $W_n\cO_{X,\et}$. Denote $$(\epsilon_*,\epsilon^*):((W_nX)_\et,W_n\cO_{X,\et})\ra ((W_nX)_\Zar,W_n\cO_X)$$ the module-theoretic functors.
Recall that every \'etale $W_nX$-scheme is of the form $W_ng:W_nU\ra W_nX$, where $g:U\ra X$ is an \'etale $X$-scheme by \Cref{etale over WnX}(1).
Now let $\cF$ be a $W_n\cO_{X,\et}$-module on $(W_nX)_{\et}$. Consider the following map (cf. \cite[p. 264]{Kato-DualitypPrimaryEtale-II})
\begin{equation}\label{transfer of module structure: etale}\tau:(W_nF_X)_*\cF\xra{}\cF,
\end{equation}
which is defined to be
\begin{align*}
((W_nF_X)_*\cF)(W_nU\xra{W_ng} W_nX)&=\cF(W_nX\times_{W_nF_X,W_nX} W_nU\xra{pr_1} W_nX)
\\
&\xra[\simeq]{W_nF_{U/X}^*}\cF(W_nU\xra{W_ng} W_nX)
\end{align*}
for any \'etale map $W_ng:W_nU\ra W_nX$ (here we use $pr_1$ to denote the first projection map of the fiber product). This is an automorphism of $\cF$ as an abelian \'etale sheaf, but it changes the $W_n\cO_{X,\et}$-module structure of $\cF$.

Define
$$K_{n,X,\et}:=\epsilon^*K_{n,X}$$
to be the complex of \'etale $W_n\cO_{X,\et}$-modules associated to the Zariski complex $K_{n,X}$ of $W_n\cO_{X}$-modules. This is still a complex of quasi-coherent sheaves with coherent cohomologies. For a proper map $f:X\ra Y$ of $k$-schemes, define
$$
\Tr_{W_nf,\et}:(W_nf)_*K_{n,X,\et}= \epsilon^*((W_nf)_*K_{n,X})\xra{\epsilon^*\Tr_{W_nf}} K_{n,Y,\et}
$$
to be the \'etale map of $W_n\cO_{Y,\et}$-modules associated to the Zariski map $\Tr_{W_nf}:K_{n,X}\ra K_{n,X}$ of $W_n\cO_{X}$-modules.
Define the Cartier operator $C'_\et$ for \'etale complexes to be the composite
\begin{align*}
{C'_\et}:K_{n,X,\et} &
\xra[\simeq]{\tau^{-1}} (W_nF_X)_*K_{n,X,\et} =\epsilon^*((W_nF_X)_*K_{n,X} )
\xra{\epsilon^*(\ref{C'_n def})}
K_{n,X,\et}.
\end{align*}
Define
$$K_{n,X,log,\et}:= \Cone(K_{n,X,\et}\xra{C'_\et-1} K_{n,X,\et})[-1].$$

We also have the sheaf-level Cartier operator. Let $X$ be a smooth $k$-scheme. Recall that by definition, $C_\et$ is the composition of the inverse of (\ref{transfer of module structure: etale}) with the module-theoretic etalization of the $W_n\cO_X$-linear map (\ref{level n cartier def, top degree}) (it has appeared in \Cref{ses log forms are inv} before):
$$C_\et:
W_n\Omega^d_{X,\et}\xra[\simeq]{\tau^{-1}} (W_nF_X)_*W_n\Omega^d_{X,\et}=\epsilon^*((W_nF_X)_*W_n\Omega^d_{X})
\xra{\epsilon^*(\ref{level n cartier def, top degree})} W_n\Omega^d_{X,\et}.$$

\bprop[{cf. \Cref{Compatibility of $C'_n$ and $C$: Prop}}]\label{Compatibility of $C'_n$ and $C$ etale}
$C'_\et$ is the natural extension of $C'$ to the small \'etale site, i.e.,
$$\epsilon_*{C'_\et}={C'}:K_{n,X}\ra K_{n,X}.$$
If $X$ is smooth, $C_\et$ is the natural extension of $C$ to the small \'etale site
$$\epsilon_*{C_\et}={C}:W_n\Omega^d_X\ra W_n\Omega^d_X.$$
And one has compatibility
$${C_\et}=\cH^{-d}({C'_\et}).$$
\eprop
\bpf
The first two claims are clear. The last claim follows from the compatibility of $C$ and $C'$ in the Zariski case (\Cref{Compatibility of $C'_n$ and $C$: Prop}).
\epf

\bprop[{cf. \Cref{Raynaud-Illusie}}]\label{Raynaud-Illusie etale}
Let $X$ be a separated scheme of finite type over $k$ with $k=\bar k$. Then
$$H^i(W_nX,K_{n,X,\et}):=R^i\Gamma(W_nX,K_{n,X,\et})\xra{C'_\et-1} H^i(W_nX,K_{n,X,\et})$$
is surjective for every $i$. In particular,
$$R^i\Gamma(W_nX,K_{n,X,log,\et})\simeq 
H^i(W_nX,K_{n,X,\et})^{C'_\et-1}.$$
\eprop
\bpf
The quasi-coherent descent from the \'etale site to the Zariski site gives
$$
R\Gamma((W_nX)_{\et},K_{n,X,\et})
=R\Gamma((W_nX)_{\Zar},K_{n,X,\Zar}).
$$
Taking the $i$-th cohomology groups, the desired surjectivity then follows from the compatibility of $C'$ and $C'_{\et}$ (\Cref{Compatibility of $C'_n$ and $C$ etale}) and the Zariski case (\Cref{Raynaud-Illusie}).
\epf

In the \'etale topology and for any perfect field $k$, the surjectivity of
$$C_\et-1: W_n\Omega^d_{X,\et}\ra  W_n\Omega^d_{X,\et}$$
is known without the need of \Cref{Raynaud-Illusie etale} (cf. \Cref{ses log forms are inv}).
For the same reasoning as in \Cref{log qis complex log}, we have
\bprop[{cf. \Cref{log qis complex log}
}]\label{log qis complex log et}
Assume $X$ is smooth of pure dimension $d$ over a perfect field $k$. Then the natural map
$$
W_n\Omega^d_{X,log,\et}[d]\xra{} K_{n,X,log,\et}$$
is a quasi-isomorphism of complexes of abelian sheaves.
\eprop

We go back to the general non-smooth case. The proper pushforward property in the \'etale setting is very similar to the Zariski case.
\bprop[Proper pushforward, {cf. \Cref{def of Tr f log}}]\label{def of Tr f log et}
For $f:X\ra Y$ proper, we have a well-defined map of complexes of \'etale sheaves 
\beq\label{Tr_W_nf,log,et}
\Tr_{W_nf,log,\et}: (W_nf)_*K_{n,X,log,\et}\ra K_{n,X,log,\et}
\eeq
given by applying $\Tr_{W_nf,\et}$ termwise.
\eprop
\bpf
The map $\tau^{-1}$ is clearly functorial with respect to any map of abelian sheaves. The rest of the proof goes exactly as in \Cref{def of Tr f log}. 
\epf

\bprop[{cf. \Cref{Kato's complex thickening iso}}]
\label{Kato's complex thickening iso etale}
Let $i:X_0\hra X$ be a nilpotent immersion. Then the natural map
$$\Tr_{W_ni,log,\et}:(W_ni)_*K_{n,X_0,log,\et}\ra K_{n,X,log,\et}$$
is a quasi-isomorphism.
\eprop
\bpf
This is is a direct consequence of the functoriality of the map $\tau^{-1}$ and  \Cref{Kato's complex thickening iso}. 
\epf

Let $i:Z\hra X$ be a closed immersion with $j:U\hra X$ being the open complement as before. Define
$$\uline\Gamma_Z(\cF):=\Ker(\cF\ra j_*j^{-1}\cF)$$
for any \'etale abelian sheaf $\cF$ on $X$, just as in the Zariski case (cf. (\ref{Gamma_Z zar def})). Replacing $Z$ (resp. $X$) by a nilpotent thickening will define the same functor as $\uline\Gamma_Z(-)$, because the definition of the functor $\uline\Gamma_Z$ only depends on the pair $(X,U)$. Recall that if $\cF=\cI$ is an injective $\Z/p^n$-sheaf,
$$0\ra \uline\Gamma_Z(\cI)\ra\cI\ra j_*j^{-1}\cI\ra 0$$
is exact. In fact, because $j_!\Z/p^n$ is a subsheaf of the constant sheaf $\Z/p^n$ on $X$, 
the map $\Hom_X(\Z/p^n,\cI)\ra \Hom_X(j_!\Z/p^n,\cI)$ is surjective. 
Since $\Hom_X(j_!\Z/p^n,\cI)=\Hom_U(\Z/p^n,j^{-1}\cI)=\Hom_X(\Z/p^n,j_*j^{-1}\cI)$,
the map $\Hom_X(\Z/p^n,\cI)\ra \Hom_X(\Z/p^n,j_*j^{-1}\cI)$ is surjective, and hence we have the claim.
This implies that for any complex $\cF^\bullet$ of \'etale $\Z/p^n$-sheaves with bounded cohomologies,
\beq\label{Localization triangle KnXlog: etale}
R\uline\Gamma_Z(\cF^\bullet)\ra\cF^\bullet\ra j_*j^{-1}\cF^\bullet\xra{+1}
\eeq
is a distinguished triangle in $D^b(X,\Z/p^n)$ (cf. (\ref{Localization triangle KnXlog: preliminary})).
\bprop[cf. \Cref{localization triangle KnXlog}]\label{Localization triangle KnXlog etale}
Let $i:Z\hra X$ be a closed immersion with  open complement $j:U\hra X$, as before. Then
\benu
\item
(Purity)
We can identify canonically the functors $$(W_ni)_*=R\uline\Gamma_Z\circ(W_ni)_*:
D^b((W_nZ)_\et,\Z/p^n)\ra D^b((W_nX)_\et,\Z/p^n).$$
The composition of this canonical identification with the trace map 
$$(W_ni)_*K_{n,Z,log,\et}
=R\uline\Gamma_Z((W_ni)_*K_{n,Z,log,\et})
\xra{\Tr_{W_ni,log,\et}}
R\uline\Gamma_Z(K_{n,X,log,\et})$$
is a quasi-isomorphism of complexes of \'etale $\Z/p^n$-sheaves.
\item
(Localization triangle)
$$(W_ni)_*K_{n,Z,log,\et}\xra{\Tr_{W_ni,log,\et}}K_{n,X,log,\et}\xra{}
(W_nj)_*K_{n,U,log,\et}\xra{+1}$$
is a distinguished triangle in $D^b((W_nX)_{\et},\Z/p^n)$.
\eenu
\eprop
\bpf
\benu
\item
One only needs to show that $(W_ni)_*=R\uline\Gamma_Z\circ(W_ni)_*$, and then the rest of the proof is the same as in \Cref{localization triangle KnXlog}(1). Let $\cI$ be an injective \'etale $\Z/p^n$-sheaf on $W_nZ$. Since $\Hom_{W_nX}(-,(W_ni)_*\cI)=\Hom_{W_nZ}((W_ni)^{-1}(-),\cI)$ and $(W_ni)^{-1}$ is exact, we know $(W_ni)_*\cI$ is an injective abelian sheaf on $(W_nX)_\et$. This implies that $R(\uline\Gamma_Z\circ(W_ni)_*)=R\uline\Gamma_Z\circ(W_ni)_*$ by the Leray spectral sequence, and thus $(W_ni)_*=R(W_ni)_*
=R(\uline\Gamma_Z\circ(W_ni)_*)=R\uline\Gamma_Z\circ(W_ni)_*$.
\item
Note that $(W_nj)_*K_{n,U,log,\et}=R(W_nj)_*K_{n,U,log,\et}$. In fact, the terms of $K_{n,U,log,\et}$ are quasi-coherent $W_n\cO_{X,\et}$-modules which are $(W_nj)_*$-acyclic in the \'etale topology (because $R^if_*(\epsilon^*\cF)=\epsilon^*(R^if_*\cF)$ for any quasi-coherent Zariski sheaf $\cF$ and any quasi-compact quasi-separated morphism $f$ \cite[\href{https://stacks.math.columbia.edu/tag/071N}{Tag 071N}]{stacks-project}.). Now the first part and the distinguished triangle (\ref{Localization triangle KnXlog: etale}) imply the claim.
\eenu
\epf

\section{Bloch's cycle complex $\Z_{X,t}^c$}
\label{section bloch's cycle complex}
Let $X$ be a separated scheme of finite type over $k$ of dimension $d$.
Let
$$\Delta^i=\Spec k[T_0,\dots, T_i]/(\sum T_j-1).$$
Define $z_0(X,i)$ to be the free abelian group generated by closed integral subschemes $Z\subset X\times \Delta^i$ that intersect all faces properly and $\dim Z=i.$
We say two closed subschemes $Z_1,Z_2$ of a scheme $Y$ \emph{intersect properly} if for every irreducible component $W$ of the schematic intersection $Z_1\cap Z_2:=Z_1\times_Y Z_2$, one has
\beq\label{dimension inequality for intersect properly}
\dim W\le \dim Z_1 +\dim Z_2 - \dim Y
\eeq
(cf. \cite[A.1]{Geisser-MotKTopCycHol}).
A subvariety of $X\times \Delta^i$ is called a \emph{face} if it is determined by some $T_{j_1}=T_{j_2}=\dots= T_{j_s}=0$ ($0\le j_1<\dots <j_s\le i$). Note that a face is Zariski locally determined by a regular sequence of $X\times \Delta^i$. Therefore the given inequality condition (\ref{dimension inequality for intersect properly}) in the definition of $z_0(X,i)$ is equivalent to the equality condition \cite[(53)]{Geisser-MotKTopCycHol}.

The above definition defines a sheaf $z_0(-,i)$ in both the Zariski and the \'etale topology on $X$ (\cite[p.270]{Bloch-CycleAndK}. See also \cite[Lemma 3.1]{Geisser-MotCohOverDedekind}). 
Define the complex of sheaves
$$\ra z_0(-,i)\xra{d} z_0(-,i-1)\ra\dots z_0(-,0)\ra 0$$
with differential map
$$d(Z)=\sum_j(-1)^j[Z\cap V(T_j)].$$
Here we mean by $V(T_j)$ the closed integral subscheme determined by $T_j$ and by $[Z\cap V(T_j)]$ the linear combination of the reduced irreducible components of the scheme theoretic intersection $Z\cap V(T_j)$ with coefficients being intersection multiplicities. $z_0(X,\bullet)$ is then a homological complex concentrated in degree $[0,\infty)$.
Labeling cohomologically, we set
$$(\Z^c_X)^i=z_0(-,-i).$$
This complex is nonzero in degrees 
$$(-\infty,0].$$
Define the higher Chow group
\begin{align*}
\CH_0(X,i)
&:= H_{i}(z_0(X,\bullet))
=H^{-i}(\Z^c_X(X))
\end{align*}
for any $i$.
The higher Chow groups with coefficients in an abelian group $A$ will be denoted
\begin{align*}
\CH_0(X,i;A)
&:=
H^{-i}(\Z^c_X(X)\otimes_\Z A).
\end{align*}
The complex $\Z^c_{X,t}$, with either $t=\Zar$ or $t=\et$, has the following functoriality properties (cf. \cite[Prop. 1.3]{Bloch-CycleAndK}). If $f:X\ra Y$ is a proper morphism, then there is a chain map $f_*:f_*\Z^c_X\ra \Z^c_Y$ by the pushforward of cycles. If $f:X\ra Y$ is a quasi-finite flat morphism, then there is a chain map $f^*:\Z^c_Y\ra f_*\Z^c_X$ by the pullback of cycles.

\section{Kato's complex of Milnor $K$-theory ${C^M_{X,t}}$}
\label{section Milnor K}
Recall that given a field $L$, the $q$-th Milnor $K$-group $K^M_q(L)$ of $L$ is defined to be the $q$-th graded piece of the graded commutative ring
$$\bigoplus_{q\ge 0}K^M_q(L)=\frac{\bigoplus_{q\ge 0}(L^*)^{\otimes q}}{(a\otimes (1-a)\mid a,1-a\in  L^*)},$$
where $(a\otimes (1-a)\mid a,1-a\in  L^*)$ denotes the two-sided ideal of the graded commutative ring $\bigoplus_{q\ge 0}(L^*)^{\otimes q}$ generated by elements of the form $a\otimes (1-a)$ with $a,1-a\in  L^*$.
The image of an element $a_1\otimes\dots\otimes a_q\in (L^*)^{\otimes q}$ in $K^M_q(L)$ is denoted by $\{a_1,\dots, a_q\}$.

If $L$ is a discrete valuation field with valuation $v$ and residue field $k(v)$, 
the group homomorphism
$$\partial_v: K^M_q(L)\ra K^M_{q-1}(k (v)),\quad
\partial_v(\{\pi_v,u_1,\dots ,u_{q-1}\})=\{\bar u_1,\dots, \bar u_{q-1}\}
$$
is called the map of the \textit{tame symbol}.
Here $\pi_v$ is a local parameter with respect to $v$,  $u_1,\dots u_{q-1}$ are units in the valuation ring of $v$, and $\bar u_1,\dots, \bar u_{q-1}$ are the images of  $u_1,\dots u_{q-1}$ in the residue field $k(v)$. This is consistent with the sign convention in \cite[p.328]{Rost-ChowCoeff}.

For every natural number $q$ and every finite field extension $L'/L$, there exists a unique group homomorphism
$$\Nm_{L'/L}:K^M_q(L')\ra K^M_q(L)$$
such that 
\benu 
\item 
For any field extensions $L\subset L'\subset L''$, one has $\Nm_{L/L}=id$ and $\Nm_{L'/L}\circ \Nm_{L''/L'}=\Nm_{L''/L}$;
\item 
Let $L(T)$ be the function field of $\A^1_k$. For every $x\in K^M_q(L(T))$ one has 
$$\sum_v\Nm_{L(v)/L}(\partial_v(x))=0,$$
where $v$ runs over all discrete valuations of $L(T)$, and $L(v)$ denotes the residue field at valuation $v$.
\eenu 
The map $\Nm_{L'/L}$ is called the \textit{norm map} associated to the finite field extension $L'/L$.

Recall the definition of a Milnor $K$-sheaf on a point $X=\Spec L$, where $L$ is any field. 
$\cK^M_{\Spec L,q,\Zar}$ is the constant sheaf associated to the abelian group $K^M_q(L)$  (without the assumption that $L$ is an infinite field, cf. \cite[Prop. 10(4)]{Kerz-MilnorKFinResFields}), and $\cK^M_{\Spec L,q,\et}$ is the \'etale sheaf associated to the presheaf
$$L'\mapsto K^M_{q}(L'); \quad L'/L \text{ finite separable.}$$
Choose a separable closure $L^\sep$ of $L$. Then the geometric stalk at the geometric point $\Spec L^\sep$ over $\Spec L$ is $\colim_{L\subset L'\subset L^\sep} K^M_{q}(L')$, which is equal to $K^M_{q}(L^\sep)$ because the filtered colimit commutes with the tensor product and the quotient. Now by Galois descent of the \'etale sheaf condition, the sheaf $\cK^M_{\Spec L,q,\et}$ is precisely
$$L'\mapsto K^M_{q}(L^\sep)^{\Gal(L^\sep/L')}; \quad L'/L \text{ finite separable.}$$
Here the Galois action is given on each factor.

Let $X$ be a separated scheme of finite type over $k$ of dimension $d$.  
Now with the topology $t=\Zar$ or $t=\et$, we have the corresponding Gersten complex of Milnor $K$-theory, denote by $C^M_{X,t}$ (the differentials $d^M$ will be introduced below):
\begin{equation}\label{Kato's Milnor complex sheaf version}
\bigoplus_{x\in X_{(d)}} \iota_{x*}
\cK^M_{x,d,t}\xra{d^M}
\dots \xra{d^M}
\bigoplus_{x\in X_{(1)}} \iota_{x*}
\cK^M_{x,1,t} \xra{d^M}
\bigoplus_{x\in X_{(0)}} \iota_{x*}
\cK^M_{x,0,t},
\end{equation}
where $\iota_x: \Spec k(x)\hra X$ the natural inclusion map.
As part of the convention, $$(C^M_{X,t})^i=\bigoplus_{x\in X_{(-i)}}\iota_{x,*}\cK^M_{x,-i,t}.$$
In other words, (\ref{Kato's Milnor complex sheaf version}) sits in degrees 
$$[-d,0].$$
It remains to introduce the differential maps.

If $t=\Zar$, the differential map $d^M$ in (\ref{Kato's Milnor complex sheaf version}) is defined in the following way.
Let $x\in X_{(q)}$ be a dimension $q$ point, and $\rho:X'\ra \bar{\{x\}}$ be the normalization of $\bar{\{x\}}$ with generic point $x'$. Define
$$
(d^M)^x_y:
K^M_{q}(x)=K^M_{q}(x')
\xra{\sum \partial^{x'}_{y'}} \bigoplus_{y'|y}K^M_{q-1}(y')
\xra{\sum \Nm_{y'/y}}K^M_{q-1}(y).
$$
Here we have used the shortened symbol
$K^M_q(x):=K^M_q(k(x)).$ 
The notation $y'|y$ means that $y'\in {X'}^{(1)}$ is in the fiber of $y$.
\beq\label{tame symbol zar}
\partial^{x'}_{y'}: K^M_{q}(x')\ra K^M_{q-1}(y')
\eeq
is the Milnor tame symbol of the discrete valuation field $k(x')$ with valuation defined by $y'$. And
\beq\label{norm zar}
\Nm_{y'/y}: K^M_{q-1}(y')\ra K^M_{q-1}(y)
\eeq
is the Milnor norm map of the finite field extension $k(y)\subset k(y')$. 
The differential $d^M$ of this complex is given by
$$d^M:=\sum_{x\in X_{(q)}}\sum_{y\in X_{(q-1)}\cap \bar{\{x\}}}
(d^M)^x_y:
\bigoplus_{x\in X_{(q)}}
K^M_{q}(x)
\ra
\bigoplus_{y\in X_{(q-1)}}
K^M_{q-1}(y).
$$

If $t=\et$, set $x\in X_{(q)}$, $y\in X_{(q-1)}\cap \bar{\{x\}}$. Denote by $\rho:X'\ra \bar{\{x\}}$ the normalization map and denote by $x'$ the generic point of $X'$. One can canonically identify the \'etale abelian sheaves 
$\cK^M_{x,q,\et}$ and $\rho_*\cK^M_{x',q,\et}$
on $\bar{\{x\}}$ (here $\cK^M_{x,q,\et}$ on $\bar{\{x\}}$ means the pushforward of the sheaf $\cK^M_{x,q,\et}$ on the point $\Spec k(x)$ via $\Spec k(x)\ra \bar{\{x\}}$),
and thus identify
$\iota_{x,*}\cK^M_{x,q,\et}$ and $\iota_{x,*}\rho_*\cK^M_{x',q,\et}$
on $X$.
Let $y'\in X'^{(1)}$ such that $\rho(y')=y$. Then the componentwise differential map
$$(d^M)^x_{y}:\iota_{x,*}\cK^M_{x,q,\et}\ra \iota_{y,*}\cK^M_{y,q-1,\et}$$
is defined to be the composition
$$(d^M)^x_{y}=\iota_{y,*}(\Nm)\circ \rho_*(\partial).$$
Here $\partial:=\sum_{y'\in X'^{(1)}\cap \rho^{-1}(y)} \partial^{x'}_{y'}$, where
\begin{equation}\label{tame symbol etale}
\partial^{x'}_{y'}:\iota_{x',*}\cK^{M}_{x',q,\et}\ra \iota_{y',*}\cK^{M}_{y',q-1,\et}
\end{equation}
on $X'$ is defined to be the sheafification of the tame symbol on the presheaf level. Indeed, the tame symbol is a map of \'etale presheaves by \cite[R3a]{Rost-ChowCoeff}.
And $\Nm:=\sum_{y'\in X'^{(1)}\cap \rho^{-1}(y)}\Nm_{y'/y}$, where
\begin{equation}\label{norm etale}
\Nm_{y'/y}: \rho_*\cK^{M}_{y',q-1,\et}\ra \cK^{M}_{y,q-1,\et}
\end{equation}
on $y$ is defined to be the sheafification of the norm map on the presheaf level. The norm map is a map of \'etale presheaves by \cite[R1c]{Rost-ChowCoeff}.

The complex $C^M_{X,t}$, either $t=\Zar$ or $t=\et$, is covariant for proper morphisms  and contravariant for quasi-finite flat morphisms  (\cite[(4.6)(1)(2)]{Rost-ChowCoeff}). The pushforward map associated to a proper morphism is induced by the Milnor norm map, and the pullback map associated to a quasi-finite flat morphism is induced by the pullback map of the structure sheaves.

\section{Kato-Moser's complex of logarithmic de Rham-Witt sheaves $\tilde \nu_{n,X,t}$}
\label{section Gersten complex of logarithmic de Rham-Witt sheaves}
Kato first defined the Gersten complex of the logarithmic de Rham-Witt sheaves in \cite[\S1]{KatoKCT-HassePrinciple}. Moser in \cite[(1.3)-(1.5)]{Moser} sheafified Kato's construction on the \'etale site and studied its dualizing properties. 
We will adopt here the sign conventions in \cite{Rost-ChowCoeff}. 

Let $Y$ be a $k$-scheme. Let $q \in \N$ be an integer. Recall that in \Cref{subsubsection Comparison between W_nOmega^d log and KnXlog}, we have defined $W_n\Omega_{Y,log,t}^q$, with either $t=\Zar$ or $t=\et$, to be the abelian subsheaf of $W_n\Omega_{Y,t}^q$ \'etale locally generated by log forms.
We will freely use %
$W_n\Omega^q_{L,log,t}$ for $W_n\Omega^q_{\Spec L,log,t}$ below.

Now let $X$ be a separated scheme of finite type over $k$ of dimension $d$. Define the Gersten complex $\tilde \nu_{n,X,t}$, in the topology $t=\Zar$ or $t=\et$, to be the complex of $t$-sheaves isomorphic to $C^M_{X,t}/p^n$ via the Bloch-Gabber-Kato isomorphism \cite[2.8]{Bloch-Kato-pAdicEtaleCohomology}:
\begin{equation}\label{tilde vnX}
0\ra
\bigoplus_{x\in X_{(d)}} \iota_{x,*}W_n\Omega_{k(x),log,t}^{d}\ra
\dots \ra
\bigoplus_{x\in X_{(1)}} \iota_{x*}W_n\Omega_{k(x),log}^{1} \ra
\bigoplus_{x\in X_{(0)}} \iota_{x,*}W_n\Omega_{k(x),log,t}^{0}
\ra 0.
\end{equation}
Here $\iota_x:\Spec k(x)\ra X$ is the natural map. We will still denote by $\partial$ the reduction of the tame symbol $\partial$ mod $p^n$ (cf. (\ref{tame symbol zar})(\ref{tame symbol etale})), but denote by $\tr$ the reduction of Milnor's norm $\Nm$ mod $p^n$ (cf. (\ref{norm zar})(\ref{norm etale})). The reason for the later notation will be clear from \Cref{Compatibility of Milnor norm and Grothendieck trace}. As part of the convention,
$$\tilde \nu_{n,X,t}^i=\bigoplus_{x\in X_{(-i)}} \iota_{x*}W_n\Omega_{k(x),log,t}^{-i},$$
i.e. $\tilde \nu_{n,X}$ is concentrated in degrees
$$[-d,0].$$

\bprop\label{Localization triangle tilde vnX}
Let $i:Z\hra X$ be a closed immersion with $j:U\hra X$ its open complement.
We have the following short exact sequence for $t=\Zar$:
$$
\xymatrix{
0\ar[r] &i_*\tilde \nu_{n,Z,\Zar}\ar[r]&\tilde \nu_{n,X,\Zar}\ar[r]& j_*\tilde \nu_{n,U,\Zar}\ar[r]&0.
}
$$
For $t=\et$, one has the localization triangle
$$
i_*\tilde \nu_{n,Z,\et}\ra \tilde \nu_{n,X,\et}\ra Rj_*\tilde \nu_{n,U,\et}\xra{+1}.
$$
\eprop
\bpf
$\tilde \nu_{n,X,\Zar}$ is a complex of flasque sheaves (therefore $Rj_*(\tilde \nu_{n,X,\Zar})=j_*\tilde \nu_{n,X,\Zar}$), and one has the sequence being short exact in this case.
If $t=\et$, the purity theorem holds \cite[Corollary on p.130]{Moser}, i.e., $i_*\tilde \nu_{n,Z,\et}=\uline\Gamma_Z(\tilde \nu_{n,X,\et})\xra{\simeq} R\uline\Gamma_Z(\tilde \nu_{n,X,\et})$.
We are done with the help of the distinguished triangle (\ref{Localization triangle KnXlog: etale}) in the \'etale topology.
\epf

Functoriality of $\tilde \nu_{n,X,t}$ is the same as that of $C^M_{X,t}$ via $d\log$. We omit the statement.

\part{The maps}
\section{Construction of the chain map $\zeta_{n,X,log,t}: {C^M_{X,t}} \xra{}  K_{n,X,log,t}$}\label{section construction of chain map zeta}
\subsection{Construction of the chain map $\zeta_{n,X,t}: {C^M_{X,t}} \xra{}  K_{n,X,t}$}
Let $x\in X_{(q)}$ be a dimension $q$ point. 
$\iota_x:\Spec k(x)\ra X$ is the canonical map and $i_x:\bar{\{x\}}\hra X$ the closed immersion.
At degree $i=-q$, and over a point $x$, we define the degree $i$ map to be $\zeta_{n,X,t}^i:=\sum_{x\in X_{(q)}}\zeta_{n,x,t}^i$, with
\begin{align}\label{zeta t}
\zeta_{n,x,t}^{i}:
(W_n\iota_x)_*\cK^{M}_{x,q,t}&\xra[]{d\log}
(W_n\iota_x)_*W_n\Omega_{k(x),log,t}^{q}
\subset
(W_n\iota_x)_*W_n\Omega_{k(x),t}^{q}\\
&=(W_ni_x)_*K_{n,\bar{\{x\}},t}^{i}
\xra{(-1)^i\Tr_{W_ni_x}}
K_{n,X,t}^{i}.
\nonumber
\end{align}
We will use freely the notation $\zeta_{n,X,t}^i$ with some of its subscript or superscript dropped.

It is worth noticing that all the maps of \'etale sheaves involved here are given by the sheafification of the respective Zariski maps on the \'etale presheaf level. So to check commutativity of a composition of such maps between \'etale sheaves, it suffices to check on the $t=\Zar$ level. Keeping the convention as before, we usually omit the subscript $\Zar$ if we are working with the Zariski topology.

\bprop\label{chain map zata:CMXt to KnXt}
Let $X$ be a separated scheme of finite type over $k$ with $k$ being a perfect field of characteristic $p>0$. For $t=\Zar$ and $t=\et$, the map
$$\zeta_{n,X,t}: {C^M_{X,t}} \xra{}  K_{n,X,t},$$
as defined termwise in (\ref{zeta t}), is a chain map of complexes of sheaves on the site $(W_nX)_t$.
\eprop
Note that we have a canonical identification $(W_nX)_t=X_t$ for both $t=\Zar$ and $t=\et$. We use $(W_nX)_t$ just for the convenience of describing the $W_n\cO_X$-structure of residual complexes appearing later.
\bpf
To check $\zeta_{n,X,t}$ is a map of complexes, it suffices to check that the diagram
$$
\xymatrix{
(C^M_{X,t})^{i}\ar[r]^{d^M_X}\ar[d]^{\zeta_{n,X,t}^{i}}
&(C^M_{X,t})^{i+1}\ar[d]^{\zeta_{n,X,t}^{i+1}}\\
(K_{n,X,t})^{i}\ar[r]^{d_X}&(K_{n,X,t}^{i+1})^{i+1}
}
$$
commutes for $t=\Zar$. To this end, it suffices to show: for each $x\in X_{(q)}$, and $y\in X_{(q-1)}$ which is a specialization of $x$, the diagram
\begin{equation}\label{CommDiag1}
\xymatrix{
(W_n\iota_x)_*\cK^M_{x,q}\ar[r]^{(d_X^M)^x_y}
\ar[dd]^{\zeta_{n,x}}&
(W_n\iota_y)_*\cK^M_{y,q-1}\ar[d]^{\zeta_{n,y}}
\\
&
(W_ni_{y,x})_{*}K_{n,\bar{\{y\}}}^{-q+1}
\ar[d]^{-\Tr_{W_ni_{y,x}}}
\\
K_{n,\bar{\{x\}}}^{-q}
\ar[r]^{d_{\bar{\{x\}}}}
&K_{n,\bar{\{x\}}}^{-q+1}
}
\end{equation}
commutes. Here $i_{y,x}:\bar{\{y\}}\hra \bar{\{x\}}$ denotes the canonical closed immersion. 

Since the definition of the differential maps in ${C^M_{X}} $ involves normalization, consider the normalization $\rho:X'\ra \bar{\{x\}}$ of $\bar{\{x\}}$, and form the cartesian square
\begin{displaymath}
\xymatrix{
\bar{\{y\}}\times_{\bar{\{x\}}} X'\ar[d]\ar@{^{(}->}[r]
&X'=\bar{\{x'\}}\ar[d]_{\rho}
\\
\bar{\{y\}}\ar@{^{(}->}[r]^{i_{y,x}}&\bar{\{x\}}.
}
\end{displaymath}
Denote the generic point of $X'$ by $x'$.
Suppose $y'$ is one of the generic points of the irreducible components of $\bar{\{y\}}\times_{\bar{\{x\}}} X'$, and denote by $Y'$ the irreducible component corresponding to $y'$.
In particular, $y'$ is a codimension 1 point in the normal scheme $X'$,
thus is regular. Because the base field $k$ is perfect, $y'$ is also a smooth point in $X'$. According to \Cref{partial Ekedahl for singular}, the degree $[-q,-q+1]$ terms of $K_{n,X'}$ are of the form
$$(W_n\iota_{x'})_*H^0_{x'}(W_n\Omega_{X'}^{q})
\xra{\delta}\bigoplus_{y'\in X'_{(q-1)}}
(W_n\iota_{y'})_*H^1_{y'}(W_n\Omega_{X'}^{q})
\xra{}\dots,$$
where $\delta$ denotes the differential map of the residual complex $K_{n,X}$. 
After localizing at a single $y'\in X'^{(1)}$ in the Zariski sense, one gets
$$(W_n\iota_{x'})_*H^0_{x'}(W_n\Omega_{X'}^{q})
\xra{\delta_{y'}}
(W_n\iota_{y'})_*H^1_{y'}(W_n\Omega_{X'}^{q})
\xra{}\dots.$$

Consider the following diagrams.
Write $\iota_{x'}: \Spec k(x')\hra X'$, $\iota_{y'}:\Spec k(y')\hra X'$ the inclusions, $i_{y',x'}:Y'=\bar{\{y'\}}\hra X'$ the closed immersion, we have a diagram
\begin{equation}\label{circle 2,3}
\xymatrix{
(W_n\iota_{x'})_*\cK^M_{x',q}
\ar[r]^{
\partial^{x'}_{y'}}\ar[d]^{d\log}&
(W_n\iota_{y'})_*\cK^M_{y',q-1}\ar[d]^{d\log}
\\
(W_n\iota_{x'})_*W_n\Omega^{q}_{k(x')}
\ar[dr]^{\delta_{y'}}&
(W_n\iota_{y'})_*W_n\Omega^{q-1}_{k(y')}
\ar[d]^{
-\Tr_{W_n(i_{y',x'})}}
\\
&
(W_n\iota_{y'})_*H_{y'}^{1}(W_n\Omega_{X'}^{q}).
}
\end{equation}
For any $y'\in \rho^{-1}(y)\subset X'^{(1)}$, we have a diagram
\begin{equation}\label{circle 4}
\xymatrix@C=5em{
(W_n\rho)_*\cK^M_{y',q-1}
\ar[r]^{\Nm_{y'/y}}\ar[d]^{d\log}&
\cK^M_{y,q-1}\ar[d]^{d\log}\\
(W_n\rho)_*W_n\Omega^{q-1}_{k(y')}\ar[r]^{\Tr_{W_n\rho}}&
W_n\Omega^{q-1}_{k(y)}.
}
\end{equation}
Write $i_{y',x'}:Y'=\bar{\{y'\}}\hra X'$, $i_{y,x}:\bar{\{y\}}\hra \bar{\{x\}}$, we have a diagram
\begin{equation}\label{circle 6}
\xymatrix@C=5em{
(W_n\rho)_*(W_n\iota_{y'})_*W_n\Omega_{k(y')}^{q-1}
\ar[r]^{\Tr_{W_n\rho}}
\ar[d]_{\Tr_{W_n(i_{y',x'})}}&
(W_n\iota_{y})_*W_n\Omega_{k(y)}^{q-1}\ar[d]^{\Tr_{W_n(i_{y,x})}}
\\
(W_n\rho)_*(W_n\iota_{y'})_*H_{y'}^{1}(W_n\Omega_{X'}^{q})
\ar[r]^{\Tr_{W_n\rho}}&
K_{n,\bar{\{x\}}}^{-(q-1)},
}
\end{equation}
and a diagram
\begin{equation}\label{circle 5}
\xymatrix@C=5em{
(W_n\iota_{x'})_*W_n\Omega^{q}_{k(x')}
\ar[r]^(.38){d_{X'}=\sum\delta_{y'}}\ar[d]_{\Tr_{W_n\rho}}^{\simeq}&
\bigoplus_{y'\in \rho^{-1}(y)}
(W_n\iota_{y'})_* H^1_{y'}(W_n\Omega^{q-1}_{X'})
\ar[d]^{\Tr_{W_n\rho}}
\\
K_{n,\bar{\{x\}}}^{-q}
\ar[r]^{d_{\bar{\{x\}}}}&
K_{n,\bar{\{x\}}}^{-(q-1)}.
}
\end{equation}
All the trace maps above are trace maps of residual complexes at a certain degree.
(\ref{circle 6}) is the degree $q-1$ part of the diagram
$$
\xymatrix{
(W_n\rho)_*(W_ni_{y',x'})_*K_{n,Y'}
\ar[r]^{\Tr_{W_n\rho}}
\ar[d]_{\Tr_{W_n(i_{y',x'})}}&
(W_ni_{y,x})_*K_{n,\bar{\{y\}}}\ar[d]^{\Tr_{W_n(i_{y,x})}}
\\
(W_n\rho)_*K_{n,X'}
\ar[r]^{\Tr_{W_n\rho}}&
K_{n,\bar{\{x\}}}
}
$$
(the trace map on top is the trace map of the restriction of $W_n\rho$ to $W_nY'$),
and thus is commutative by the functoriality of the Grothendieck trace map with respect to composition of morphisms (\Cref{trace for residual}(4)).
(\ref{circle 5}) is simply the degree $-q$ to $-q+1$ terms of the trace map $\Tr_{W_n\rho}: (W_n\rho)_*K_{n,X'}\ra K_{n,\bar{\{x\}}}$, thus is also commutative. 
It remains to check the commutativity of (\ref{circle 2,3}) and (\ref{circle 4}). And these are \Cref{Compatibility of the Tame symbol and the residue symbols zhiqiande lemma} and \Cref{Compatibility of Milnor norm and Grothendieck trace}.

One notices that diagram (\ref{CommDiag1}) decomposes into the four diagrams (\ref{circle 2,3})-(\ref{circle 5}):
$$
\xymatrix@C=3em{
(W_n\rho)_*(W_n\iota_{x'})_*\cK^M_{x',q}
\ar[r]^{\bigoplus_{y'|y}
\partial^{x'}_{y'}}\ar[d]^{d\log}&
\bigoplus_{y'|y}
(W_n\rho)_*(W_n\iota_{y'})_*\cK^M_{y',q-1}
\ar[d]^{d\log}
\ar[r]^(.62){\sum_{y'/y}\Nm_{y'/y}}&
(W_n\iota_{y})_*\cK^M_{y,q-1}
\ar[d]^{d\log}
\\
(W_n\rho)_*(W_n\iota_{x'})_*W_n\Omega^{q}_{k(x')}
\ar[dd]^{\Tr_{W_n\rho}}
\ar@{}[ur]|{(\ref{circle 2,3})}
\ar[dr]^{\bigoplus_{y'|y}\delta_{y'}}&
\bigoplus_{y'|y}
(W_n\rho)_*(W_n\iota_{y'})_*W_n\Omega^{q-1}_{k(y')}
\ar[d]^{\bigoplus_{y'|y}
-\Tr_{W_n(i_{y',x'})}}
\ar[r]^(.62){\Tr_{W_n\rho}}
\ar@{}[ur]|(.65){(\ref{circle 4})}
&
(W_n\iota_{y})_*W_n\Omega^{q-1}_{k(y)}
\ar[d]^{-\Tr_{W_n(i_{y,x})}}
\\
&
\bigoplus_{y'|y}
(W_n\rho)_*(W_n\iota_{y'})_*H_{y'}^{1}(W_n\Omega_{X'}^{q})
\ar[r]^(.62){\Tr_{W_n\rho}}
\ar@{}[ur]|(.65){(\ref{circle 6})}&
K_{n,\bar{\{x\}}}^{-(q-1)}
\ar@{=}[d]
\\
K_{n,\bar{\{x\}}}^{-q}
\ar[rr]^{d_{\bar{\{x\}}}}
\ar@{}[urr]|{(\ref{circle 5})}
&&
K_{n,\bar{\{x\}}}^{-(q-1)}
}
$$
Here by symbol $y'|y$ we mean that $y'\in \rho^{-1}(y)$. Notice that we have added a minus sign to both vertical arrows of (\ref{circle 6}) in the corresponding square above, but this does not affect its commutativity.
Since one can canonically identify
$$(W_n\rho)_*(W_n\iota_{x'})_*\cK^M_{x',q} \text{\quad with \quad}(W_n\iota_{x})_*\cK^M_{x,q},$$
to show the commutativity of the diagram (\ref{CommDiag1}), it only remains to show \Cref{Compatibility of the Tame symbol and the residue symbols zhiqiande lemma} and \Cref{Compatibility of Milnor norm and Grothendieck trace}.
\epf

\ble\label{Compatibility of the Tame symbol and the residue symbols zhiqiande lemma}
For an integral normal scheme $X'$, with $x'\in X'$ being the generic point and $y'\in X'^{(1)}$ being a codimension 1 point, the diagram (\ref{circle 2,3}) is commutative.
\ele
\bpf
Given a $y'\in X'^{(1)}$ lying over $y$,
the abelian group $K^M_{q}(x')$ is generated by
$$\{\pi',u_1,\dots ,u_{q-1}\} \text{ and }
\{v_1,\dots,v_{q-1},v_{q}\},$$
where $u_1,\dots  ,u_{q-1}, v_1,\dots,v_{q-1},v_{q}\in \cO_{X',y'}^*$, and $\pi'$ is a chosen uniformizer of the discrete valuation ring $\cO_{X',y'}$. It suffices to check the commutativity for these generators.

In the first case, the left-bottom composition gives
\begin{align*}
(\delta_{y'}\circ d\log)
(\{\pi',u_1,\dots ,u_{q-1}\})&=
\delta_{y'}(d\log[\pi']_n  d\log[u_1]_n \dots  d\log[u_{q-1}]_n)\\
&=
\left[
\txt{
$d[\pi']_n  d\log[u_1]_n \dots  d\log[u_{q-1}]_n$\\
$[\pi']_n$
}\right].
\end{align*}
The last equality above is given by \cite[A.1.2]{CR11}.
Here we have used the fact that $[\pi']$ is a regular element in $W_nX'$, since $\pi'$ is regular in $X'$.
The top-right composition gives
\begin{align*}
(-\Tr_{W_n(i_{y',x'})}\circ &d\log  \circ \partial^{x'}_{y'}) (\{\pi',u_1,\dots ,u_{q-1}\})
\\
&=
(-\Tr_{W_n(i_{y',x'})}\circ d\log)
\{\bar u_1,\dots,\bar u_{q-1}\}
\\
&=
-\Tr_{W_n(i_{y',x'})}
(d\log[\bar u_1]_n \dots  d\log[\bar u_{q-1}]_n)
\\
&=
\left[
\txt{
$d[\pi']_n
d\log[\bar u_1]_n \dots  d\log[\bar u_{q-1}]_n
$\\
$[\pi']_n$
}\right].
\end{align*}
The last equality is given by \cite[A.2.12]{CR11}.
So the diagram (\ref{circle 2,3}) is commutative in this case.

In the second case, since $\partial^{x'}_{y'}(\{v_1,\dots,v_{q}\})=0$, we need to check the left-bottom composite also gives zero. In fact,
\begin{align*}
(\delta_{y'}\circ d\log)(\{v_1,\dots,v_{q}\})&
=
\delta_{y'}(d\log[v_1]_n \dots  d\log [v_{q}]_n)
\\
&=
\left[
\txt{
$[\pi']_n\cdot d\log[v_1]_n \dots  d\log [v_{q}]_n$\\
$[\pi']_n$
}\right]
\\&=0.
\end{align*}
The second equality is due to \cite[A.1.2]{CR11}. The last equality is because, in a small neighborhood $V$ of $y'$, the element $[\pi']_n\cdot d\log[v_1]_n \dots  d\log [v_{q}]_n\in W_n\Omega_{V}^q$ lies in the $W_n\cO_V$-submodule $[\pi']_n\cdot W_n\Omega_{V}^q$. 
\epf

\ble[Compatibility of Milnor norm and Grothendieck trace]\label{Compatibility of Milnor norm and Grothendieck trace} 
Let $F/E$ be a finite field extension with both fields $E$ and $F$ being of transcendence degree $q-1$ over $k$. Suppose there exists a finite morphism $g$ between integral separated finite type $k$-schemes, such that $F$ is the function field of the source of $g$ and $E$ is the function field of the target of $g$, and the field extension $F/E$ is induced via the map $g$. Then the following diagram commutes
$$
\xymatrix@C=4em{
K^M_{q-1}(F)\ar[r]^{\Nm_{F/E}}\ar[d]^{d\log}&
K^M_{q-1}(E)\ar[d]^{d\log}\\
W_n\Omega_{F}^{q-1} \ar[r]^{\Tr_{W_ng}}&W_n\Omega_{E}^{q-1}.
}
$$
Here the norm map $\Nm_{F/E}$ denotes the norm map from Milnor $K$-theory, and $\Tr_{W_ng}$ denotes the Grothendieck trace map associated to the finite morphism $g$.
\ele
\brmk
The compatibility of the trace map with the norm and the pushforward of cycles in various settings has been known by the experts, and many definitions and properties of the trace map in the literature reflect this viewpoint. 
But since we have not found a proof of the compatibility of the Milnor norm with the trace map defined via the Grothendieck duality theory, we include a proof here.
\ermk
\bpf 
We start the proof by some reductions.
Since both $\Nm_{F/E}$ and $\Tr_{F/E}$ are independent of the choice of towers of simple field extensions, without loss of generality, one can suppose $F$ is a finite simple field extension over $E$. Now $F=E(a)=\frac{E[T]}{f(T)}$ for some monic irreducible polynomial $f(T)\in E[T]$ with $a\in F$ being one of its roots. This realizes $\Spec F$ as an $F$-valued point $P$ of $\P^1_E$, namely,
$$
\xymatrix{
*+[l]{\Spec F=P\,}\ar@{^{(}->}[r]^{i_P}\ar[dr]_{g}&\P^1_E\ar[d]^{\pi}\\
&\Spec E.
}
$$
All the three morphisms on above are morphisms of finite type (although not between schemes of finite type over $k$), so it makes sense to talk about the associated trace maps for residual complexes.
But for the particular residual complexes we are interested in, we need to enlarge the schemes involved to schemes of finite type over $k$, while preserving the morphism classes (e.g., closed immersion, smooth morphism, etc) of the morphisms between them.

To this end, take $Y$ to be any separated smooth connected scheme of finite type over $k$ with $E$ being the function field.
Since $\P^1_E$ is the generic fiber of $Y\times_k \P^1_k$, 
by possibly shrinking $Y$ to an affine neighborhood $\Spec B$ of $pr_1(P)$ (here $pr_1: Y\times_k \P^1_k\ra Y$ is the first projection map) one can extend the above diagram to the following:
$$\xymatrix@C=3em{
*+[l]{\Spec F\in W}
\ar@{^{(}->}[r]^{i_W}\ar[dr]_g&\P^1_Y\ar[d]^\pi\\
&*+[r]{Y=\Spec B\ni \Spec E.}
}$$
Here $W:=\bar{\{P\}}^{\P^1_Y}$ is the closure of the point $P$ in $\P^1_Y$. This is a commutative diagram of finite type $k$-schemes.
In particular, it makes sense to talk about the residual complexes $K_{n,Y}, K_{n,W}$ and $K_{n,\P^1_Y}$. 

Now it remains to show the commutativity of the following diagram
\begin{equation}\label{norm trace}
\xymatrix@C=4em{
K^M_{q-1}(E(a))\ar[r]^{\Nm_{E(a)/E}}\ar[d]^{d\log}&
K^M_{q-1}(E)\ar[d]^{d\log}\\
W_n\Omega_{E(a)}^{q-1} \ar[r]^{\Tr_{W_ng}}&W_n\Omega_{E}^{q-1},
}
\end{equation}
where $\Tr_{W_ng}$ denotes the trace map for residual complexes $\Tr_{W_ng}:(W_ng)_*K_{n,W}\ra K_{n,Y}$ at degree $-(q-1)$.

We do induction on $[E(a):E]$.
If $[E(a):E]=1$, then both the Grothendieck trace $\Tr_{W_ng}:W_n\Omega^{q-1}_{E(a)/k}\ra W_n\Omega^{q-1}_{E/k}$ and the norm map $\Nm_{E(a)/E}:K^M_{q-1}(E(a))\ra K^M_{q-1}(E)$ are the identity, therefore the claim holds. Now the induction step. Suppose the diagram (\ref{norm trace}) commutes for $[E(a):E]\le r-1$. We will need to prove the commutativity for $[E(a):E]=r$.

First note that $\Tr_{W_ng}:(W_ng)_*K_{n,W}\ra K_{n,Y}$
naturally decomposes into
\beq\label{Trg=TriW circ Trpi}
(W_ng)_*K_{n,W}\xra{(W_n\pi)_*\Tr_{W_ni_PW}}
(W_n\pi)_*K_{n,\P^1_Y}
\xra{\Tr_{W_n\pi}} K_{n,Y}.
\eeq
by \Cref{trace for residual}(4).
$H^1_P(W_n\Omega^q_{\P^1_Y})$ is a direct summand of the degree $-(q-1)$ part of $K_{n,\P^1_Y}$. One can canonically identify
\beq\label{local P1Y=local P1E}
H^1_P(W_n\Omega^q_{\P^1_Y})=H^1_P(W_n\Omega^q_{\P^1_E}),
\eeq
via pulling back along the natural map $\P^1_E\hra \P^1_Y$.
Thus on degree $-(q-1)$ and at the point $P$, the map (\ref{Trg=TriW circ Trpi}) is canonically identified with
$$W_n\Omega_{E(a)}^{q-1}
\xra{\Tr_{W_ni_PW}}
H^1_P(W_n\Omega_{\P^1_E}^{q})
\xra{\Tr_{W_n\pi}}
W_n\Omega_{E}^{q-1}.$$

Consider the diagram
$$
\xymatrix@C=4em{
K^M_{q}(E(T))
\ar@{}[rdd]|{(-1)}
\ar[r]^{\partial_P}\ar[dd]^{d\log}& K^M_{q-1}(E(a))\ar[r]^{\Nm_{E(a)/E}}\ar[d]^{d\log}&
K^M_{q-1}(E)\ar[dd]^{d\log}
\\
&W_n\Omega_{E(a)}^{q-1}
\ar[d]^{\Tr_{W_ni_PW}}\ar[dr]^{\Tr_{W_ng}}&
\\
W_n\Omega_{E(T)}^{q}
\ar[r]^{\delta_P}&
H^1_P(W_n\Omega_{\P^1_E}^{q})
\ar[r]^(0.6){\Tr_{W_n\pi}}&
W_n\Omega_{E}^{q-1}.
}
$$
We have used the identification (\ref{local P1Y=local P1E}) in this diagram. We have seen that the left square is commutative up to sign $-1$, as a special case of \Cref{Compatibility of the Tame symbol and the residue symbols zhiqiande lemma} (i.e. take normal scheme $X'=\P^1_E$ and $y':=P=\Spec F$). 
Since $\partial_P$ is surjective, to show the commutativity of the trapezoid on the right, it suffices to show that the composite square is commutative up to $-1$. 
For any element
$$s:=\{s_1,\dots,s_{q-1}\}\in K^M_{q-1}(E(a)),$$
one can always find a lift
$$\tilde s:=\{f,\tilde s_1,\dots,\tilde s_{q-1}\}\in K^M_{q}(E(T)),$$
such that each of the $s_i=s_i(T)$ is a polynomial of degree $\le r-1$ 
(e.g. decompose $E(a)$ as a $r$-dimensional $E$-vector space $E(a)=\bigoplus_{j=0}^{r-1} Ea^j$ and suppose $s_i=\sum_{j=0}^{r-1} b_{i,j}a^j$ with $b_{i,j}\in E$, then $\tilde s_i=\tilde s_i(T)=\sum_{j=0}^{r-1} b_{i,j}T^j$ satisfies the condition), 
and $\partial_P(\tilde s)=s$. Denote by $$y_{i,1},\dots, y_{i,a_i}\quad (1\le i\le q-1)$$  the closed points of $\P^1_E$ corresponding to the irreducible factors of the polynomials $\tilde s_1,\dots,\tilde s_{q-1}$. Note that the local section $\tilde s_{i,l}$ cutting out $y_{i,l}$ is by definition an irreducible factor of $\tilde s_i$, and therefore $\deg \tilde s_{i,l}< r$ for all $i$ and all $l$.

We claim that 
\begin{equation}\label{sum is zero delta_P circ Tr}
\sum_{y\in (\P^1_E)_{(0)}} (\Tr_{W_n\pi})_y\circ \delta_{y} =0:W_n\Omega_{E(T)/k}^{q}\ra W_n\Omega_{E/k}^{q-1}.
\end{equation}
In fact,
\beq\label{ses Witt} 
0\ra W_n\Omega_{\P^1_E}^{q}\ra W_n\Omega_{E(T)}^{q}\ra \bigoplus_{y\in (\P^1_E)_{(0)}} (W_n\iota_y)_*H^1_y(W_n\Omega_{\P^1_E}^{q})\ra 0
\eeq
is an exact sequence \cite[1.5.9]{CR12}, where $\iota_y:y\hra \P^1_E$ is the natural inclusion of the point $y$. Taking the long exact sequence with respect to the global section functor, one arrives at the following diagram with the row being a complex 
$$\xymatrix{
W_n\Omega_{E(T)}^{q}
\ar[r]^(.3){\delta} &
\bigoplus_{y\in (\P^1_E)_{(0)}} H^1_y(W_n\Omega_{\P^1_E}^{q})
\ar[r]^{} \ar[dr]_{\sum_{y}(\Tr_{W_n\pi})_y}&
H^1(\P^1_E,W_n\Omega_{\P^1_E}^{q})\ar[d]^{\Tr_{W_n\pi}}
\\
&&W_n\Omega_{E/k}^{q-1}.
}$$
The trace maps on left of the above are induced from the degree $0$ part of $\Tr_{W_n\pi}:(W_n\pi)_*K_{n,\P^1_Y}\ra K_{n,Y}$. 
The trace map on the right of the above is induced also by $\Tr_{W_n\pi}:(W_n\pi)_*K_{n,\P^1_Y}\ra K_{n,Y}$, while the global cohomology group is calculated via (\ref{ses Witt}), 
i.e., one uses the last two terms of (\ref{ses Witt}) as an injective resolution of the sheaf $W_n\Omega_{\P^1_E}^q$, and then  $\Tr_{W_n\pi}:(W_n\pi)_*K_{n,\P^1_Y}\ra K_{n,Y}$ induces the map of complexes (sitting in degrees $[-1,0]$) on global sections, 
and then the map of cohomologies on degree $0$ gives our trace map $H^1(\P^1_E,W_n\Omega_E^q)\ra W_n\Omega^{q-1}_E$ on the right.
From the construction of these trace maps, the diagram on above is by definition commutative. Therefore (\ref{sum is zero delta_P circ Tr}) holds.

One notices that $
\delta_{y} \circ d\log(\tilde s)=0$ unless $y\in\{p,y_{1,1},\dots, y_{q-1,a_{q-1}},\infty\}$. Now we calculate
\begin{align*}
(\Tr_{W_ng}&\circ d\log)(s)\\
&=(\Tr_{W_ng}\circ d\log\circ \partial_P)(\tilde s)\\
&=-((\Tr_{W_n\pi})_P\circ \delta_P\circ d\log)(\tilde s)
\quad\text{(\Cref{Compatibility of the Tame symbol and the residue symbols zhiqiande lemma})}\\
&=\sum_{y\in\{y_{1,1},\dots, y_{q-1,a_{q-1}},\infty\}} ((\Tr_{W_n\pi})_y\circ\delta_{y} \circ d\log)(\tilde s)
\text{\quad (\ref{sum is zero delta_P circ Tr})}\\
&=-\sum_{y\in\{y_{1,1},\dots, y_{q-1,a_{q-1}},\infty\}}  (d\log \circ \Nm_{E(k(y))/E}\circ \partial_y)(\tilde s)\\
&\text{\qquad\qquad\qquad\qquad\qquad\qquad\qquad\qquad\qquad (induction hypothesis)}\\
&=(d\log\circ \Nm_{E(a)/E}\circ\partial_P)(\tilde s)
\text{\quad (\cite[2.2 (RC)]{Rost-ChowCoeff})}\\
&=(d\log\circ \Nm_{E(a)/E})(s).
\end{align*}
This finishes the induction.
\epf

\subsection{Functoriality of $\zeta_{n,X,t}:{C^M_{X,t}}\ra K_{n,X,t}$}

Let $k$ denote a perfect field of positive characteristic $p$.

\bprop[Proper pushforward]\label{zeta compatibility proper with push}
$\zeta$ is compatible with proper pushforward. I.e., for $f:X\ra Y$ a proper map, the following diagram is commutative
\begin{equation*}
\xymatrix@C=7ex{
(W_nf)_*{C^M_{X,t}} \ar[d]^{f_*}\ar[r]^{\zeta_{n,X,t}}& (W_nf)_*K_{n,X,t}\ar[d]^{f_*}\\
C_{Y,t}^M \ar[r]^{\zeta_{n,Y,t}} & K_{n,Y,t}.
}
\end{equation*}
Here $f_*$ on the left denotes the pushforward map for Kato's complex of Milnor $K$-theory (cf. \Cref{section Milnor K}), and $f_*$ on the right denotes the Grothendieck trace map $\Tr_{W_nf,t}$ for residual complexes.
\eprop

\bpf
We only need to prove the proposition for $t=\Zar$ and for degree $i\in[-d,0]$. Then by the very definition of the $\zeta$ map and the compatibility of the trace map with morphism compositions \Cref{trace for residual}(4), it suffices to check the commutativity at points $x\in X_{(q)}$, $y\in Y_{(q)}$, where $q=-i$:
$$
\xymatrix{
K^M_{q}(x)\ar[r]^{d\log}\ar[d]_{f_*}& W_n\Omega_{k(x)}^{q}
\ar[d]^{f_*}
\\
K^M_{q}(y)\ar[r]^{d\log}&
W_n\Omega_{k(y)}^{q}.
}
$$
\begin{enumerate}
\item
If $y\neq f(x)$, both pushforward maps are zero maps, therefore we have the desired commutativity.
\item
If $y=f(x)$, by definition of $\zeta$ and the pushforward maps, 
we need to show the commutativity of the following diagram for the finite field extension $k(y)\subset k(x)$
$$
\xymatrix{
K^M_{q}(x)\ar[r]^{d\log}\ar[d]_{\Nm_{k(x)/k(y)}}& W_n\Omega_{k(x)}^{q}
\ar[d]^{\Tr_{W_nf}}
\\
K^M_{q}(y)\ar[r]^{d\log}&
W_n\Omega_{k(y)}^{q}.
}
$$
This is precisely \Cref{Compatibility of Milnor norm and Grothendieck trace}. 
\end{enumerate}
\epf

\bprop[\'Etale pullback]\label{zeta compatibility etale pullback}
$\zeta$ is compatible with \'etale pullbacks. I.e., for $f:X\ra Y$ an \'etale morphism, the following diagram is commutative
\begin{equation*}
\xymatrix@C=4em{
C_{Y,t}^M \ar[d]^{f^*}\ar[r]^{\zeta_{n,Y,t}}& K_{n,Y,t}\ar[d]^{f^*}\\
(W_nf)_*{C^M_{X,t}} \ar[r]^{\zeta_{n,X,t}} &
(W_nf)_* K_{n,X,t}.
}
\end{equation*}
Here $f^*$ on the left denotes the pullback map for Kato's complex of Milnor $K$-theory (cf. \Cref{section Milnor K}), and $f^*$ on the right denotes the pullback map for residual complexes (\ref{f* Kn,X}).
\eprop

\bpf
It suffices to prove the proposition for $t=\Zar$. Take $y\in Y_{(q)}$. Consider the cartesian diagram
$$\xymatrix{
*+[l]{X\times_Y\bar{\{y\}}=:W}\ar[r]^{f|_W}\ar@{^(->}[d]^{i_W}&
\bar{\{y\}}\ar@{^(->}[d]^{i_{y}}
\\
X\ar[r]^f&Y.
}\qquad\qquad$$
Then the desired diagram at point $y$ decomposes in the following way at degree $-q$:
$$
\xymatrix@C=4em{
K^M_{q}(y)\ar[r]^(.4){d\log}\ar[d]_{f^*}&
W_n\Omega_{k(y)}^{q}=K_{n,\bar{\{y\}}}^{-q}
\ar[d]^{(f|_W)^*}\ar[r]^(.65){\Tr_{W_ni_y}}&
K_{n,Y}^{-q}\ar[d]^{f^*}
\\
\bigoplus_{x\in W_{(q)}}K^M_{q}(x)\ar[r]^(.4){d\log}&
\bigoplus_{x\in W_{(q)}}W_n\Omega_{k(x)}^{q}=
K_{n,W}^{-q}\ar[r]^(.65){\Tr_{W_ni_W}}&
K_{n,X}^{-q}.
}$$
The left square commutes because both $f^*$ and $(f\mid_W)^*$ are induced by the natural map $f^*:\cO_Y\ra f_*\cO_X$.
The right square commutes due to \Cref{pull-push compatibility KnX zar}.
\epf

\subsection{Extend to $K_{n,X,log,t}$}\label{section of extending phi to phi log}

Recall the complex $K_{n,X,log,t}:=\Cone (K_{n,X,t}\xra{C'_t-1}K_{n,X,t})[-1]$, i.e.,
$$K_{n,X,log,t}^i=K_{n,X,t}^{i}\oplus K_{n,X,t}^{i-1}.$$
Notice that
\begin{equation}\label{K to Klog map of graded algebras}
K_{n,X,t}\ra K_{n,X,log,t},\quad a\mapsto (a,0)
\end{equation}
is not a chain map. Nevertheless,
\bprop
We keep the same assumptions as in \Cref{chain map zata:CMXt to KnXt}. The chain map $\zeta_{n,X,t}: {C^M_{X,t}} \xra{}  K_{n,X,t}$
composed with (\ref{K to Klog map of graded algebras}) gives a chain map
$$\zeta_{n,X,log,t}:= (\ref{K to Klog map of graded algebras}) \circ \zeta_{n,X,t}:
{C^M_{X,t}} \xra{}  K_{n,X,log,t}$$
of complexes of abelian sheaves on $(W_nX)_t$.
\eprop
We will also use the shortened notation $\zeta_{log,t}$ for $\zeta_{n,X,log,t}$. If $t=\Zar$, the subscript $\Zar$ will also be omitted.
\bpf
Given $x\in X_{(q)}$, we prove commutativity of the following diagram
\begin{equation*}
\xymatrix@C=4em{
\iota_{x*}\cK^{M}_{x,q,t}\ar[r]^(.4){d\log}\ar@{=}[d]&
\iota_{x*}W_n\Omega_{k(x),log,t}^{q}\ar@{^{(}->}[d]
\ar[r]^{\Tr_{W_ni_x,t}}&
K_{n,X,t}^{-q}
\ar[r]^{C'_{X,t}-1}&
K_{n,X,t}^{-q}\ar@{=}[d]
\\
\iota_{x*}\cK^{M}_{x,q,t}\ar[r]^(.4){d\log}&
\iota_{x*}W_n\Omega_{k(x),t}^{q}
\ar[r]^{C'_{\bar{\{x\}},t}-1}&
(i_{x,*}K_{n,\bar{\{x\}},t})^{-q}
\ar[r]^(.6){\Tr_{W_ni_x,t}}&
K^{-q}_{n,X,t}.
}
\end{equation*}
The left square naturally commutes. The right square also commutes, because $C'$ is compatible with the Grothendieck trace map $\Tr_{W_ni_x}$ (the proofs of \Cref{def of Tr f log} and \Cref{def of Tr f log et} give the case for $t=\Zar$ and $t=\et$, respectively).
Now because $C'_{\bar{\{x\}},t}-1:W_n\Omega_{k(x),t}^q\ra W_n\Omega_{k(x),t}^q$, which is identified with $C_{\bar{\{x\}},t}-1$ as a result of \Cref{Compatibility of $C'_n$ and $C$: Prop} and \Cref{Compatibility of $C'_n$ and $C$ etale}, annihilates $W_n\Omega_{k(x),log,t}^{q}$, the composite of the second row is zero. Thus the composite of the first row is zero.
This yields a unique chain map
$$\zeta_{n,X,log,t}:{C^M_{X,t}} \ra K_{n,X,log,t},$$
i.e., on degree $i=-q$, we have $\zeta_{n,X,log,t}=\sum_{x\in X_{(q)}} \zeta_{n,x,log,t}$ with
\begin{align*}
\zeta_{n,x,log,t}^i:\cK^M_{x,q,t}&\ra K_{n,X,log,t}^i=K_{n,X,t}^i\oplus K_{n,x,t}^{i-1},\\
s=\{s_1\dots, s_{q}\}&\mapsto (\zeta_{n,X,t}^i(s),0).
\end{align*}
\epf

As a direct corollary of \Cref{zeta compatibility proper with push} and \Cref{zeta compatibility etale pullback}, one has the following proposition.
\bprop[Functoriality]\label{Functoriality Milnor to KnXlog}
\benu
\item
$\zeta_{log,t}$ is compatible with proper pushforward. I.e., for $f:X\ra Y$ a proper map, the following diagram of complexes is commutative
\begin{equation*}
\xymatrix@C=5em{
(W_nf)_*{C^M_{X,t}} \ar[d]^{f_*}\ar[r]^{\zeta_{n,X,log,t}}& (W_nf)_*K_{n,X,log,t}\ar[d]^{f_{*}}\\
C_{Y,t}^M \ar[r]^{\zeta_{n,Y,log,t}} & K_{n,Y,log,t}.
}
\end{equation*}
Here $f_*$ on the left denotes the pushforward map for Kato's complex of Milnor $K$-theory (cf. \Cref{section Milnor K}), and $f_*$ on the right denotes $\Tr_{W_nf,log,t}$ as defined in \Cref{def of Tr f log} and \Cref{def of Tr f log et}.
\item
$\zeta_{log,t}$ is compatible with \'etale pullbacks. I.e., for $f:X\ra Y$ an \'etale morphism, the following diagram of complexes is commutative
\begin{equation*}
\xymatrix@C=5em{
C_{Y,t}^M \ar[d]^{f^*}\ar[r]^{\zeta_{n,Y,log,t}}& K_{n,Y,log,t}\ar[d]^{f^*}\\
(W_nf)_*{C^M_{X,t}} \ar[r]^{\zeta_{n,X,log,t}} &
(W_nf)_* K_{n,X,log,t}.
}
\end{equation*}
Here $f^*$ on the left denotes the pullback map for Kato's complex of Milnor $K$-theory (cf. \Cref{section Milnor K}), and $f^*$ on the right denotes the pullback map defined in \Cref{f* KnXlog}.
\eenu
\eprop

\subsection{The map $\bar \zeta_{n,X,log,t}: {C^M_{X,t}} /p^n\simeq \tilde \nu_{n,X,t} \xra{}  K_{n,X,log,t}$ is a quasi-isomorphism}
Since $\zeta_{n,X,t}$ is termwise defined via the $d\log$ map, it annihilates $p^n C^M_{X,t}$. Therefore $\zeta_{n,X,log,t}$ annihilates $p^nC^M_{X,t}$ as well, and it induces a chain map
$$\bar \zeta_{n,X,log,t}:{C^M_{X,t}} /p^n\ra K_{n,X,log,t}.$$
Since the $d\log$ map induces an isomorphism of complexes ${C^M_{X,t}} /p^n\simeq \tilde \nu_{n,X,t}$, to show $\bar \zeta_{n,X,log,t}$ is a quasi-isomorphism, it is equivalent to show
$$ \bar \zeta_{n,X,log,t}:  \tilde\nu_{n,X,t}\ra K_{n,X,log,t}$$ is a quasi-isomorphism.

\ble\label{vertical arrow right}
Suppose $X$ is separated smooth over the perfect field $k$ of characteristic $p>0$. Then for any level $n$,
$$
\bar \zeta_{n,X,log,\et}:
\tilde \nu_{n,X,\et}\xra{} K_{n,X,log,\et};
$$
is a quasi-isomorphism. If we moreover have $k=\bar k$, then
$$
\bar \zeta_{n,X,log,\Zar}:
\tilde \nu_{n,X,\Zar}\xra{} K_{n,X,log,\Zar}
$$ is also a quasi-isomorphism.
\ele

\bpf
This is a local problem, thus it suffices to prove the statement for each connected component of $X$. Therefore we assume $X$ is of pure dimension $d$ over $k$. Then for any level $n$, we have a quasi-isomorphism (\cite[Cor 1.6]{GrosSuwa-LaConjDeGersten}) 
$$W_n\Omega_{X,log,t}^d[d]
\xra{\simeq} \tilde \nu_{n,X,t}.$$
We also have
\begin{align*}
W_n\Omega^d_{X,log,\et}[d]&\xra{\simeq}K_{n,X,log,\et} \quad
\text{(by \Cref{log qis complex log et}), and}
\\
W_n\Omega^d_{X,log,\Zar}[d]&\xra{\simeq}K_{n,X,log,\Zar}
\quad\text{if $k=\bar k$ (by \Cref{log qis complex log}).}
\end{align*}
On degree $-d$, we have a diagram
\begin{displaymath}
\xymatrix@C=4em{
\displaystyle
\tilde \nu_{n,X,t}^{-d}=
\bigoplus_{x\in X^{(0)}}
(W_n\iota_x)_*W_n\Omega^d_{k(x),log,t}
\ar[r]^{\bar\zeta_{n,x,log,t}^{-d}}&
\displaystyle
K_{n,X,log,t}^{-d}=
\bigoplus_{x\in X^{(0)}} (W_n\iota_x)_*H_x^0(W_n\Omega^d_{X,t})
\\
W_n\Omega^d_{X,log,t}\ar[r]^{(-1)^{d}}
\ar[u]&
W_n\Omega^d_{X,log,t}\ar[u]
}
\end{displaymath}
which is naturally commutative, due to the definition of $\bar\zeta_{n,X,log,t}$. It induces quasi-isomorphisms as stated in the lemma.
\epf

\bthm\label{Main theorem}
Let $X$ be a separated scheme of finite type over a perfect field $k$ of characteristic $p>0$. Then the chain map
$$
\bar \zeta_{n,X,log,\et}:\tilde \nu_{n,X,\et}\xra{} K_{n,X,log,\et}
$$
is a quasi-isomorphism. Moreover if $k=\bar k$,
$$
\bar \zeta_{n,X,log,\Zar}:\tilde \nu_{n,X,\Zar}\xra{} K_{n,X,log,\Zar}
$$
is also a quasi-isomorphism.
\ethm
\bpf
One can assume that $X$ is reduced. In fact, the complex $\tilde \nu_{n,X,t}$ is defined to be the same complex as $\tilde \nu_{n,X_{\mathrm{red}},t}$ (see (\ref{tilde vnX})), and we have a quasi-isomorphism $K_{n,X_{\mathrm{red}},log,t}\xra{\simeq} K_{n,X,log,t}$ given by the trace map, according to \Cref{Kato's complex thickening iso} and \Cref{Kato's complex thickening iso etale}. One notices that $\bar \zeta_{n,X_{\mathrm{red}},log,t}$ is compatible with $\bar \zeta_{n,X,log,t}$ because of the functoriality of the map $\zeta_{log,t}$ with respect to proper maps \Cref{Functoriality Milnor to KnXlog}(1). As long as we have a quasi-isomorphism
$$\bar \zeta_{n,X_{\mathrm{red}},log,t}:\tilde \nu_{n,X_{\mathrm{red}},t}\xra{}
K_{n,X_{\mathrm{red}},log},$$
we will get automatically that
$$\bar \zeta_{n,X,log,t}:
\tilde \nu_{n,X_{\mathrm{red}},t}=
\tilde \nu_{n,X,t}\xra{\bar \zeta_{n,X_{\mathrm{red}},log,t}}
K_{n,X_{\mathrm{red}},log}\xra{\simeq}
K_{n,X,log,t}$$
is a quasi-isomorphism.

Now we do induction on the dimension of the reduced scheme $X$. Suppose $X$ is of dimension $d$, and suppose $\bar \zeta_{n,Y,log,t}$ is a quasi-isomorphism for schemes of dimension $\le d-1$. Now decompose $X$ into the singular part $Z$ and the smooth part $U$
$$U\xhra{j} X\xhla{i} Z.$$
Then $Z$ has dimension $\le d-1$.
Consider the following diagram in the derived category of complexes of $\Z/p^n$-modules
\beq\label{triangle tilde vnXt to KnXlog}
\xymatrix@C=4em{
i_*\tilde \nu_{n,Z,t}\ar[r]\ar[d]^{i_*\bar\zeta_{n,Z,log,t}}& \tilde \nu_{n,X,t}\ar[r]\ar[d]^{\bar\zeta_{n,X,log,t}}&
Rj_*\tilde \nu_{n,U,t} \ar[d]^{Rj_*\bar\zeta_{n,U,log,t}}_{}\ar[r]^{+1}&
i_*\tilde \nu_{n,Z,t}[1]
\ar[d]^{i_*\bar\zeta_{n,Z,log,t}[1]}
\\
i_*K_{n,Z,log,t}\ar[r]^{\Tr_{W_ni,log}}&
K_{n,X,log,t}\ar[r] &
Rj_*K_{n,U,log,t}\ar[r]^{+1}&
i_*K_{n,Z,log,t}[1],
}
\eeq
where the two rows are distinguished triangles coming from \Cref{localization triangle KnXlog}, \Cref{Localization triangle KnXlog etale} and \Cref{Localization triangle tilde vnX}.
We show that the three squares in (\ref{triangle tilde vnXt to KnXlog}) are commutative in the derived category. The left square is commutative because of \Cref{Functoriality Milnor to KnXlog}(1).
The middle square is induced from the diagram
\beq\label{triangle tilde vnXt to KnXlog: middle}
\xymatrix{
\tilde \nu_{n,X,t}\ar[r]\ar[d]^{\bar\zeta_{n,X,log,t}}&
j_*\tilde \nu_{n,U,t} \ar[d]^{j_*\bar\zeta_{n,U,log,t}}
\\
K_{n,X,log,t}\ar[r] &
j_*K_{n,U,log,t}
}\eeq
of chain complexes. Let $x\in X_{(q)}$. If $x\in X_{(q)}\cap U$, both $\tilde \nu_{n,X,t} \ra j_*\tilde \nu_{n,U,t}$ and $K_{n,X,log,t}\ra j_*K_{n,U,log,t}$ give identity maps at $x$, therefore the square (\ref{triangle tilde vnXt to KnXlog: middle}) commutes in this case.
If $x\in X_{(q)}\cap Z$, both of these give the zero map at $x$, therefore the square (\ref{triangle tilde vnXt to KnXlog: middle}) is also commutative. The right square of (\ref{triangle tilde vnXt to KnXlog}) can be decomposed in the following way (cf. (\ref{Localization triangle KnXlog: preliminary}) and (\ref{Localization triangle KnXlog: etale})):
$$\xymatrix@C=4em{
Rj_*\tilde \nu_{n,U,t} \ar[d]^{Rj_*\bar\zeta_{n,U,log,t}}_{}\ar[r]^{+1}
&
R\uline\Gamma_Z(\tilde \nu_{n,X,t})[1]
\ar[d]^{R\uline\Gamma_Z(\bar\zeta_{n,X,log,t})[1]}
&\ar[l]_(.4){i_*}^(.4){\simeq}
i_*\tilde \nu_{n,Z,t}[1]
\ar[d]^{i_*\bar\zeta_{n,Z,log,t}[1]}
\\
Rj_*K_{n,U,log,t}\ar[r]^{+1}
&
R\uline\Gamma_Z(K_{n,X,log,t})[1]
&\ar[l]_(.4){i_*}^(.4){\simeq}
i_*K_{n,Z,log,t}[1].
}$$
The map $i_*$ on the first row is induced by the norm map of Milnor $K$-theory. It is clearly an isomorphism of complexes if $t=\Zar$. It is a quasi-isomorphism if $t=\et$ due to the purity theorem \cite[p.130 Cor.]{Moser}. The map $i_*$ on the second row is induced from $\Tr_{W_ni,log,t}$ as defined in \Cref{def of Tr f log} and \Cref{def of Tr f log et}, and it is an isomorphism due to \Cref{localization triangle KnXlog}(1) if $t=\Zar$, and \Cref{Localization triangle KnXlog etale} if $t=\et$.
The first square commutes by the naturality of the $+1$ map. The second commutes because of the compatibility of $\zeta_{log,t}$ with the proper pushforward \Cref{Functoriality Milnor to KnXlog}(1). We thus deduce that the right square of (\ref{triangle tilde vnXt to KnXlog}) commutes.

Now consider over any perfect field $k$ for either of the two cases:
\begin{enumerate}
\item
$t=\et$ and $k$ is a perfect field, or
\item
$t=\Zar$ and $k=\bar k$.
\end{enumerate}
The left vertical arrow of (\ref{triangle tilde vnXt to KnXlog}) is a quasi-isomorphism because of the induction hypothesis. The third one counting from the left is also a quasi-isomorphism because of \Cref{vertical arrow right}. Thus so is the second one. 
\epf

\section{Combine $\psi_{X,t}:\Z^c_{X,t}\ra {C^M_{X,t}}$  with $\zeta_{n,X,log,t}: {C^M_{X,t}} \xra{}  K_{n,X,log,t}$}
\label{section combine psi with zeta}

\subsection{The map $\psi_{X,t}:\Z^c_{X,t}\ra C^M_{X,t}$}
In \cite[2.14]{Zhong14}, Zhong constructed a map of abelian groups $\psi_{X,t}(X):\Z^c_{X}(X)\ra C^M_{X,\Zar}(X)$ based on the Nesterenko-Suslin-Totaro isomorphism \cite[Thm. 4.9]{NesterenkoSuslin}\cite{Totaro-MilnorKisSimplestPart}. 
It is straightforward to check that Zhong's construction induces a well-defined map  of complexes $\psi:=\psi_{X,t}:\Z^c_{X,t}\ra C^M_{X,t}$ of sheaves for both $t=\Zar$ and $t=\et$.
Zhong in \cite[2.15]{Zhong14} proved that $\psi$ is covariant with respect to proper morphisms, and contravariant with respect to quasi-finite flat morphisms.

\subsection{The map $\bar\zeta_{n,X,log,t}\circ \bar \psi_{X,t}: \Z^c_{X,t}/p^n\xra{\simeq} K_{n,X,log,t}$ is a quasi-isomorphism}
\label{subsection Zhong's qis combined with zeta}

In {\cite[2.16]{Zhong14}} Zhong proved that
$\psi_{X,\et}$
combined with the Bloch-Gabber-Kato isomorphism \cite[2.8]{Bloch-Kato-pAdicEtaleCohomology}, gives a quasi-isomorphism
$$\bar \psi_{X,\et}: \Z^c_{X,\et}/p^n\xra{\simeq} \tilde \nu_{n,X,\et}.$$
In the proof, Zhong actually showed that these two complexes of sheaves on each section of the big Zariski site over $X$ are quasi-isomorphic. Therefore by restriction to the small Zariski site, we have
$$\bar \psi_{X,\Zar}: \Z^c_{X,\Zar}/p^n\xra{\simeq} \tilde \nu_{n,X,\Zar}.$$

Combining Zhong's quasi-isomorphism with \Cref{Main theorem}:
\bthm\label{Main theorem 2}
Let $X$ be a separated scheme of finite type over $k$ with $k$ being a perfect field of positive characteristic $p$.
Then the following composition
$$
\bar\zeta_{n,X,log,\et}\circ \bar \psi_{X,\et}: \Z^c_{X,\et}/p^n\xra{\simeq} K_{n,X,log,\et},
$$
is a quasi-isomorphism. If moreover $k=\bar k$, then the following composition
$$
\bar\zeta_{n,X,log,\Zar}\circ \bar\psi_{X,\Zar}:
\Z^c_{X,\Zar}/p^n\xra{\simeq} K_{n,X,log,\Zar},
$$
is also a quasi-isomorphism.
\ethm

\brmk 
From the construction of the maps $\bar\zeta_{n,X,log,t}$ and $\bar \psi_{X,t}$, we can describe explicitly their composite map.
We write here only the Zariski case, and the \'etale case is just given by the Zariski version on the small \'etale site and then doing the \'etale sheafification.

Let $U$ be a Zariski open subset of $X$. Let $Z\in (\Z^c_{X,\Zar})^i(U)=z_0(U,-i)$ be a prime cycle.
\bit
\item
If $i\in[-d,0]$ and $\dim p_U(Z)=-i$, set $q=-i$. Then $Z$ as a cycle of dimension $q$ in $U\times \Delta^{q}$, is dominant over some $u=u(Z)\in U_{(q)}$ under the projection $p_U: U\times \Delta^{q}\ra U$. 
By slight abuse of notation, we denote by
$T_0,\dots, T_{q}\in 
k(Z)$ the pullbacks of the corresponding coordinates via $Z
\hra U\times \Delta^{q}$.
Since $Z$
intersects all faces properly, $T_0,\dots, T_{q}\in k(Z)^*$.
Thus $\{\frac{-T_0}{T_{q}},\dots,\frac{-T_{q-1}}{T_{q}}\} \in K^M_{q}(k(Z))$ is well-defined. 
Take the Zariski closure of $\Spec k(Z)$ in $U\times \Delta^{q}$, and denote it by $Z'$. Then $p_U$ maps $Z'$ to $\bar{\{u\}}^U=\bar{\{u\}}^X\cap U$. 
Denote by $i_u: \bar{\{u\}}^X\hra X$ the closed immersion, and denote the composition
$$Z'\xra{p_U} \bar{\{u\}}^U\hra \bar{\{u\}}^X\xhra{i_u} X$$
by $h$.
The map $h$ is clearly generically finite, then there exists an open neighborhood $V$ of $u$ in $X$ such that the restriction $h:h^{-1}(V)\ra V$ is finite. Then $W_nh: W_n(h^{-1}(V))\ra W_nV$ is also finite. Therefore it makes sense to consider the trace map $\Tr_{W_nh}$ near the generic point of $Z'$. Similarly, it makes sense to consider the trace map $\Tr_{W_np_U}$ near the generic point of $Z'$.
Then we calculate
\begin{align*}
\zeta_{log}(\psi(Z))&=
(-1)^i\Tr_{W_ni_u} (d\log (\Nm_{k(Z)/k(u(Z))}\{\frac{-T_0}{T_{q}},\dots,\frac{-T_{q-1}}{T_{q}}\}))
\\
&=(-1)^i\Tr_{W_ni_u}(\Tr_{W_np_U}
d\log\{\frac{-T_0}{T_{q}},\dots,\frac{-T_{q-1}}{T_{q}}\})
\quad(\text{\Cref{Compatibility of Milnor norm and Grothendieck trace}})
\\
&=(-1)^i\Tr_{W_nh}(
\frac{T_qdT_0-T_0dT_q}{T_0T_q}
 \cdots  
\frac{T_qdT_{q-1}-T_{q-1}dT_q}{T_{q-1}T_q}
)
\end{align*}
Here in the last step we have used the functoriality of the trace map with respect to composition of morphisms (\Cref{trace for residual}(4)).
\item
If $i\notin [-d,0]$ or $\dim p_U(Z)\neq -i$, we have $\zeta_{log}(\psi(Z))=0$.
\eit
\ermk

Combining the functoriality of Zhong's map $\psi$ with \Cref{Functoriality Milnor to KnXlog}, one arrives at the following proposition.
\bprop[Functoriality]
\label{Functoriality cycle to KnXlog}
The composition $\bar\zeta_{n,X,log,t}\circ \bar\psi_{X,t}:\Z^c_{X,t}/p^n\xra{\simeq} K_{n,X,log,t}$ is covariant with respect to proper morphisms, and contravariant with respect to \'etale morphisms for both $t=\Zar$ and $t=\et$.
\eprop

\part{Applications}
\section{De Rham-Witt analysis of $\tilde \nu_{n,X,t}$ and $K_{n,X,log,t}$}\label{section De Rham-Witt analysis}
Let $X$ be a separated scheme of finite type over $k$ of dimension $d$. In this section we will use terminologies as defined in \cite[\S1]{CR12}, such as Witt residual complexes, etc.

Recall that Ekedahl defined a map of complexes of $W_n\cO_X$-modules (cf. \cite[Def. 1.8.3]{CR12})
$$\uline p:=\uline p_{\{K_{n,X}\}_n}:R_*K_{n-1,X,t}\ra K_{n,X,t}.$$ 
By abuse of notations, we denote by $R:W_{n-1}X\hra W_nX$ the closed immersion induced by the restriction map on the structure sheaves $R:W_n\cO_X\ra W_{n-1}\cO_X$.
\ble\label{uline p for Kn,X,log lemma}
The map $\uline p:R_*K_{n-1,X,t}\ra K_{n,X,t}$ induces a map of complexes of abelian sheaves
\begin{equation}\label{uline p for Kn,X,log}
\uline p:K_{n-1,X,log,t}\ra K_{n,X,log,t}
\end{equation}
by applying $\uline p$ on each summand.
\ele
\bpf
It suffices to show that $C'_t:K_{n,X,t}\ra K_{n,X,t}$ commutes with $\uline p$ for both $t=\et$ and $t=\Zar$. 
For $t=\et$, $C'_\et$ is the composition of $\tau^{-1}:K_{n,X,\et}\ra (W_nF_X)_*K_{n,X,\et}$ and $\epsilon^*(C'_\Zar):(W_nF_X)_*K_{n,X,\et}\ra K_{n,X,\et}$. Since $\tau^{-1}$ is functorial with respect to any map of abelian sheaves, we know that
$$\xymatrix{
R_*K_{n-1,X,\et}\ar[r]^(.4){\tau^{-1}}\ar[d]^{\uline p}&
\ar[d]^{\uline p}
(W_nF_X)_*R_*K_{n-1,X,\et}\\
K_{n,X,\et}\ar[r]^(.4){\tau^{-1}}&
(W_nF_X)_*K_{n,X,\et}
}$$
is commutative, thus it suffices to prove the proposition for $t=\Zar$. 
That is, it suffices to show the diagrams (\ref{uline p log well def diag.1}) and (\ref{uline p log well def diag.2}) commute:
\beq\label{uline p log well def diag.1}
\xymatrix@C=8em{
R_*K_{n-1,X}
\ar[r]_(.44){\simeq}^(.44){R_*(\ref{(1.0.2)})}
\ar[d]^{\uline p}&
R_*(W_{n-1}F_X)^\triangle K_{{n-1},X}
\ar[d]^{\uline p_{\{(W_nF_X)^\triangle K_{n,X}\}_n}}
\\
K_{n,X}\ar[r]_(.4){\simeq}^(.4){(\ref{(1.0.2)})}&
(W_nF_X)^\triangle K_{n,X},
}
\eeq
\beq\label{uline p log well def diag.2}
\xymatrix@C=3em{
(W_{n}F_X)_*R_*(W_{n-1}F_X)^\triangle K_{{n-1},X}
\ar[r]^{\simeq}
\ar[d]^{(W_nF_X)_*\uline p_{\{(W_nF_X)^\triangle K_{n,X}\}_n}}&
R_*(W_{n-1}F_X)_*(W_{n-1}F_X)^\triangle K_{{n-1},X}
\ar[r]^(.65){R_*\Tr_{W_{n-1}F_X}}
&
R_*K_{{n-1},X}\ar[d]^{\uline p}
\\
(W_nF_X)_*(W_nF_X)^\triangle K_{n,X}
\ar[rr]^(.6){\Tr_{W_nF_X}}
&&
K_{n,X}.
}
\eeq
Here $\uline p:=\uline p_{\{K_{n,X}\}_n}$ is the lift-and-multiplication-by-$p$ map associated to the Witt residual complex $\{K_{n,X}\}_n$, while $\uline p_{\{(W_nF_X)^\triangle K_{n,X}\}_n}$ denotes the one associated to Witt residual system $\{(W_nF_X)^\triangle K_{n,X}\}_n$ (cf. \cite[1.8.7]{CR12}). By definition, the map
$$\uline p_{\{(W_nF_X)^\triangle K_{n,X}\}_n}: R_*(W_{n-1}F_X)^\triangle K_{n-1,X}\ra (W_{n}F_X)^\triangle K_{n,X}$$
is given by the adjunction map of
$$
(W_{n-1}F_X)^\triangle K_{n-1,X}
\xra[\simeq]{(W_{n-1}F_X)^\triangle (^a\uline p)}
(W_{n-1}F_X)^\triangle R^\triangle K_{n,X}
\simeq
R^\triangle (W_{n}F_X)^\triangle K_{n,X},
$$
where $^a\uline p$ is the adjunction of $\uline p$ for residual complexes (cf. \cite[Def. 1.8.3]{CR12}). The second diagram (\ref{uline p log well def diag.2}) commutes because the trace map $\Tr_{W_nF_X}$ induces a well-defined map between Witt residual complexes \cite[Lemma 1.8.9]{CR12}.

It remains to show the commutativity of (\ref{uline p log well def diag.1}). According to the definition of $\uline p_{(W_nF_X)^\triangle K_{n,X}}$ in \cite[1.8.7]{CR12}, we are reduced to show the adjunction square commutes:
$$\xymatrix@C=7em{
R^\triangle K_{n,X}
\ar[r]^(.44){R^\triangle(\ref{(1.0.2)})}&
R^\triangle(W_nF_X)^\triangle K_{n,X}
\ar[r]^(.44){\simeq} &
(W_{n-1}F_X)^\triangle R^\triangle K_{n,X}
\\
K_{n-1,X}\ar[u]^{^a\uline p}
\ar[rr]^{(\ref{(1.0.2)})}&&
(W_{n-1}F_X)^\triangle K_{n-1,X}
\ar[u]^{(W_{n-1}F_X)^\triangle(^a \uline p)}.
}$$
And this is $(W_{n-1}\pi)^\triangle$ applied to the following diagram
$$\xymatrix@C=7em{
R^\triangle W_nk
\ar[r]^(.44){R^\triangle(\ref{(1.0.1)})}&
R^\triangle(W_nF_k)^\triangle W_nk
\ar[r]^(.44){\simeq} &
(W_{n-1}F_k)^\triangle R^\triangle W_nk
\\
W_{n-1}k\ar[u]^{^a\uline p}
\ar[rr]^{(\ref{(1.0.1)})}&&
(W_{n-1}F_k)^\triangle W_{n-1}k
\ar[u]^{(W_{n-1}F_k)^\triangle(^a \uline p)}.
}$$
We are reduced to show its commutativity.
Notice that this diagram is over $\Spec W_{n-1}k$, where the only possible filtration is the one-element set consisting of the unique point of $\Spec W_{n-1}k$. 
This means that the Cousin functor associated to this filtration sends any dualizing complex to itself, and the map $^a\uline p$ in the sense of a map either between residual complexes \cite[Def. 1.8.3]{CR12} or between dualizing complexes \cite[Def. 1.6.3]{CR12} actually agree. 

Now we start the computation.
Formulas for (\ref{(1.0.1)}) and for $^a\uline p$ (in the sense of a map between dualizing complexes)  are explicitly given in \Cref{subsection on Def of KnXlog} and \cite[1.6.4(1)]{CR12}, respectively. 
Label the source and the target of $W_nF_k$ by $\Spec W_nk_1$ and $\Spec k_2$ respectively.
Take $a\in W_{n-1}k_1$. 
Denote $\bar{W_nF_k}:(\Spec W_nk_1,W_nk_1)\ra (\Spec W_nk_2,(W_nF_{k})_* (W_nk_1))$, and $\bar R: (\Spec W_{n-1}k_i, W_{n-1}k_i)\ra (\Spec W_{n}k_i, R_*W_{n-1}k_i)$ $(i=1,2)$ the natural maps of ringed spaces. 
Now the down-right composition $((W_{n-1}F_k)^\triangle(^a\uline p))\circ (\ref{(1.0.1)})$ equals to the Cousin functor $E_{(W_{n-1}F_k)^\triangle R^\triangle Z^\bullet(W_nk)}$ applied to the following composition
\begin{align*}
W_{n-1}k_1 &\xra[\simeq]{(\ref{(1.0.1)})}
\bar{W_{n-1}F_k}^*\Hom_{W_{n-1}k_2}((W_{n-1}F_{k})_*(W_{n-1}k_1),W_{n-1}k_2)
\\
&\xra[\simeq]{
^a\uline p}
\bar{W_{n-1}F_k}^*\Hom_{W_{n-1}k_2}((W_{n-1}F_{k})_*(W_{n-1}k_1),
\bar R^*\Hom_{W_nk_2}(R_*W_{n-1}k_2,W_{n}k_2)),
\\
a&\mapsto 
[(W_{n-1}F_{k})_*1\mapsto (W_{n-1}F_{k})^{-1}(a)]
\\&
\mapsto
[(W_{n-1}F_{k})_*1\mapsto
[R_*1\mapsto \uline p (W_{n-1}F_{k})^{-1}(a)]].
\end{align*}
And $R^\triangle(\ref{(1.0.1)})\circ (^a\uline p)$ equals to the Cousin functor $E_{(W_{n-1}F_k)^\triangle R^\triangle Z^\bullet(W_nk)}$ applied the following composition
\begin{align*}
W_{n-1}k_1 &\xra[\simeq]{^a\uline p}
\bar R^* \Hom_{W_nk_1}(R_*W_{n-1}k_1,W_nk_1)
\\
&\xra{
(\ref{(1.0.1)})}
\bar R^* \Hom_{W_nk_1}(R_*W_{n-1}k_1,
\bar{W_{n}F_k}^*\Hom_{W_{n}k_2}((W_{n}F_{k})_*(W_{n}k_1),W_{n}k_2)
),
\\
a&\mapsto
[R_*1\mapsto \uline p (a)]
\\
&\mapsto
[R_*1\mapsto [(W_nF_k)_*1\mapsto (W_nF_k)^{-1}\uline p (a)]].
\end{align*}
It remains to identify $\uline p ((W_{n-1}F_k)^{-1}a)$ and $(W_nF_k)^{-1}\uline p (a)$. And this is straightforward: write $a=\sum_{i=0}^{n-2} V^i[a_i] \in W_{n-1}k_1$,
\begin{align}\label{uline p commutes with WnFk^-1}
(W_nF_k)^{-1}\uline p (a)
&=\sum_{i=0}^{n-2}(W_nF_k)^{-1}\uline p (V^i[a_i])
=\sum_{i=0}^{n-2}(W_nF_k)^{-1}(V^{i+1}[a_i^p])\\
&
=\sum_{i=0}^{n-2}(V^{i+1}[a_i])
=\uline p \sum_{i=0}^{n-2} (V^i[a_i^{1/p}])
=\uline p ((W_{n-1}F_k)^{-1}a).
\nonumber
\end{align}
Hence we finish the proof.
\epf

However we don't naturally have a restriction map $R$ between residual complexes. Nevertheless, we can use the quasi-isomorphism $\bar \zeta_{n,X,log,t}:\tilde \nu_{n,X,t}\xra{\simeq} K_{n,X,log,t}$ to build up a map \begin{equation}\label{R for Kn,X,log}
R:K_{n,X,log,t}\ra K_{n-1,X,log,t}
\end{equation}
in the derived category $D^b(X,\Z/p^n)$. For this we will need to show that $\uline p$ and $R$ induce chain maps for $\tilde \nu_{n,X,t}$. This should be well-known to experts, we add here again due to a lack of reference.

\ble[]\label{uline p and R are chain maps for tilde v_nX}
$$\uline p: \tilde \nu_{n,X,t}\ra \tilde \nu_{n+1,X,t},\quad R: \tilde \nu_{n+1,X,t}\ra \tilde \nu_{n,X,t}$$
given by $\uline p$ and $R$ termwise, are well defined maps of complexes for both $t=\Zar$ and $t=\et$.
\ele

\bpf
It suffices to prove for $t=\Zar$.
Let $x\in X_{(q)}$ be a point of dimension $q$. Let $\rho:X'\ra \bar{\{x\}}$ be the normalization of $\bar{\{x\}}$. Let $x'$ be the generic point of $X'$ and $y'\in X'^{(1)}$ be a codimension $1$ point. Denote $y:=\rho(y')$.
It suffices to check the commutativity of the following diagrams in (1) and (2).
\begin{enumerate}
\item
Firstly,
$$
\xymatrix{
W_n\Omega_{x',log}^q\ar[r]^{\partial}\ar[d]^{\uline p}&
W_n\Omega_{y',log}^{q-1}\ar[d]^{\uline p}
\\
W_{n+1}\Omega_{x',log}^q\ar[r]^{\partial}&
W_{n+1}\Omega_{y',log}^{q-1},
}
\quad
\xymatrix{
W_{n+1}\Omega_{x',log}^q\ar[r]^{\partial}\ar[d]^{R}&
W_{n+1}\Omega_{y',log}^{q-1}\ar[d]^{R}
\\
W_{n}\Omega_{x',log}^q\ar[r]^{\partial}&
W_{n}\Omega_{y',log}^{q-1}.
}
$$
Notice that $p=\uline p\circ R$. Suppose $\pi'$ is a uniformizer of discrete valuation ring $\cO_{X',y'}$ and $u_1,\dots,u_q$ are invertible elements in $\cO_{X',y'}$. Calculate
\begin{align*}
\uline p(\partial (d\log [\pi']_n &d\log[u_2]_n\dots d\log [u_q]_n))\\
&=\uline p(d\log [u_2]_n\dots d\log [u_q]_n)\\
&=p(d\log [u_2]_{n+1}\dots d\log [u_q]_{n+1})\\
&=p(\partial(d\log [\pi']_{n+1} d\log[u_2]_{n+1}\dots d\log [u_q]_{n+1}))\\
&=\partial (p(d\log [\pi']_{n+1} d\log[u_2]_{n+1}\dots d\log [u_q]_{n+1}))\\
&=\partial(\uline p(d\log [\pi']_n d\log[u_2]_n\dots d\log [u_q]_n)),
\end{align*}
and
\begin{align*}
\uline p(\partial (d\log [u_1]_{n} &d\log[u_2]_{n}\dots d\log [u_q]_{n}))\\
&=0\\
&=p(\partial(d\log [u_1]_{n+1} d\log[u_2]_{n+1}\dots d\log [u_q]_{n+1}))\\
&=\partial(p(d\log [u_1]_{n+1} d\log[u_2]_{n+1}\dots d\log [u_q]_{n+1}))\\
&=\partial(\uline p(d\log [u_1]_{n} d\log[u_2]_{n}\dots d\log [u_q]_{n})).
\end{align*}
This proves the first diagram. Now the second.
\begin{align*}
R(\partial (d\log [\pi']_{n+1} &d\log[u_2]_{n+1}\dots d\log [u_q]_{n+1}))\\
&=R(d\log [u_2]_{n+1}\dots d\log [u_q]_{n+1})\\
&=d\log [u_2]_{n}\dots d\log [u_q]_{n}\\
&=\partial(d\log [\pi']_{n} d\log V[u_2]_{n}\dots d\log V [u_q]_{n})\\
&=\partial(R(d\log [\pi']_{n+1} d\log[u_2]_{n+1}\dots d\log [u_q]_{n+1})),
\end{align*}
and
\begin{align*}
R(\partial (d\log [u_1]_{n+1} &d\log[u_2]_{n+1}\dots d\log [u_q]_{n+1}))\\
&=0\\
&=\partial(d\log [u_1]_{n} d\log[u_2]_{n}\dots d\log [u_q]_{n})\\
&=\partial(R(d\log [u_1]_{n+1} d\log[u_2]_{n+1}\dots d\log [u_q]_{n+1})).
\end{align*}
\item
Secondly,
$$
\xymatrix{
W_n\Omega_{y',log}^{q-1}\ar[r]^{\tr}\ar[d]^{\uline p}&
W_n\Omega_{y,log}^{q-1}\ar[d]^{\uline p}
\\
W_{n+1}\Omega_{y',log}^{q-1}\ar[r]^{\tr}&
W_{n+1}\Omega_{y,log}^{q-1},
}
\quad
\xymatrix{
W_{n+1}\Omega_{y',log}^{q-1}\ar[r]^{\tr}\ar[d]^{R}&
W_{n+1}\Omega_{y,log}^{q-1}\ar[d]^{R}
\\
W_{n}\Omega_{y',log}^{q-1}\ar[r]^{\tr}&
W_{n}\Omega_{y,log}^{q-1}.
}
$$
Notice that $\rho:X'\ra \bar{\{x\}}^X$ can be restricted to a map from $\bar{\{y'\}}^{X'}$ to $\bar{\{y\}}^X$ ($\bar{\{x\}}^X$ denotes the closure of $x$ in $X$, and similarly for $\bar{\{y'\}}^{X'}$, $\bar{\{y\}}^X$). Furthermore, $y'$ (resp. $y$) belongs to the smooth locus of $\bar{\{y'\}}^{X'}$ (resp.  $\bar{\{y\}}^X$), and there $\uline p$ and $R$ come from the restriction of the $\uline p$ and $R$ on the respective smooth locus.
The map $\tr$, induced by Milnor's norm map, agrees with the Grothendieck trace map $\Tr_{W_n\rho}$ due to \Cref{Compatibility of Milnor norm and Grothendieck trace}. And according to compatibility of the Grothendieck trace map  with the Witt system structure (i.e. de Rham-Witt system structure with zero differential) on canonical sheaves \cite[4.1.4(6)]{CR12}, we arrive at the desired commutativity.
\end{enumerate}

\epf

\ble\label{uline p F ses}
Assume either
\bit
\item
$t=\Zar$ and $k=\bar k$, or
\item
$t=\et$ and $k$ being a perfect field of characteristic $p> 0$.
\eit
Then we have the following short exact sequence
\begin{equation}\label{pR ses for tilde vn,X}
0\ra \tilde \nu_{i,X,t}\xra{\uline p^j} \tilde \nu_{i+j,X,t}\xra{R^i} \tilde \nu_{j,X,t}\ra 0,
\end{equation}
in the category of complexes of sheaves over $X_t$, and a distinguished triangle
\begin{equation}\label{pR triangle for Kn,X,log}
K_{i,X,log,t}\xra{\uline p^{j}}
K_{i+j,X,log,t}\xra{R^i}
K_{j,X,log,t}\xra{+1}
\end{equation}
in the derived category $D^b(X_t,\Z/p^n)$.
\ele
\bpf
\begin{enumerate}
\item
Because of \Cref{uline p and R are chain maps for tilde v_nX}, it suffices to show
$$0\ra W_i\Omega^q_{x,log,t}\xra{\uline p^j} W_{i+j}\Omega^q_{x,log,t}\xra{R^i} W_j\Omega^q_{x,log,t}\ra 0$$
is short exact for any given point $x\in X_{(q)}$. And this is true for $t=\et$ because of \cite[Lemme 3]{CTSS-TorDanGroupeDeChow}. And for $t=\Zar$, one further needs $R^1\epsilon_*W_n\Omega_{x,log,\et}^q=0$ for any $x\in X_{(q)}$ if $k=\bar k$, which is proved in \cite[Cor. 2.3]{Suwa-NoteOnGerstenForlogSheaves}.
\item
Now it suffices to show that $\uline p$ and $R$ for the system $\{K_{n,X,log,t}\}_n$ are compatible with $\uline p$ and $R$ of the system $\{\tilde \nu_{n,X,t}\}_n$, via the quasi-isomorphism $\bar \zeta_{n,X,log,t}$. The compatibility for $R$ is clear by definition. It remains to check the compatibility for $\uline p$. Because $\bar\zeta_{n,X,log,t}=(\ref{K to Klog map of graded algebras})\circ \bar\zeta_{n,X,t}$, it suffices to check compatibility of $\uline p:\tilde \nu_{n-1,X,t}\ra \tilde \nu_{n,X,t}$ with $\uline p: K_{n-1,X,t}\ra K_{n,X,t}$ via $\bar\zeta_{n,X,t}$. At a given degree $-q$ and given point $x\in X_{(q)}$, the map $\bar\zeta_{n,X,t}:\tilde \nu_{n,X,t}\ra K_{n,X,t}$ factors as
$$(W_n\iota_x)_*W_n\Omega_{x,log,t}^q\xra{} (W_ni_x)_*W_n\Omega_{x,t}^q
=(W_ni_x)_*K^q_{n,\bar{\{x\}},t}\xra{(-1)^q\Tr_{W_ni_x,t}} K^q_{n,X,t}.$$
The first arrow is the inclusion map and is naturally compatible with $\uline p$. The compatibility of $\uline p$ via the trace map is given in \cite[Lemma 1.8.9]{CR12}. 
\end{enumerate}
\epf

\section{Higher Chow groups of zero cycles}
\label{section Higher Chow groups of zero cycles}
Let $k$ be a perfect field of characteristic $p>0$. In the whole \Cref{section Higher Chow groups of zero cycles}, $X$ denotes a separated scheme of finite type over $k$ of dimension $d$ with structure map $\pi:X\ra k$.
\subsection{First properties}
\bprop\label{triangle prop ZcX KnX KnX}
There is a distinguished triangle
$$\Z^c_{X,\et}/p^n\ra K_{n,X,\et}\xra{C'_\et-1} K_{n,X,\et}\xra{+1}$$
in the derived category $D^b(X_\et,\Z/p^n)$.
If $k=\bar k$, one also has the Zariski counterpart. Namely, we have a distinguished triangle
\beq\label{triangle prop ZcX KnX KnX: zar}
\Z^c_{X}/p^n\ra K_{n,X}\xra{C'-1} K_{n,X}\xra{+1}
\eeq
in the derived category $D^b(X,\Z/p^n)$.

In particular, if $k=\bar k$ and $X$ is Cohen-Macaulay of pure dimension $d$, then $\Z^c_{X}/p^n$ is concentrated at degree $-d$, and the triangle (\ref{triangle prop ZcX KnX KnX: zar}) becomes
$$\Z^c_{X}/p^n\ra W_n\omega_X[d]\xra{C'-1} W_n\omega_X[d]\xra{+1}$$
in this case. Here $W_n\omega_X$ is the only non-vanishing cohomology sheaf of $K_{n,X}$ (if $n=1$, $W_1\omega_X=\omega_X$ is the usual dualizing sheaf on $X$).
\eprop
\bpf
This is direct from the main result \Cref{Main theorem 2} and \Cref{surjectivities 1-C C^-1 -1 etc CM}.
\epf
\bprop\label{CH_0=C'-1 invariant}
Assume $k=\bar k$. Then higher Chow groups of zero cycles equal the $C'$-invariant part of the cohomology groups of Grothendieck's coherent dualizing complex, i.e.,
$$\CH_0(X,q;\Z/p^n)
=H^{-q}(W_nX,K_{n,X})^{C'-1}.$$
\eprop
\bpf
This follows directly from the \Cref{Raynaud-Illusie}
and the main result \Cref{Main theorem 2}.
\epf

\bcor[Relation with $p$-torsion Poincar\'e duality]\label{combine with JSS} 
There is an isomorphism 
$$K_{n,X,log,\et} \simeq R\pi^!(\Z/p^n)$$
in $D^b(X_\et, \Z/p^n)$,
where $R\pi^!$ is the extraordinary inverse image functor defined in \cite[Expos\'e XVIII, Thm 3.1.4]{SGA4-3}.
\ecor
\bpf 
This follows directly from the main theorem \Cref{Main theorem} and \cite[Thm. 4.6.2]{JSS}.
\epf

\bcor[Affine vanishing]
\label{Vanishing higher Chow}
Suppose $X$ is affine and Cohen-Macaulay of pure dimension $d$. Then
\benu
\item
If $t=\Zar$ and $k=\bar k$, $$\CH_0(X,q,\Z/p^n)
=0$$ for $q\neq d$.
\item
If $t=\et$, $$R^{-q}\Gamma(X_\et,\Z^c_{X}/p^n)
=0$$
for $q\neq d,d-1$. If one further assumes $k=\bar k$ or smoothness, the possible non-vanishing occurs only in degree $q=d$.
\eenu
\ecor
\bpf
If $X$ is Cohen-Macaulay of pure dimension $d$, $W_nX$ is also Cohen-Macaulay of pure dimension $d$ by Serre's $S_k$-criterion, and $K_{n,X,t}$ is concentrated at degree $-d$ for all $n$ \cite[3.5.1]{Conrad-GDBC}. 
Now Serre's affine vanishing theorem implies $H^{-q}(W_nX,K_{n,X,t})=0$ for $q\neq d$. This implies that $R^{-q}\Gamma(W_nX,K_{n,X,log,t})=0$ unless $q=d,d-1$. With the given assumptions, \Cref{Main theorem 2} implies that $\CH_0(X,q,\Z/p^n)=R^{-q}\Gamma(X_\et,\Z^c_{X}/p^n)=0$ unless $q=d,d-1$.
If one also assumes $k=\bar k$, \Cref{CH_0=C'-1 invariant} gives the vanishing result for $q=d-1$.

If $X$ is smooth, $C_\et-1:W_n\Omega_{X,\et}^d\ra W_n\Omega_{X,\et}^d$ is surjective by \cite[1.6(ii)]{GrosSuwa-AbelJacobi} (see (\ref{C-1 ses})). By the compatibility of $C_\et$ and $C'_\et$ \Cref{Compatibility of $C'_n$ and $C$ etale}, one deduces that $C'-1:\cH^{-d}(K_{n,X,\et})\break\ra \cH^{-d}(K_{n,X,\et})$ is surjective.
\epf

Generalizing Bass's finiteness conjecture for $K$-groups (cf. \cite[IV.6.8]{Weibel-KBook}), the finiteness of higher Chow groups in various arithmetic settings appears in the literature. 
The following result was first proved by Geisser \cite[\S5, eq. (12)]{Geisser} using the finiteness result from the \'etale cohomology theory, and here we deduce it as a corollary of our main theorem, which essentially relies on the finiteness of coherent cohomologies on a proper scheme. We remark that Geisser's result is more general than ours in that he allows arbitrary torsion coefficients.

\bcor[Finiteness, Geisser]\label{Finiteness of higher Chow}
Assume $k=\bar k$. Let $X$ be proper over $k$. Then for any $q$,
$$\CH_{0}(X,q;\Z/p^n)\quad 
$$
is a finite $\Z/p^n$-module.
\ecor
\bpf
According to \Cref{Main theorem 2}, $R^{-q}\Gamma(X,\Z^c_{X}/p^n)
=R^{-q}\Gamma(X,K_{n,X,log})$. 
Thus it suffices to show that for every $q$, $R^{-q}\Gamma(X,K_{n,X,log})$ is a finite $\Z/p^n$-module. 
First of all, since the cohomology group $R^{-q}\Gamma(X,K_{n,X,log})$ is the $C'$-invariant part of $R^{-q}\Gamma(X,K_{n,X})$ by \Cref{Raynaud-Illusie} and \Cref{Raynaud-Illusie etale},  $R^{-q}\Gamma(X,K_{n,X,log})$ is a module over the invariant ring $(W_nk)^{1-W_nF_X^{-1}}=\Z/p^n$. 
Because $X$ is proper, $R^{-q}\Gamma(X,K_{n,X})$ is a finite $W_nk$-module by the local-to-global spectral sequence. Then \Cref{1-T invariant is a finite module} gives us the result.

Alternatively, we can also do induction on $n$.
In the $n=1$ case, because $R^{-q}\Gamma(X,K_{X,log})$ is the $C'$-invariant part of the finite dimensional $k$-vector space $H^{-q}(X,K_{X})$ again by \Cref{Raynaud-Illusie} and \Cref{Raynaud-Illusie etale}, it is a finite $\F_p$-module by $p^{-1}$-linear algebra \Cref{1-T surjective prop in appendix}. 
The desired result then follows from the long exact sequence associated to (\ref{pR triangle for Kn,X,log}) by induction on $n$.
\epf

We refer to \Cref{def of sigma-linear generalized} and \Cref{def of sigma-linear generalized rmk}(2) for the definition of the semisimplicity and the notation $(-)_{\mathrm{ss}}$ in this context.
\bcor[Semisimplicity]
\label{semisimplicity cor}
Assume $k=\bar k$. Let $X$ be proper over $k$. Then for any $q$,
$$H^{-q}(W_nX,K_{n,X})_{\mathrm{ss}}=\CH_0(X,q;\Z/p^n)\otimes_{\Z/p^n} W_nk.$$
\ecor 
\bpf 
Since $X$ is proper, $H^{-q}(W_nX,K_{n,X})$ is a finite $W_nk$-module for any $q$.
Then according to \Cref{Mss=M1-T otimes Wnk},
$$H^{-q}(W_nX,K_{n,X})_{\mathrm{ss}}=H^{-q}(W_nX,K_{n,X})^{C'-1}\otimes_{\Z/p^n} W_nk.$$
The claim now follows from \Cref{CH_0=C'-1 invariant}.
\epf

\subsection{\'Etale descent}
The results \Cref{etale descent tilde vnX}, \Cref{etale descent ZcX} in this subsection are well-known to experts.
\bprop[Gros-Suwa]\label{etale descent tilde vnX}
Assume $k=\bar k$. Then one has a canonical isomorphism
$$\tilde \nu_{n,X,\Zar}=\epsilon_*\tilde \nu_{n,X,\et}
\xra{\simeq}
R\epsilon_*\tilde \nu_{n,X,\et}$$
in the derived category $D^b(X,\Z/p^n)$.
\eprop
\bpf
If $k=\bar k$, terms of the \'etale complex $\tilde \nu_{n,X,\et}$ are $\epsilon_*$-acyclic according to \cite[3.16]{GrosSuwa-AbelJacobi}.
\epf

The \'etale descent of Bloch's cycle complex with $\Z$-coefficients is shown in \cite[Thm 3.1]{Geisser}, assuming the Beilinson-Lichtenbaum conjecture which is now proved by Rost and Voevodsky. Hence the \'etale descent holds conjecture-free. 
Note that one can also deduce the mod $p^n$ version as a corollary of \Cref{etale descent tilde vnX} via Zhong's quasi-isomorphism in \Cref{subsection Zhong's qis combined with zeta} (which is dependent on the main result of Geisser-Levine \cite[1.1]{GeisserLevine}). 
\bprop[Geisser-Levine]\label{etale descent ZcX}
Assume $k=\bar k$. Then one has a canonical isomorphism
$$\Z^c_{X,\Zar}/p^n=\epsilon_*\Z^c_{X,\et}/p^n
\xra{\simeq}R\epsilon_*\Z^c_{X,\et}/p^n.$$
in the derived category $D^b(X,\Z/p^n)$.
As a result,
$$\CH_0(X,q;\Z/p^n)\simeq R^{-q}\Gamma(X_\et, \Z^c_{X,\et}/p^n).$$
\eprop
\bpf
Clearly, we have the compatibility
$$\xymatrix{
\Z^c_{X,\Zar}/p^n\ar[d]^{\psi_{X,\Zar}}
\ar[r]^{\simeq}
&R\epsilon_*\Z^c_{X,\et}/p^n
\ar[d]^{R\epsilon_*\psi_{X,\et}}
\\
\tilde \nu_{n,X,\Zar}
\ar[r]^{\simeq}&
R\epsilon_*\tilde \nu_{n,X,\et}.
}$$
Thus
$\Z^c_{X,\Zar}/p^n\xra[\simeq]{\bar\psi_{X,\Zar}}
\tilde \nu_{n,X,\Zar}=\epsilon_*\tilde \nu_{n,X,\et}
\xra[\simeq]{\text{\Cref{etale descent tilde vnX}}}
R\epsilon_*\tilde \nu_{n,X,\et}
\xra[\simeq]{R\epsilon_*\bar\psi_{X,\et}}
R\epsilon_*\Z^c_{X,\et}/p^n.$
\epf

\bcor\label{etale descent vanishing?}
Assume $k=\bar k$. Suppose $X$ is affine and Cohen-Macaulay of pure dimension $d$. Then
$$R^i\epsilon_*(\Z^c_{X,\et}/p^n)
=R^i\epsilon_*\tilde \nu_{n,X,\et}=0,\quad i\neq -d.$$
\ecor
\bpf
This is a direct consequence of \Cref{etale descent ZcX}, \Cref{etale descent tilde vnX} and \Cref{Vanishing higher Chow}.
\epf

\subsection{Birational geometry and rational singularities}\label{subsection Birational geometry and rational singularities}

Recall the following definition of \textit{resolution-rational singularities}, which are more often called \textit{rational singularities} before in the literature, but here we follow the terminology from \cite{Kovacs-RatSing} (see also \Cref{remark RatSing}(1) and \cite[(9.12.1)]{Kovacs-RatSing}).

\bdefn[{cf. \cite[9.1]{Kovacs-RatSing}}]
\label{def RatSing}
An integral $k$-scheme $X$ is said to have \emph{resolution-rational singularities}, if
\benu
\item
there exists a birational proper morphism $f:\tilde X\ra X$ with $\tilde X$ smooth (such a $f$ is called a \emph{resolution of singularities} or simply a \emph{resolution} of $X$), and
\item\label{cohomological condition for rational singularities}
$R^if_*\cO_{\tilde X}=R^if_*\omega_{\tilde X}=0$ for $i\ge 1$. And $f_*\cO_{\tilde X}=\cO_X$.
\eenu
Such a map $f:\tilde X\ra X$ is called a \emph{rational resolution} of $X$. 
\edefn 
Note that the cohomological condition (\ref{cohomological condition for rational singularities}) is equivalent to the following condition
\bit
\item[(2')]
$\cO_X
\simeq Rf_*\cO_{\tilde X}$, $f_*\omega_{\tilde X}\simeq Rf_*\omega_{\tilde X}$ in the derived category of abelian Zariski sheaves.
\eit

\brmk 
\label{remark RatSing}
\benu 
\item 
According to \cite[8.2]{Kovacs-RatSing}, on integral $k$-schemes of pure dimension, our definitions for resolution-rational singularities and for rational resolutions  are the same as the ones in \cite[9.1]{Kovacs-RatSing}. 
\item 
Necessary conditions for an integral $k$-scheme to have resolution-rational singularities,
are that the scheme is normal and Cohen-Macaulay. 
The normality statement follows from the equality $f_*\cO_{\tilde{X}}=\cO_X$, and the Cohen-Macaulay statement is a standard result, see e.g. \cite[8.3]{Kovacs-RatSing}. 
\item 
According to \cite[9.6]{Kovacs-RatSing}, resolution-rational singularities are pseudo-rational. By definition \cite[1.2]{Kovacs-RatSing}, a $k$-scheme $X$ is said to have \emph{pseudo-rational singularities}, if it is normal Cohen-Macaulay, and for every normal scheme $X'$, every projective birational morphism $f:X'\ra X$, the composition $f_*\omega_{X'}\ra Rf_*\omega_{X'}\xra{\Tr_f} \omega_X$ is an isomorphism.
\eenu 
\ermk

\bcor\label{invariance under rational resolution}
Let $X$ and $Y$ be integral $k$-schemes of pure dimensions which have pseudo-rational singularities.
Suppose there are proper birational $k$-morphisms $f:Z\ra X$ and $g:Z\ra Y$ where $Z$ is a normal Cohen-Macaulay scheme. Then we have
$$R^{-q}\Gamma(X_\et,\Z^c_{\tilde X,\et}/p^n)= R^{-q}\Gamma( Y_\et,\Z^c_{Y,\et}/p^n)
$$
for all $q$ and all $n\ge 1$. 
If we assume furthermore $k=\bar k$, we also have 
$$\CH_0(X, q, \Z/p^n) = \CH_0(Y, q, \Z/p^n)$$
for all $q$ and all $n\ge 1$. 
\ecor
\brmk\label{invariance under rational resolution rmk} 
\benu 
\item 
In particular, since for any rational resolution of singularities $f:\tilde X\ra X$, $\tilde X$ and $X$ are properly birational as $k$-schemes (i.e., take $Z$ to be $\tilde X$), one can compute the higher Chow groups of zero cycles of $X$ via those of $\tilde X$.
\item 
Deleting the pseudo-rational singularities assumption (in particular, we relax the Cohen-Macaulay assumptions on $X$ and $Y$), the proof still passes through with the following assumption: $X$ and $Y$ are linked by a chain of proper birational maps and each of these maps has its trace map being a quasi-isomorphism between the residual complexes. Such a proper birational map is called a \textit{cohomological equivalence} in \cite[8.4]{Kovacs-RatSing}.
\item 
If normal Macaulayfications for integral varieties exist (e.g., conjecture \cite[1.13]{Kovacs-RatSing} is true), the assumption of \Cref{invariance under rational resolution} can be weakened to the following: \textit{Let $X$ and $Y$ be integral $k$-schemes of pure dimensions which have pseudo-rational singularities and are properly birational, i.e., there are proper birational $k$-morphisms $f:Z\ra X$ and $g:Z\ra Y$ with $Z$ being some integral scheme.}
In fact, we can replace $Z$ by a normal Cohen-Macaulay scheme by the following process.
Using Chow's Lemma \cite[4.1]{Kovacs-RatSing}, we know that there exist projective birational morphisms $f':Z_1\ra Z$ and $g':Z_2\ra Z$ such that the compositions $Z_1\xra{f'} Z\xra{f}X$ and $Z_2\xra{g'} Z\xra{g} Y$ are also birational and projective. Let $U\subset Z$ be an open dense subset such that $f'$ and $g'$ restricted to the preimage of $U$ are isomorphisms. Take $Z'$ be the Zariski closure of the image of the diagonal of $U$ in $Z_1\times_Z Z_2$ with the reduced scheme structure. Then the two projections $Z'\ra Z_1$ and $Z'\ra Z_2$ are also projective and birational.
This means that by replacing $Z'$ with $Z$, $f$ with $Z'\ra X$ and $g$ with $Z'\ra Y$, we can assume our $f:Z\ra X$, $g:Z\ra Y$ to be projective birational and our $Z$ to be integral.
Using normal Macaulayfication \cite[1.13]{Kovacs-RatSing} we can additionally assume that $Z$ is normal Cohen-Macaulay. 

In particular, since the conjecture \cite[1.13]{Kovacs-RatSing} is known to be true for varieties of dimension at most 4 over algebraically closed fields (cf. \cite[1.14(iii)]{Kovacs-RatSing}), one can state \Cref{invariance under rational resolution} with this weakened assumption in this case.
\eenu 
\ermk 
\bpf
Note that $f$ and $g$ are pseudo-rational modifications by \cite[9.7]{Kovacs-RatSing}.
Suppose that $X$ is of pure dimension $d$. Then so is $Z$. Now \cite[8.6]{Kovacs-RatSing} implies that the trace map of $f$ induces an isomorphism
$$\Tr_f:Rf_* K_{Z,t} \xra{\simeq}K_{X,t}$$
in $D^b(X_t,\Z/p)$. Thus
$$\Tr_{f,log}: Rf_*K_{Z,log,t}\xra{\simeq}K_{X,log,t}$$
is also an isomorphism in $D^b(X_t,\Z/p)$. Consider the diagram
\beq\label{invariance under rational resolution:diag1}
\xymatrix{
f_*K_{Z,log,t}\ar[r]^{\uline p^{n-1}}\ar[d]^{\Tr_{f,log}}&
f_*K_{n,Z,log,t}\ar[r]^{R}\ar[d]^{\Tr_{W_nf,log}}&
f_*K_{n-1,Z,log,t}\ar[r]^{+1}\ar[d]^{\Tr_{W_{n-1}f,log}}&
f_*K_{Z,log,t}[1]\ar[d]^{\Tr_{f,log}[1]}
\\
K_{X,log,t}\ar[r]^{\uline p^{n-1}}&
K_{n,X,log,t}\ar[r]^{R}&
K_{n-1,X,log,t}\ar[r]^{+1}&
K_{X,log,t}[1]
}
\eeq
in $D^b(X_t,\Z/p)$. The first row is $Rf_*$ applied to the triangle (\ref{pR triangle for Kn,X,log}) on $Z$. The second row is the triangle (\ref{pR triangle for Kn,X,log}) on $X$. The left square commutes on the level of complexes by the compatibility of the trace map with $\uline p$ \cite[1.8.9]{CR12}. To prove commutativity of the middle square in the derived category, it suffices to show the square
$$\xymatrix{
f_*\tilde\nu_{n,Z,t}\ar[r]^R\ar[d]^{f_*}&
f_*\tilde\nu_{n-1,Z,t}\ar[d]^{f_*}
\\
\tilde\nu_{n,X,t}\ar[r]^R&
\tilde\nu_{n-1,X,t}
}$$
commutes on the level of complexes. Since the vertical maps $f_*$ for Kato-Moser complexes are $\tr$ 
(cf. \S\ref{section Gersten complex of logarithmic de Rham-Witt sheaves}), which are by definition the reduction of the norm maps for Milnor $K$-theory, they agree with the Grothendieck trace maps $\Tr_{W_nf}$, $\Tr_{W_{n-1}f}$ by \Cref{Compatibility of Milnor norm and Grothendieck trace}. 
And according to the compatibility of $R$ with the Grothendieck trace maps \cite[4.1.4(6)]{CR12}, we arrive at the desired commutativity.
The right square in (\ref{invariance under rational resolution:diag1}) commutes by the naturality of the "$+1$" map.
With all these commutativities we conclude that the vertical maps in (\ref{invariance under rational resolution:diag1}) define a map of triangles. By induction on $n$ we deduce that
$$\Tr_{W_nf,log}: Rf_*K_{n,Z,log,t} \xra{\simeq}K_{n,X,log,t}$$
is an isomorphism in $D^b(X_t,\Z/p^n)$ for every $n$. The main result \Cref{Main theorem 2} thus implies
$$R^{-q}\Gamma(Z_\et,\Z^c_{\tilde X,\et}/p^n)= R^{-q}\Gamma( X_\et,\Z^c_{X,\et}/p^n)$$
for all $q$ and $n$. If $k=\bar k$, the same theorem also implies that
$$\CH_0(Z, q, \Z/p^n) = \CH_0(X, q, \Z/p^n)$$
for all $q$ and $n$.

Now replacing $f$ by $g$ everywhere in the above argument and we get the result.
\epf

\subsection{Galois descent}
\bcor\label{Galois descent cor}
Let $X$ and $Y$ be separated schemes of finite type over $k$ of dimension $d$.
Let $f:Y\ra X$ be a finite \'etale Galois map with Galois group $G$. Then $$R^{-d}\Gamma(X_{\et},\Z^c_X/p^n)=R^{-d}\Gamma(Y_{\et},\Z^c_Y/p^n)^G.$$
If $k=\bar k$, we also have
$$\CH_0(X,d;\Z/p^n) = \CH_0(Y,d;\Z/p^n)^G.$$
\ecor 
\bpf 
The pullback $f^*$ induces two canonical maps
$$f^*:\Z^c_{X,\et}\ra (f_*\Z^c_{Y,\et})^G, \quad f^*: K_{n,X,log,\et}\ra (f_*K_{n,Y,log,\et})^G.$$
Both of them are isomorphisms of complexes, because each term of these complexes is an \'etale sheaf.
Because of the contravariant functoriality with respect to \'etale morphisms (\Cref{Functoriality cycle to KnXlog} and \cite[2.15]{Zhong14}), $\bar\zeta_{log}\circ\bar\psi$ is $G$-equivariant. That is, the diagram
$$\xymatrix@C=4em{
\Z^c_{X,\et}/p^n\ar[r]^{\bar\zeta_{log}\circ\bar\psi}\ar[d]^{f^*}&
K_{n,X,log,\et}\ar[d]^{f^*}\\
(f_*\Z^c_{Y,\et}/p^n)^G\ar[r]^{\bar\zeta_{log}\circ\bar\psi}&
(f_*K_{n,Y,log,\et})^G
}$$
commutes.

Apply $R^{-d}\Gamma(X_{\et},-)$ to the isomorphism $f^*:K_{n,X,log,\et}\ra (f_*K_{n,Y,log,\et})^G$, one gets
$$R^{-d}\Gamma(X_{\et},K_{n,X,log,\et})=R^{-d}\Gamma(X_{\et},(f_*K_{n,Y,log,\et})^G).$$
Consider the local-to-global spectral sequence associated to the right hand side of this equality, there is only one non-zero term in the $E_\infty$-page with total degree $-d$ (which is a term in the $E_2$-page), thus we have 
$$R^{-d}\Gamma(X_{\et},(f_*K_{n,Y,log,\et})^G)=H^0(X_{\et},\cH^{-d}((f_*K_{n,Y,log,\et})^G)).$$
Because $(-)^{G}$ commutes with taking kernels and with $H^0$, we have 
$$H^0(X_{\et},\cH^{-d}((f_*K_{n,Y,log,\et})^G))=H^0(X_{\et},\cH^{-d}(f_*K_{n,Y,log,\et}))^G.$$
Because $f_*$ preserves kernels, we have 
$$H^0(X_{\et},\cH^{-d}(f_*K_{n,Y,log,\et}))^G=H^0(Y_\et, \cH^{-d}(K_{n,Y,log,\et}))^G.$$
Again by the observation from the spectral sequence, this means
$$H^0(Y_\et, \cH^{-d}(K_{n,Y,log,\et}))^G=R^{-d}\Gamma(Y_{\et},K_{n,Y,log,\et})^G.$$
Since $\bar\zeta_{log}\circ\bar\psi$ is $G$-equivariant, the main theorem \Cref{Main theorem 2} implies
$$R^{-d}\Gamma(X_{\et},\Z^c_X/p^n)=R^{-d}\Gamma(Y_{\et},\Z^c_Y/p^n)^G.$$
If $k=\bar k$, \Cref{etale descent ZcX} implies
$$\CH_0(X,d;\Z/p^n) = \CH_0(Y,d;\Z/p^n)^G.$$
\epf

\part*{Appendix}
\begin{appendices}
\section{Semilinear algebra}\label{appendix: sigma-linear algebra}

\bdefn\label{def of sigma-linear}
Let $k$ be a perfect field of positive characteristic $p$, and $V$ be a finite dimensional $k$-vector space. A \emph{$p$-linear map} (resp. $p^{-1}$-linear map) on $V$ is a map $T : V \ra V$, such that
\begin{align*}
&T(v + w) = T(v) + T(w), \quad T(cv) = c^pT(v),\quad v,w\in V, c\in k.\\
\text{(resp. }
&T(v + w) = T(v) + T(w), \quad T(cv) = c^{-p}T(v),\quad v,w\in V, c\in k.
\text{)}
\end{align*}
We say a map $T:V\ra V$ is \textit{semilinear} if it is either $p$-linear or $p^{-1}$-linear. A semilinear map $T:V\ra V$ is called \emph{semisimple}, if $\oIm T=V$. 
\edefn

\brmk
Let $T$ be a semilinear map.
\begin{enumerate}
\item
Note that
$$\{c \in k \mid c^p = c\}=\F_p= \{c \in k \mid c^{-p} = c\}.$$
The fixed points of $T$
$$V^{1-T}:=\{v\in V\mid T(v)=v\}$$
is a $\F_p$-vector space. 
\item
There is a descending chain of $k$-vector subspaces of $V$
$$\oIm T\supset \oIm T^2\supset\dots\supset \oIm T^n\supset \dots.$$
Since $V$ is finite dimensional, it becomes stationary for some large $N\in \N$.
Define
$$V_{\mathrm{ss}}:=\bigcap_{n\ge 1} \oIm (T^n)= \oIm (T^N)= \oIm (T^{N+1})=\dots. $$
Obviously,
\begin{enumerate}
\item
$V_{\mathrm{ss}}$ is a $k$-vector subspace of $V$ that is stable under $T$.
$T$ is semisimple on $V_{\mathrm{ss}}$.
\item
$V^{1-T}\subset V_{\mathrm{ss}}.$
\end{enumerate}
\end{enumerate}
\ermk

The proof of the following result is given in \cite{SGA7II} for $p$-linear maps, but an analogous proof also works for $p^{-1}$-linear maps.

\bprop[{\cite[Expos\'{e} XXII, Cor. 1.1.10, Prop. 1.2]{SGA7II}}]\label{1-T surjective prop in appendix}
Suppose $k$ is a separably closed field of positive characteristic $p$. 
Then
$$1-T:V\ra V$$
is surjective. And
$$V_{\mathrm{ss}}\simeq V^{1-T}\otimes_{\F_p} k,$$
which in particular means $V^{1-T}$ is a finite dimensional $\F_p$-vector space with $\dim_{\F_p} V^{1-T}=\dim_k V_{\mathrm{ss}}$.
\eprop

We generalize the definition of a semilinear map.
\bdefn\label{def of sigma-linear generalized}
Let $k$ be a perfect field of positive characteristic $p$, and 
let $W_nk$ be the ring of the $n$-th truncated Witt vectors of $k$. Let $M$ be a finitely generated $W_nk$-module. 
A \emph{$p$-linear map} (resp. \emph{$p^{-1}$-linear map}) on $M$ is a map $T : M \ra M$, such that
\begin{align*}
&T(m + m') = T(m) + T(m'), \quad T(cm) = W_nF_k(c)T(m),\quad
m,m'\in M, c\in W_nk.\\
\text{(resp. }
&T(m + m') = T(m) + T(m'), \quad T(cm) = W_nF_k^{-1}(c)T(m),\quad m,m'\in M, c\in W_nk.
\text{)}
\end{align*}
Here $F_k$ denotes the $p$-th power Frobenius on the field $k$. We say that $T$ is \textit{semilinear} if it is either $p$-linear or $p^{-1}$-linear in this sense.
 A semilinear map $T:V\ra V$ is called \emph{semisimple}, if $\oIm T=V$. 
\edefn

\brmk 
\label{def of sigma-linear generalized rmk}
Let $T$ be a semilinear map in the sense of \Cref{def of sigma-linear generalized}.
\benu 
\item 
Write $\sigma=W_nF_k$ (resp. $\sigma=W_nF_k^{-1}$). Then
$$(W_nk)^{1-\sigma}:= \{c \in W_nk \mid \sigma(c) = c\}=\Z/p^n$$
for both cases.
The fixed points of $T$
$$M^{1-T}:=\{m\in M\mid T(m)=m\}$$
is a $\Z/p^n$-module. 
\item 
As in the case of vector spaces,
$$\oIm T\supset \oIm T^2\supset\dots\supset \oIm T^n\supset \dots$$
is a descending chain of $W_nk$-submodules of $M$. It becomes stationary for some large $N\in\N$, because $M$ as a finitely generated $W_nk$-module is artinian.
Define the $W_nk$-submodule of $M$ 
$$M_{\mathrm{ss}}:=\bigcap_{n\ge 1} \mathrm{Im}(T^n)= \mathrm{Im}(T^N)= \mathrm{Im}(T^{N+1})=\dots.$$
Then
\benu 
\item 
$M_{\mathrm{ss}}$ is a $W_nk$-submodule of $M$ that is stable under $T$. $T$ is semisimple on $M_{\mathrm{ss}}$.
\item 
$M^{1-T}\subset M_{\mathrm{ss}}$.
\item 
$(M/p)_{\mathrm{ss}}=M_{\mathrm{ss}}/p\subset M/p$.
\eenu
\eenu 
\ermk

\bprop\label{1-T surjective prop in appendix for W_n}
Let $k$ be a separably closed field of positive characteristic $p$.
Then
$$1-T:M\ra M$$ is surjective.
\eprop
\bpf
Take $m\in M$. Because $M$ is finitely generated as a $W_nk$-module, $M/pM$ is a finite dimensional $k$-vector space. Then \Cref{1-T surjective prop in appendix} implies that there exists a $m'\in M$, such that
$(1-T)(m')-m\in pM$. That is, there exists a $m_1\in M$ such that
$$(1-T)(m')=m+pm_1.$$
Do the same process with $m_1$ instead of $m$, one gets a $m'_1\in M$ and a $m_2\in M$ such that
$$(1-T)(m'_1)=m_1+pm_2.$$
Thus
$$(1-T)(m'-pm'_1)=m-p^2m_2.$$
Repeat this process. After finitely many times, because $p^n=0$ in $W_nk$,
$$(1-T)(m'-pm'_1+\dots+(-1)^{n-1}p^{n-1}m'_{n-1})=m.$$
\epf

\bprop\label{1-T invariant is a finite module}
Let $k$ be a separably closed field of positive characteristic $p$.
Then 
\benu 
\item 
$M^{1-T}/(pM)^{1-T}= (M/p)^{1-T}.$
\item 
$M^{1-T}$ is a finite $\Z/p^n$-module.
\eenu 
\eprop

\bpf
Since $W_nk$ is of $p^n$-torsion, we know that $p^mM=0$ for some $m\le n$. Do induction on the smallest number $m$ such that $p^mM=0$. If $m=1$, the first claim is trivial, and $M=M/p$ is actually a finite dimensional $k$-vector space, thus the second claim follows from \Cref{1-T surjective prop in appendix}. 

Now we assume $m>1$.
Note that $T$ induces a semilinear map on $pM$ and $pM$ is a finite $W_nk$-module, so by \Cref{1-T surjective prop in appendix for W_n} the map $1-T: pM\ra pM$ is surjective.
Now we have the two rows on the bottom of the following diagram being exact:
$$\xymatrix{
&0&0&0&
\\
0\ar[r]&  M^{1-T}/(pM)^{1-T} \ar[u]\ar[r]& M/p \ar[u]\ar[r]^{1-T}& M/p \ar[u]\ar[r]& 0
\\
0\ar[r]&  M^{1-T}\ar[u] \ar[r]& M \ar[u]\ar[r]^{1-T}& M \ar[u]\ar[r]& 0
\\
0\ar[r]&  
(pM)^{1-T}\ar[u]\ar[r]& pM\ar[u]\ar[r]^{1-T}& 
pM \ar[u]\ar[r]& 0.
\\
&0\ar[u]&0\ar[u]&0\ar[u]&
}$$
The vertical maps between the last two rows are natural inclusions, and the first row is the cokernels of these inclusion maps. The snake lemma implies that the first row is exact, which means that 
$$M^{1-T}/(pM)^{1-T}= (M/p)^{1-T}.$$ 
This is a finite $\Z/p^n$-module by the case $m=1$.
On the other hand, since $p^{m-1}\cdot pM=0$,  the induction hypothesis applied to the $W_nk$-module $pM$ gives $(pM)^{1-T}$ is a finite $\Z/p^n$-module. Now the vertical exact sequence on the left gives that $M^{1-T}$ is a finite $\Z/p^n$-module.
\epf

\bprop\label{Mss=M1-T otimes Wnk}
Let $k$ be a separably closed field of positive characteristic $p$. Then we have an identification of $W_nk$-modules
$$M_{\mathrm{ss}}\simeq M^{1-T}\otimes_{\Z/p^n} W_nk.$$
\eprop 

\bpf 
For the finite dimensional $k$-vector space $M/p$, \Cref{1-T surjective prop in appendix} tells us that
$$(M/p)_{\mathrm{ss}}\simeq (M/p)^{1-T}\otimes_{\F_p} k$$
In other words, there exists $m_1,\dots, m_d\in M$ ($d=\dim_{\F_p} M/p$), such that $m_1+pM,\dots,m_d+pM\in (M/p)^{1-T}$ generate $(M/p)_{\mathrm{ss}}$ as a $k$-vector space. 
Because of \Cref{1-T invariant is a finite module}(1), one can choose $m_1,...,m_d\in M^{1-T}$.
Since $M$ is a finite generated $W_nk$-module, $M_{\mathrm{ss}}$ as a submodule is also finitely generated over $W_nk$. 
Note moreover that $(M/p)_{\mathrm{ss}}=M_{\mathrm{ss}}/p$.
Apply Nakayama's lemma, $m_1,...,m_d\in M^{1-T}$ generate $M_{\mathrm{ss}}$ as an $W_nk$-module.
\epf

\end{appendices}

\medskip



\bigskip 

Bergische Universit\"at Wuppertal, 
Gau\ss strasse 20, D-42119 Wuppertal, Germany

\textit{Email address:} \url{renfei@uni-wuppertal.de
}
\end{document}